\documentclass[10pt,a4paper, reqno]{amsart}
\usepackage{amsfonts}
\usepackage[utf8]{inputenc}
\usepackage[foot]{amsaddr}
\usepackage{amsthm}
\usepackage{amsmath}
\usepackage{amscd}
\usepackage{t1enc}
\usepackage[mathscr]{eucal}
\usepackage{mathtools}
\usepackage{indentfirst}
\usepackage{graphicx}
\usepackage{graphics}
\usepackage{pict2e}
\usepackage{enumerate}
\usepackage{enumitem}
\usepackage{tcolorbox}
\usepackage{epic}
\numberwithin{equation}{section}
\usepackage[margin=2.9cm]{geometry}
\usepackage{epstopdf} 
\usepackage{float}
\usepackage{chngcntr}
\counterwithout{equation}{section}

\usepackage{dsfont}
\usepackage{nomencl}
\makenomenclature
\usepackage{hyperref}

\renewcommand{\nompreamble}{The next list describes several notations that will be later used within the body of the document}

\theoremstyle{plain}
\newtheorem{Th}{Theorem}[section]
\newtheorem{Lemma}[Th]{Lemma}
\newtheorem{Cor}[Th]{Corollary}
\newtheorem{Prop}[Th]{Proposition}

\theoremstyle{definition}
\newtheorem{Def}[Th]{Definition}

\newtheorem{?}[Th]{Problem}

\newcommand{\supp}{\textup{supp }}

\DeclareMathOperator*{\argmin}{arg\,min}

\usepackage{fancyhdr}
\renewcommand\headrulewidth{1.0pt}
\makeatletter
\def\headrule{{\if@fancyplain\let\headrulewidth\plainheadrulewidth\fi
\hrule\@height\headrulewidth\@width\headwidth
\vskip 2pt% 2pt between lines
\hrule\@height.5pt\@width\headwidth% lower line with .5pt line width
\vskip-\headrulewidth
\vskip-1.5pt}}
\makeatother

\providecommand{\keywords}[1]
{
  \small	
  \textbf{\textit{Keywords}} 
}

\begin{document}

\title{Analysis of simultaneous inpainting and geometric separation based on sparse decomposition}
\author{Van Tiep Do$^{*\ddagger}$, Ron Levie$^{*\dagger}$, Gitta Kutyniok$^{\dagger\S}$} 
       

%    author two information

\address{ $^{*}$ Department of Mathematics, Technische Universit{\"a}t Berlin, 10623 Berlin, Germany}

\address{ $^{\dagger}$ Department of Mathematics, Ludwig-Maximilians-Universit{\"a}t M{\"u}nchen, Munich, Germany} 

\address{ $^{\ddagger}$ Vietnam National University, 334 Nguyen Trai, Thanh Xuan, Hanoi} 

\address{ $^{\S}$ Department of Physics and Technology, University of Troms{\o}}

\email{tiepdv@math.tu-berlin.de, levie@math.lmu.de,  kutyniok@math.lmu.de}
\thanks{}


% \date{\today}

\dedicatory{}

\begin{abstract}

Natural images are often the superposition of various parts of different geometric characteristics. 
For instance, an image might be a mixture of cartoon and texture structures. In addition, images are often given with missing data. In this paper we develop a method for simultaneously decomposing an image to its two underlying parts and inpainting the missing data. Our separation--inpainting method is based on an $l_1$ minimization approach, using two dictionaries, each sparsifying one of the image parts but not the other.
We introduce a comprehensive convergence analysis of our method, in a general setting, utilizing the concepts of joint concentration, clustered sparsity, and cluster coherence. 
As the main application of our theory, we consider the problem of separating and inpainting an image to a cartoon and texture parts.
\end{abstract}

\maketitle

\section{Introduction}	

A digital image typically has two or more distinct constituents, e.g.,  it might contain a cartoon and texture components. A key question is whether we can separate the image into its two components. This task is of interest in many applications, such as compression and restoration \cite{5, 28, 57}. 
The separation problem is underdetermined and seems to be impossible to stably solve. However, if we have prior knowledge about the types of geometric components underlying the image, the separation task is possible,  as shown in previous papers \cite{4,23,27,32,36}. Some approaches commonly used in this context for image separation are variational methods, e.g., \cite{28, 37} and PDE-based separation methods \cite{37}. 

More recently, compressed sensing approaches showed that $l_1$ minimization can stably and precisely solve this problem both theoretically \cite{ 10, 26,31,35,40, 42} and empirically   \cite{4,27}. The core idea is to use multiple dictionaries, each sparsely representing one geometric part of the image. As opposed to thresholding approaches like \cite{22}, which are harder to analyze, the $l_1$ minimization approach comes with strong theoretical results.

Another classical problem in imaging science is to restore corrupted or missing parts of images, namely \emph{image inpainting}. Images are often damaged due to different factors including improper storage, chemical processing, or losses of image data during transmission.  The problem of image inpainting 
 is of interest in numerous applications, from the restoration of scratched photos and corrupted images, to the removal of selected objects. Image inpainting is an ill-posed problem, so  some prior information on the missing part is required here as well.
 One approach that has been effectively applied to this problem is total variational \cite{44,48,50,51}. Variational approaches work well on piecewise smooth images, but generally perform poorly on images that contain a superposition of cartoon and texture. On the other hand, local statistical analysis and prediction have been shown to perform well at inpainting the texture content \cite{46,49}. 
 \textcolor{black}{Another class of methods for image inpainting (and denoising) are transform-based methods. The key idea is to use  sparse representation systems such as wavelets, curvelets, shearlets, and framlets, to solve optimization problem in the transform domain. Among these methods, iterative thresholding algorithms have shown good performance for image inpainting \cite{45,59,60}.} %\textcolor{blue}{\Large You did not fix the problem that we discussed! Why is this another class of problems?! The problems that you discussed before also use dictionaries/frames! We agreed that this is not another class, so why did you delete my remarks without fixing the text?} 
 In our analysis, we focus on a  compressed sensing driven approach \textcolor{black}{via $l_1$ minimization}, based on the prior information that the components of the image are sparsely represented by two representation systems.  
 
\begin{figure}[H] 
  \includegraphics[width=0.335\textwidth]{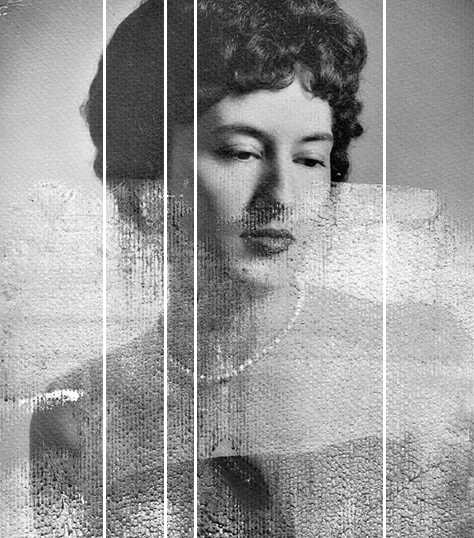}
  \caption{\label{FG1} 
   A corrupted photo, where the noise is a texture part together with missing stripes, and the clean image is the cartoon part.}
\end{figure}

 For an illustration of texture and cartoon parts of an image, we present in  Figure \ref{FG1} a corrupted photo with missing part and region covered by texture.

\subsection{Separation and inpainting through cluster coherence}

The problem of inpainting and separation of cartoon from texture was proposed in various papers, for instance \cite{24,33, 45}. The general approach is to find two dictionaries, each sparsely representing one image part, and not sparsely representing the other. For example, texture is sparsely represented by a Gabor system, and cartoon by curvelets/shearlets. Then, in the separation and inpainting algorithm, the image is decomposed to a sparse representation based on the  combined Gabor-curvelet system, using some compressed sensing approach.
However, in the above-mentioned papers, there is no theoretical analysis of the success of the proposed methods. \textcolor{black}{In \cite{24}  and \cite{45}, the separation and inpainting  problem is formulated as the minimization of the sum of the two $l_1$ norms of the two dictionary coefficients of the image. Each of these papers proposed an algorithm for the minimization, but neither proved that the $l_1$ minimization problem actually separates the two components. In our paper, we consider a similar setting to the above papers, and give a full approximation analysis of the method, with convergence guarantees. Our paper is the first that proves  that the $l_1$ minimization problem indeed inpaints and separates the two components of the image.}

Our analysis is based on the theoretical machinery of joint concentration and cluster coherence. These notions were used in the past for analyzing separation problems and inpainting problems, but never (to the best of \textcolor{black}{our} knowledge) for the simultaneous problem. Given two dictionaries, both joint concentration and cluster coherence 
 are definitions that quantify the ability of each dictionary to sparsely represent signals of one type, but not signals of another type.
Some papers use joint concentration to do inpainting \cite{19, 25}, some to do separation \cite{10,13}.  Proving theoretical results is typically easier when using joint concentration. However, checking joint concentration on dictionaries in practice is difficult. Thus the notion of cluster coherence was introduced in \cite{10}, which is strictly stronger than joint concentration, and makes checking the assumptions of the theory easier in practical examples.  Papers like \cite{10,13, 19,25} use the notions of cluster coherence either for separation or for inpainting, but not simultaneously.

In this paper, we modify and extend the definitions of joint concentration to accommodate a simultaneous analysis of separation and inpainting. We then further extend the notion of cluster coherence to fit the modified definitions of joint concentration. With the new general definitions, we consider the problem of  inpainting and separating texture from cartoon. We show that the universal shearlet systems \cite{25} and Gabor systems \cite{1}
 satisfy the cluster coherence requirement, thus proving the success of our inapinting--separation method.

\subsection{Our contribution}
We summarize our main contributions as follows. 
\begin{itemize}
    \item 
    We modify the joint concentration and cluster coherence definitions to accommodate the simultaneous separation-inpainting problem (Section \ref{S3}).
    \item Using the modified notions of joint concentration and cluster coherence  we prove that the $l_1$ sparse decomposition method successfully inpaints and separates the two parts of an image in a generic setting (Theorem \ref{TR1}).
    \item We use the general theory to inpaint and separate cartoon from texture in images. Here, the sparse decomposition method is based on a Gabor basis, sparsifying texture, and a universal shearlet dictionary, sparsifying cartoons (Section \ref{S4}). The two systems are shown to satisfy the cluster coherence property, thus proving the success of the inpainting-separation cartoon from texture method (Theorem \ref{TR10}).
\end{itemize}

\section{Simultaneous inpainting and separation approach} 

In this section we summarize our general theory of simultaneous inpainting and separation.

\subsection{The separation problem} 

The task is to extract the two components $\mathcal{C}$ and  $\mathcal{T}$ from the observed image $f$, where we assume that 
\begin{equation} \label{EQ1001}
f= \mathcal{C} + \mathcal{T}.
\end{equation} 
Here, only $f$ is given, and the components $\mathcal{C}$ and $\mathcal{T}$ are unknown to us. 

To solve this underdetermined problem, we assume that each component can be sparsely represented by some dictionary, that cannot sparsely represent the other component. In our theory, dictionaries are modeled as frames \cite{38}.
\begin{Def}
Let $I$ be a discrete index set. 
A sequence $\Phi = \{ \phi_i \}_{i \in I}$ in a separable Hilbert space $\mathcal{H}$ is called a frame for $\mathcal{H}$ if there exist constants $0 < A \leq B < \infty$  such that 
$$ A \| f \|_2^2 \leq \sum_{i \in I} | \langle f, \phi_i \rangle |^2 \leq B \|f \|_2^2, \quad  \forall f \in \mathcal{H},$$
where $A, B$ are called lower and upper frame bound. If $A$ and $B$ can be chosen to be equal, we call it $(A-)$ \emph{tight frame}. If A=B=1, $\{ \phi_i \}_{i \in I}$ is called a \emph{Parseval frame}.
\end{Def}

In our context, the Hilbert space $\mathcal{H}$ is the space of signals/images, and $\Phi$ is the dictionary.
By abusing notion, we also denote by $\Phi$ the \emph{synthesis operator}
$$ \Phi : l_2(I) \rightarrow \mathcal{H}, \quad \Phi(\{ a_i \}_{i \in I} ) = \sum_{i \in I }a_i \phi_i ,$$
which synthesizes an image $f\in\mathcal{H}$ from given coefficients $\{ a_i \}_{i \in I}$.
We denote by $\Phi^*$  the \emph{analysis operator}
$$ \Phi^* : \mathcal{H}  \rightarrow l_2(I), \quad \Phi^*(f ) = (\langle f, \phi_i \rangle)_{i \in I} .$$
The analysis operator is interpreted as the transform that computed the different dictionary coefficients of images $f$.

\subsection{The inpainting problem} 

To model the inpainting problem, we suppose that there is some missing data in the image $f$. Given the Hilbert space $\mathcal{H}$, we assume that $\mathcal{H} = \mathcal{H}_K \oplus \mathcal{H}_M$, where the subspaces $\mathcal{H}_K, \mathcal{H}_M$ denote the \emph{known part} and \emph{missing part} respectively. Let $P_{M}$ and $P_{K}$ denote the orthogonal projection of $\mathcal{H}$ upon these two subspaces, respectively.  In the inpainting problem, we assume that we are only given the image in the known part $P_K f$, and the goal is to reconstruct $f$. The inpainting problem is solved by considering a dictionary $\{\phi_i\}_{i\in I}$ that sparsely represents $f$ using a known subset of indices $\Lambda \subset I$, but does not sparsely represent $P_M f$ using the indices $\Lambda$. This idea is formalized in Definition \ref{De11}.
\subsection{The simultaneous separation--inpainting problem} 

In the simultaneous problem, the goal is to extract $\mathcal{C}$ and $\mathcal{T}$, satisfying (\ref{EQ1001}), given only the image in the known part $P_Kf$. In our approach, we consider two Parseval frames $\Phi_1$ and $\Phi_2$ in $\mathcal{H}$ that sparsely represent their respective component $\mathcal{C}$ and $\mathcal{T}$, but do not sparsely represent the other component.  This is formalized through Definitions \ref{Def4} and \ref{De10} for the joint concentration approach, and  \ref{Def4} and \ref{Def1} for the cluster coherence approach. We moreover suppose that the index subsets $\Lambda_1$ and $\Lambda_2$ of $\Phi_1$ and $\Phi_2$ respectively can sparsely represent $\mathcal{C}$ and $\mathcal{T}$, but cannot sparsely represent $P_M\mathcal{T}$ and $P_M\mathcal{C}$.  This is formalized through Definitions \ref{Def4}  and \ref{De11} for the joint concentration approach, and  \ref{Def4} and \ref{Def1} for the cluster coherence approach. 

We consider the following algorithm, for simultaneously inpaint and separate \textsc{(Inp-Sep)} geometric components, based on $l_1$ minimization.

\begin{tcolorbox}[colframe=white, colback=gray!10]
  \parbox{\textwidth}{
    \textbf{Algorithm}  \textsc{(Inp-Sep)}
    \vskip 1.5mm
    INPUT: corrupted signal $P_Kf \in \mathcal{H}_K,$ two Parseval frames $\{ \Phi_1 \}_{i \in I} $ and $ \{ \Phi_2 \}_{j \in J} $. 
    \vskip 1.5mm
    COMPUTE: $f^\star = (\mathcal{C}^\star, \mathcal{T}^\star)$ where 
   \begin{equation} (\mathcal{C}^\star, \mathcal{T}^\star ) = \argmin_{x_1, x_2} \| \Phi_1^* x_1 \|_1 + \| \Phi_2^* x_2 \|_1, \quad \textcolor{black}{s.t.} \quad  P_K(x_1 + x_2) = P_K(f). \label{EQ1000} \end{equation}
    \vskip 0.5mm
   OUTPUT: recovered components $\mathcal{C}^\star, \mathcal{T}^\star.$
      } 
\end{tcolorbox} 

In Section \ref{S3}, we provide a theoretical analysis for the success of algorithm (\ref{EQ1000}).

\section{General separation and inpainting theory} \label{S3}

In this section, we introduce a theory in which we can prove the success of Algorithm \textsc{(Inp-Sep)} for a general separation and inpainting problem.

\subsection{Joint concentration}

We propose a sufficient condition for the success of the \textsc{(Inp-Sep)} algorithm, based on the notion of joint concentration. Joint concentration was first introduced in \cite{10}. We present two slightly different notions of joint concentration, modified for our needs.

\begin{Def} \label{De10} Let $\Phi_1, \Phi_2$ be two Parseval frames. Given two sets of coefficients $\Lambda_1, \Lambda_2$, define the \emph{mixed joint concentration} $\kappa_1 = \kappa_1(\Lambda_1, \Lambda_2)$ by
\begin{equation}
\kappa_1= \kappa_1(\Lambda_1, \Lambda_2) = \sup_{x,y \in \mathcal{H}} \frac{\| \mathds{1}_{\Lambda_1} \Phi_1^* x\|_1 + \| \mathds{1}_{\Lambda_2} \Phi_2^* y\|_1  }{\| \Phi_1^* y\|_1 + \| \Phi_2^* x\|_1}.
\end{equation}
\end{Def}

Given an $l_2$ normalized signal $f$, it is common to interpret the $l_1$ norm of $f$ as a measure of spread (opposite of sparsity). Thus, the mixed joint concentration with respect to the coefficient sets $\Lambda_1$ and $\Lambda_2$ quantifies the extent to which signals can have most of their energy supported and well spread in $\Lambda_j$ while being sparsely concentrated in the other frame $\Phi_k$, for $j\neq k$.
Bounding $\kappa_1$ from above is one of the conditions that ensure that Algorithm \textsc{(Inp-Sep)} succeeds in the separation task. 
Next, we introduce another joint concentration notion which will be used to prove the success of the inpainting method. 
\begin{Def} \label{De11}
Let $\Phi_1, \Phi_2$ be two Parseval frames. Given two sets of coefficients $\Lambda_1, \Lambda_2$, define the \emph{joint concentration of the missing part} $\kappa_2 = \kappa_2(\Lambda_1, \Lambda_2)$ by
$$  \kappa_2= \kappa_2(\Lambda_1, \Lambda_2)= \sup_{\substack{x,y \in \mathcal{H};  \\ P_K x = P_K y }} \frac{\| \mathds{1}_{\Lambda_1} \Phi_1^* (x-y) \|_1 + \| \mathds{1}_{\Lambda_2} \Phi_2^* (x-y)\|_1  }{\| \Phi_1^*x\|_1 + \| \Phi_2^* y\|_1}.$$
\end{Def}
The joint concentration of the missing part quantifies the extent to which signals which coincide on the known part can have most of the energy of their difference supported and well spread in  $\Lambda_1,\Lambda_2$ while begin sparse in $\Phi_1$ and $\Phi_2$. The joint concentration 
$\kappa_2$ is thus used to encode the geometric relation between the missing part $\mathcal{H}_M$ and expansions in $\Phi_1$ and $\Phi_2$. 
Bounding $\kappa_2$ from above is another one of the conditions that ensure the success of Algorithm \textsc{(Inp-Sep)}.

\subsection{Recovery guarantee through joint concentration}
 The sufficient condition for recovery is based on finding sets of coefficients $\Lambda_1$ and $\Lambda_2$ for which the joint concentrations are small, and $\Lambda_1$ and $\Lambda_2$ capture most of the energy of the ground truth separated components $\mathcal{C}$ and $\mathcal{T}$.
For this, we recall the following definition from \cite{19}.
\begin{Def} \label{Def4}
Fix $\delta > 0$ . Given a Hilbert space $\mathcal{H}$ with a Parseval frame $\Phi, f \in \mathcal{H} $ is \emph{$\delta-$relatively sparse} in $\Phi$ with respect to $\Lambda$ if $\| \mathds{1}_{\Lambda^c} \Phi^* f \|_1 \leq \delta , $ where $\Lambda^c $ denotes $X \setminus \Lambda.$
\end{Def}
 We now prove our first separation result.
\begin{Prop} \label{PR20}
For $\delta_1, \delta_2 > 0,$ fix $\delta= \delta_1 + \delta_2,$ and suppose that $f \in \mathcal{H}$ can be decomposed as $f= \mathcal{C} +  \mathcal{T} $ so that each component $\mathcal{C}, \mathcal{T}$ is $\delta_1, \delta_2-$relatively sparse in $\Phi_1$ and $\Phi_2$ with respect to $\Lambda_1$ and $\Lambda_2$, respectively. Let $(\mathcal{C}^\star, \mathcal{T}^\star)$ solve \textsc{(Inp-Sep)}. If we have $\kappa_1 + \kappa_2 < \frac{1}{2}$, then  
\begin{equation} \label{CT1}
\| \mathcal{C}^\star - \mathcal{C} \|_2 + \| \mathcal{T}^\star - \mathcal{T} \|_2 \leq \frac{2\delta}{1 - 2 (\kappa_1 + \kappa_2)}.
\end{equation} 
\end{Prop}
\begin{proof}
First, we set
$$  \kappa= \kappa(\Lambda_1, \Lambda_2) := \sup_{P_Ku= P_Kv} \frac{\| \mathds{1}_{\Lambda_1} \Phi_1^* u\|_1 + \| \mathds{1}_{\Lambda_2} \Phi_2^* v\|_1  }{\| \Phi_1^* u\|_1 + \| \Phi_2^* v\|_1},$$
and $x= \mathcal{C}^\star - \mathcal{C}, \, y = \mathcal{T} - \mathcal{T}^\star$. Since $\Phi_1$ and $\Phi_2$ are two Parseval frames, thus
\begin{eqnarray*} 
\|  \mathcal{C}^\star - \mathcal{C} \|_2 + \| \mathcal{T}^\star - \mathcal{T} \|_2 &=& \| \Phi_1^*(\mathcal{C}^\star - \mathcal{C} ) \|_2 + \| \Phi_2^*(\mathcal{T}^\star - \mathcal{T} ) \|_2 \\
&\leq &  \| \Phi_1^*(\mathcal{C}^\star - \mathcal{C} ) \|_1 + \| \Phi_2^*(\mathcal{T}^\star - \mathcal{T} ) \|_1 \\
& = &  \| \Phi_1^*(x) \|_1 + \| \Phi_2^*(y) \|_1: = S.
\end{eqnarray*}
Now we invoke the relation $P_K(\mathcal{C}^\star + \mathcal{T}^\star ) = P_K(f) = P_K(\mathcal{C} + \mathcal{T})$ and hence $P_K(x) = P_K(y)$. By the definition of $\kappa,$ we have
\begin{eqnarray*}
S &=&\| \mathds{1}_{\Lambda_1} \Phi_1^* x \|_1 +\| \mathds{1}_{\Lambda_2} \Phi_2^* y \|_1 + \| \mathds{1}_{\Lambda_1^c} \Phi_1^* (\mathcal{C}^\star - \mathcal{C}) \|_1 +\| \mathds{1}_{\Lambda_2^c} \Phi_2^* (\mathcal{T}^\star - \mathcal{T})\|_1\\
 & \leq & \kappa S + \| \mathds{1}_{\Lambda_1^c} \Phi_1^* \mathcal{C}^\star \|_1 +\| \mathds{1}_{\Lambda_1^c} \Phi_1^* \mathcal{C} \|_1 + \| \mathds{1}_{\Lambda_2^c} \Phi_2^* \mathcal{T}^\star \|_1 +\| \mathds{1}_{\Lambda_2^c} \Phi_2^* \mathcal{T}\|_1 \\
 & \leq & \kappa S + \| \mathds{1}_{\Lambda_1^c} \Phi_1^* \mathcal{C}^\star \|_1 + \| \mathds{1}_{\Lambda_2^c} \Phi_2^* \mathcal{T}^\star \|_1 +\delta \\
  & = & \kappa S + \delta +  \| \Phi_1^* \mathcal{C}^\star \|_1 + \| \Phi_2^* \mathcal{T}^\star \|_1  - \| \mathds{1}_{\Lambda_1} \Phi_1^* \mathcal{C}^\star \|_1 - \| \mathds{1}_{\Lambda_2} \Phi_2^* \mathcal{T}^\star \|_1.
\end{eqnarray*}
We note that $(\mathcal{C}^\star, \mathcal{T}^\star) $ is a minimizer of \textsc{(Inp-Sep)}. Thus
\begin{eqnarray*}
\| \Phi_1^* \mathcal{C}^\star \|_1 + \| \Phi_2^* \mathcal{T}^\star \|_1  &\leq & \| \Phi_1^* \mathcal{C} \|_1 + \| \Phi_2^* \mathcal{T} \|_1 .
\end{eqnarray*}
Therefore,
\begin{eqnarray*}
S &\leq & \kappa S + \delta + \| \Phi_1^* \mathcal{C} \|_1 + \| \Phi_2^* \mathcal{T} \|_1 -\| \mathds{1}_{\Lambda_1} \Phi_1^* \mathcal{C}^\star \|_1 - \| \mathds{1}_{\Lambda_2} \Phi_2^* \mathcal{T}^\star \|_1  \\
&\leq & \kappa S + \delta + \| \Phi_1^* \mathcal{C} \|_1 + \| \Phi_2^* \mathcal{T} \|_1 + \|\mathds{1}_{\Lambda_1} \Phi_1^* x\|_1 + \|\mathds{1}_{\Lambda_2} \Phi_2^* y\|_1-\| \mathds{1}_{\Lambda_1} \Phi_1^* \mathcal{C} \|_1  \\
&& - \| \mathds{1}_{\Lambda_2} \Phi_2^* \mathcal{T} \|_1  \\
&\leq & \kappa S + 2 \delta + \kappa S = 2\kappa S + 2 \delta .
\end{eqnarray*} 
Thus,  $$ S \leq \frac{2 \delta}{1 - 2 \kappa}.$$
The last inequality comes from the fact that
\begin{eqnarray*}
    \frac{\| \mathds{1}_{\Lambda_1} \Phi_1^* x\|_1 + \| \mathds{1}_{\Lambda_2} \Phi_2^* y\|_1  }{\| \Phi_1^* x\|_1 + \| \Phi_2^* y\|_1} &\leq& \frac{\| \mathds{1}_{\Lambda_1} \Phi_1^* y\|_1 +
    \| \mathds{1}_{\Lambda_2} \Phi_2^* x\|_1  }{\| \Phi_1^* x\|_1 + \| \Phi_2^* y\|_1} + \\
    && \quad
  \frac{\| \mathds{1}_{\Lambda_1} \Phi_1^* (x-y) \|_1 
    + \| \mathds{1}_{\Lambda_2} \Phi_2^* (x-y)\|_1  }{\| \Phi_1^* x\|_1 + \| \Phi_2^* y\|_1}.
\end{eqnarray*} 
This leads to  $\kappa \leq \kappa_1 + \kappa_2$. Finally, we obtain (\ref{CT1})
$$ \| \mathcal{C}^\star - \mathcal{C} \|_2 + \| \mathcal{T}^\star - \mathcal{T} \|_2 \leq \frac{2\delta}{1 - 2 (\kappa_1 + \kappa_2)}. $$
\end{proof} 
 %\textcolor{black}{On the same line, a theoretical guarantee for data completion via geometric separation appeared in \cite{58} based on  notion of joint concentration, however, the construction goes against the philosophy of finding a frame that sparsely describes one part and not the other.  In our analysis, we  exploit the separate task of separation and inpainting corresponding to the role of $\kappa_1$ and $\kappa_2$. In addition, the theoretical analysis is not meaningful if we can not make it a concrete example to use in applications. We therefore try to use a new approach to transfer our theory from notions of joint concentration to notion of cluster coherence which we will define next.}
 %
 %
%
 %
 %\textcolor{blue}{\Large [Ron: Is this the paper of Emily that we discussed? Didn;t we realize that the requirement is not meaningful. This is very different from ``hard to check.'' Hard to check is not an excuse... Please summarize shortly why the requirement in their paper is meaningless.]} %
 
  \textcolor{black}{We note that \cite[Theorem 1]{58} also proposed a theoretical guarantee for the inpainting and separation problem based on joint concentration and relative sparsity. However, while the derived theorem is mathematically sound, we did not find any way to use it in practice. The reconstruction bound in \cite[Theorem 1]{58} appears artificially similar to our Proposition \ref{PR20}, but the  components $x$ and $y$ are analyzed via the frames that sparsify them in the denominator of the mixed joint concentration, so the denominator is small.  Hence, the formulation in \cite[Theorem 1]{58} goes against the philosophy of finding a frame that sparsely describes one part and not the other. In contrast, in our formulation, each $x$ and $y$ is analyzed via the frame that sparsifies the other component, and the denominator can be made large by choosing appropriate frames.}

 \subsection{Recovery guarantee through cluster coherence}

 Deriving bounds for joint concentrations can be difficult in practice. Our goal in this subsection is to replace the joint concentrations by easier to compute terms.
 For that, we recall the notion of cluster coherence. We then prove that joint concentrations can be bounded by cluster coherence terms. Deriving bounds to cluster coherence terms is generally easier than bounding joint concentrations. The following definition is taken from \cite{10}.
\begin{Def} \label{Def1}
Given two Parseval frames $\Phi_1 = (\Phi_{1i})_{i \in I}  $ and $\Phi_2 = (\Phi_{2j})_{j \in  J} $. Then the \emph{cluster coherence} $\mu_c(\Lambda, \Phi_1; \Phi_2) $ of $\Phi_1$ and $\Phi_2$ with respect to the index set $\Lambda \subset I $ is defined by
$$\mu_c(\Lambda, \Phi_1; \Phi_2)=  \max_{j\in J} \sum_{i \in \Lambda} | \langle \phi_{1i}, \phi_{2j} \rangle |$$
in case this maximum exists.
\end{Def}

We allow applying a projection $P_M$ on one or both of the frames in Definition \ref{Def1}. For example, similarly to \cite{19}, we use the notation $\mu_c(\Lambda,P_M \Phi_1; P_M \Phi_2)$ to denote 
$$\mu_c(\Lambda,P_M \Phi_1; P_M \Phi_2)=  \max_{j\in J} \sum_{i \in \Lambda} | \langle P_M \phi_{1i},P_M \phi_{2j} \rangle |.$$

Next, we present two lemmas which bound the joint concentrations by  corresponding cluster coherences. 

\begin{Lemma} \label{LM20}
We have
$$ \kappa_1(\Lambda_1, \Lambda_2) \leq \max \{ \mu_c(\Lambda_1, \Phi_1; \Phi_2), \mu_c(\Lambda_2, \Phi_2; \Phi_1)\} .$$
\end{Lemma}
\begin{proof} 
For $x, y \in \mathcal{H}$, we set $\alpha_1= \Phi_1^*y, \alpha_2= \Phi_2^* x $. Then invoking the fact that $\Phi_1$ and $\Phi_2$ are two Parseval frames, hence $x = \Phi_2 \Phi_2^* x = \Phi_2 \alpha_2$ and $y = \Phi_1 \Phi_1^* y = \Phi_1 \alpha_1$, we have
\begin{eqnarray*}
&&\| \mathds{1}_{\Lambda_1} \Phi_1^* x\|_1 + \| \mathds{1}_{\Lambda_2} \Phi_2^* y\|_1 \\
& =& \| \mathds{1}_{\Lambda_1} \Phi_1^* \Phi_2 \alpha_2 \|_1 + \| \mathds{1}_{\Lambda_2} \Phi_2^* \Phi_1 \alpha_1 \|_1 \\
& \leq &  \sum_{i \in \Lambda_1} \Big( \sum_j | \langle \Phi_{1i}, \Phi_{2j} \rangle | |\alpha_{2j}| \Big) +  \sum_{j \in \Lambda_2} \Big( \sum_i | \langle \Phi_{1i}, \Phi_{2j} \rangle | |\alpha_{1i}| \Big) \\
&  = &  \sum_{j} \Big( \sum_{i \in \Lambda_1} | \langle \Phi_{1i}, \Phi_{2j} \rangle | \Big) |\alpha_{2j}|  +  \sum_{i} \Big( \sum_{j \in \Lambda_2}| \langle \Phi_{1i}, \Phi_{2j} \rangle | \Big) |\alpha_{1i}| \\
& \leq & \mu_c(\Lambda_1, \Phi_1; \Phi_2) \| \alpha_2 \|_1 + \mu_c(\Lambda_2, \Phi_2; \Phi_1) \| \alpha_1 \|_1 \\
& \leq & \max \{ \mu_c(\Lambda_1, \Phi_1; \Phi_2), \mu_c(\Lambda_2, \Phi_2; \Phi_1)\} (\| \alpha_2 \|_1 + \| \alpha_1 \|_1 ) \\
& =  & \max \{ \mu_c(\Lambda_1, \Phi_1; \Phi_2), \mu_c(\Lambda_2, \Phi_2; \Phi_1)\} (\|\Phi_2^* x \|_1 + \| \Phi_1^* y \|_1 ).
\end{eqnarray*}
Thus, we obtain $\kappa_1(\Lambda_1, \Lambda_2) \leq \max \{ \mu_c(\Lambda_1, \Phi_1; \Phi_2), \mu_c(\Lambda_2, \Phi_2; \Phi_1)\}.$
\end{proof}
\begin{Lemma} \label{LM21}
We have
\begin{eqnarray*}
\kappa_2(\Lambda_1, \Lambda_2) & \leq &  \max \{ \mu_c(\Lambda_1, P_M\Phi_1; P_M \Phi_1) +  \mu_c(\Lambda_2, P_M \Phi_2; P_M \Phi_1), \\
\quad &&
 \mu_c(\Lambda_2, P_M\Phi_2; P_M \Phi_2) +  \mu_c(\Lambda_1, P_M \Phi_1; P_M \Phi_2) \} \\
& = & \max \{ \mu_c(\Lambda_1, P_M\Phi_1; \Phi_1) +  \mu_c(\Lambda_2, P_M \Phi_2; \Phi_1), \\
\quad &&
\mu_c(\Lambda_2, P_M\Phi_2; \Phi_2) +  \mu_c(\Lambda_1, P_M \Phi_1; \Phi_2) \}.
 \end{eqnarray*}
\end{Lemma}
\begin{proof}
For each $h\in \mathcal{H}_M$ and $x \in \mathcal{H}$, set $\alpha = \Phi_1^*(x+h)$ and $\beta= \Phi_2^* x$ then we have
\begin{eqnarray*}
h &=& x + h - x \\
&=& \Phi_1 \Phi_1^*(x+h)  -  \Phi_2 \Phi_2^*(x) \\
&=& \Phi_1 \alpha - \Phi_2 \beta.
\end{eqnarray*}
Note that $P_K$ is an orthogonal projection, and hence
$h = P_M^* P_M h = P_M^* P_M \Phi_1 \alpha  -  P_M^* P_M  \Phi_2 \beta.$ Therefore, we obtain
\begin{eqnarray*}
\| \mathds{1}_{\Lambda_1} \Phi_1^* h \|_1 & = & \| \mathds{1}_{\Lambda_1} \Phi_1^* P_M^* P_M \Phi_1 \alpha -  \mathds{1}_{\Lambda_1} \Phi_1^* P_M^* P_M  \Phi_2 \beta \|_1  \\
&=& \| \mathds{1}_{\Lambda_1} (P_M \Phi_1)^* ( P_M \Phi_1) \alpha   -  \mathds{1}_{\Lambda_1} (P_M\Phi_1)^* (P_M\Phi_2) \beta \|_1  \\
& \leq &  \| \mathds{1}_{\Lambda_1} (P_M \Phi_1)^* ( P_M \Phi_1) \alpha \|_1   +  \| \mathds{1}_{\Lambda_1} (P_M\Phi_1)^* (P_M\Phi_2) \beta\|_1  \\
& \leq & \sum_{i \in \Lambda_1} \Big( \sum_j | \langle P_M\Phi_{1i}, P_M \Phi_{1j} \rangle  | |\alpha_j|  \Big ) + \\
\quad && \sum_{i \in \Lambda_1} \Big( \sum_j | \langle P_M\Phi_{1i}, P_M \Phi_{2j} \rangle | |\beta_j|  \Big )  \\
& \leq & \sum_{j} \Big( \sum_{i \in \Lambda_1} | \langle P_M\Phi_{1i}, P_M \Phi_{1j} \rangle  | \Big )  |\alpha_j|  + \\
\quad && \sum_{j} \Big( \sum_{i \in \Lambda_1} | \langle P_M\Phi_{1i}, P_M \Phi_{2j} \rangle | \Big ) |\beta_j|   \\
& \leq & \mu_c(\Lambda_1, P_M\Phi_1; P_M \Phi_1) \| \alpha \|_1  +  \mu_c(\Lambda_1, P_M \Phi_1; P_M \Phi_2) \| \beta \|_1 \\
& = & \mu_c(\Lambda_1, P_M\Phi_1; P_M \Phi_1) \| \Phi_1^*(x+h) \|_1  +  \mu_c(\Lambda_1, P_M \Phi_1; P_M \Phi_2) \| \Phi_2^*(x) \|_1.
\end{eqnarray*}
Similarly, we have 
$$ \| \mathds{1}_{\Lambda_2} \Phi_2^* h \|_1 \leq \mu_c(\Lambda_2, P_M\Phi_2; P_M \Phi_1) \| \Phi_1^*(x+h) \|_1  +  \mu_c(\Lambda_2, P_M \Phi_2; P_M \Phi_2) \| \Phi_2^*(x) \|_1 .$$ 
This leads to 
\begin{eqnarray*}
\| \mathds{1}_{\Lambda_{1}} \Phi_1^* h \|_1 + \| \mathds{1}_{\Lambda_2} \Phi_2^* h \|_1 & \leq &  \max \{ \mu_c(\Lambda_1, P_M\Phi_1; P_M \Phi_1) +  \mu_c(\Lambda_2, P_M \Phi_2; P_M \Phi_1),  \\
\quad &&
 \mu_c(\Lambda_2, P_M\Phi_2; P_M \Phi_2) +  \mu_c(\Lambda_1, P_M \Phi_1; P_M \Phi_2) \}  \\
\quad && ( \| \Phi_1^*(x+h) \|_1  +   \| \Phi_2^*(x) \|_1).
\end{eqnarray*}

Finally,
\begin{eqnarray*}
\kappa_2(\Lambda_1, \Lambda_2) & \leq &  \max \{ \mu_c(\Lambda_1, P_M\Phi_1; P_M \Phi_1) +  \mu_c(\Lambda_2, P_M \Phi_2; P_M \Phi_1), \\
\quad &&
 \mu_c(\Lambda_2, P_M\Phi_2; P_M \Phi_2) +  \mu_c(\Lambda_1, P_M \Phi_1; P_M \Phi_2) \} \\
& = & \max \{ \mu_c(\Lambda_1, P_M\Phi_1; \Phi_1) +  \mu_c(\Lambda_2, P_M \Phi_2; \Phi_1), \\
\quad &&
\mu_c(\Lambda_2, P_M\Phi_2; \Phi_2) +  \mu_c(\Lambda_1, P_M \Phi_1; \Phi_2) \}.
\end{eqnarray*}
\end{proof}

We can now formulate a general guarantee for the success of Algorithm \textsc{(Inp-Sep)}, based on cluster coherence.
\begin{th} \label{TR1} Let 
$\Phi_1, \Phi_2$ be two Parseval frames for a Hilbert space $\mathcal{H}$. For $\delta_1, \delta_2 > 0,$ fix $\delta= \delta_1 + \delta_2,$ and suppose that $f \in \mathcal{H}$ can be decomposed as $f= \mathcal{C} +  \mathcal{T} $ so that each component $\mathcal{C}, \mathcal{T}$ is $\delta_1, \delta_2-$relatively sparse in $\Phi_1$ and $\Phi_2$ with respect to $\Lambda_1, \Lambda_2$ respectively. Let $(\mathcal{C}^\star, \mathcal{T}^\star)$ solve \textsc{(Inp-Sep)}. If we have $\mu_c< \frac{1}{2}$, then 
\begin{equation} \label{CT10}
\| \mathcal{C}^\star - \mathcal{C} \|_2 + \| \mathcal{T}^\star - \mathcal{T} \|_2 \leq \frac{2\delta}{1 - 2 \mu_c},
\end{equation} 
where $\mu_c = \max  \{ \mu_c(\Lambda_1, P_M\Phi_1; \Phi_1) +  \mu_c(\Lambda_2, P_M \Phi_ 2; \Phi_1), \mu_c(\Lambda_2, P_M\Phi_2; \Phi_2) +  \mu_c(\Lambda_1, P_M \Phi_1; \Phi_2)  \}  + \max  \{ \mu_c(\Lambda_1, \Phi_1; \Phi_2), \mu_c(\Lambda_2, \Phi_2; \Phi_1)  \} $.
\end{th}
\begin{proof}
The bound (\ref{CT10}) holds as a consequence of Lemmas \ref{LM20}, \ref{LM21} and Proposition \ref{PR20}.
\end{proof}

Let us interpret this estimate. The relative sparsity $\delta$ and the cluster coherence $\mu_c$ depend on the geometric sets of indices $\Lambda_1, \Lambda_2$, which are not used at all in Algorithm \textsc{(Inp-Sep)}. Hence, $\Lambda_1, \Lambda_2$ are analytic tools that we can choose arbitrarily, and are used only for deriving theoretical bounds.  If we choose very large $\Lambda_1, \Lambda_2$, we get very small relative sparsities $\delta_1, \delta_2$, but we might loose control over the cluster coherence $\mu_c$.  Therefore, choosing appropriate $\Lambda_1$ and $\Lambda_2$ is an important step when applying our theory in concrete examples. 

\section{Mathematical models of texture and cartoon} 

We now focus on the specific problem of inpainting and separating texture from cartoon. In this section we define our models of texture and cartoon parts, and of the missing part.

\subsection{Model of texture} 
 One of the earliest qualitative texture description that corresponds to human visual perception is: coarseness, contrast, directionality, line-likeness, regularity, roughness  \cite{41}. 
 From a mathematical modeling point of view, many definitions of texture were proposed in the past, for instance \cite{13, 30, 39, 43}.
Our model for texture is inspired by \cite{13}, and based on an expansion of Gabor frame elements. 
To motivate our definition of texture we offer the following discussion. A very restrictive definition of texture would be a periodical signal with a short period. The Fourier transform of a periodical signal is a signal supported on a delta train on a regular grid in the frequency domain. To relax the hard periodicity condition, suppose that we allow to perturb the locations of the points of the regular frequency grid. This would result in a signal supported on a sparse set of Fourier coefficients, which is exactly our definition of texture. Such a definition indeed produces images that look like a repeating pattern, but not in a strictly periodical fashion (see Figure \ref{FG2-3}-left for example). 

Before formally defining texture, we recall the \emph{Schwartz functions} or the \emph{rapidly decreasing functions}
\begin{equation} \label{EQ44}
     \mathbb{S}(\mathbb{R}^2):= \Big \{ f \in C^\infty(\mathbb{R}^2)  \mid   \forall K, N \in \mathbb{N}_0: \sup_{x \in \mathbb{R}^2} (1 + |x|^2 )^{-N/2} \sum_{|\alpha |\leq K} | D^\alpha f(x) | < \infty \Big \}. 
\end{equation}
We define the Fourier transform and inverse Fourier transform  for $f,F \in  \mathbb{S}(\mathbb{R}^2)$ by 
\textcolor{black}{$$ \hat{f}(\xi)=\mathcal{F}[f](\xi)= \int_{\mathbb{R}^2} f(x) e^{-2\pi i  x^T \xi} dx,$$}
$$ \check{F}(x) =\mathcal{F}^{-1}[F](x)= \int_{\mathbb{R}^2} F(\xi) e^{2\pi i \xi^T x} dx,$$
which can be extended to a  \textcolor{black}{well defined} Fourier transform and inverse Fourier transform for functions in $L^2(\mathbb{R}^2)$ (cf.\;\cite{56}). 
Using a window function $g:\mathbb{R}^2\rightarrow\mathbb{R}$ which localizes a texture patch in the spatial domain, we now introduce our model for texture. 
\begin{Def} \label{Df3}
Let $g \in L^2(\mathbb{R}^2)$ be a window with $\hat{g} \in C^{\infty}(\mathbb{R}^2),$ and frequency support $\supp \hat{g} \in [-1,1]^2$ satisfying the partition of unity condition
\begin{equation} \label{EQ139}
     \sum_{n \in \mathbb{Z}^2}|\hat{g}(\xi +n)|^2=1, \quad \xi \in \mathbb{R}^2. 
\end{equation}
For $s>0$, we define the $L^2-$ normalized scaled version of $g$ by 
$ g_s(x) = s \cdot g(sx).$ Let $I_T \subseteq \mathbb{Z}^2,$ be a subset of Fourier elements. A  texture is defined by
\begin{equation} \label{CT8}
\mathcal{T}_s(x) = \sum_{n \in I_T} d_{n} g_s(x) e^{2 \pi i x^T sn },
\end{equation}
where $(d_{n})_{n \in \mathbb{Z}^2}$ denotes  a sequence of complex numbers. 
\end{Def}

\textcolor{black}{A real valued texture image can be generated via Definition \ref{Df3} by choosing complex coefficients $(d_{n})_{n \in \mathbb{Z}^2}$ with the appropriate symmetry in the frequency domain ($\bar{d}_n=d_{-n}$).}

In Subsection \ref{The separation and inpainting theorem for cartoon and texture}   we add a restriction on the size of the index set $I_T$. For now, we just mention that $I_T$ is a small/sparse set in some sense. 
We also remark that by Definition \ref{Df3},  we have $\hat{g}_s(\xi) = s^{-1} \cdot \hat{g}(s^{-1} \xi)$,   $\supp \hat{g}_s \subseteq [-s, s]^2$ and the partition of unity condition (\ref{EQ139}) now reads
\begin{equation} \label{EQ140}
\sum_{n \in \mathbb{Z}^2} | \hat{g}_s(\xi + sn)|^2 = s^{-2}, \quad \xi \in \mathbb{R}^2. 
\end{equation}

\subsection{The local cartoon patch}  
In \cite{52, 53}, cartoon functions are defined as
\begin{equation}
\label{EQ:Cart}
     \mathcal{C} = f_0 + f_1 \cdot \mathds{1}_{B_{\tau}}, 
\end{equation}
where $f_0, f_1 \in L^2(\mathbb{R}^2) \cap C^\beta(\mathbb{R}^2), \beta \in (0, + \infty)$,  with compact support and $B_{\tau}$ denotes the interior of a closed, non-intersecting curve $\tau$ in $C^{\beta}(\mathbb{R}^2)$. 
In the separation algorithm of cartoon and texture, we analyze the input image using a \emph{ Gabor frame} and the so-called \emph{universal shearlet frame} (see Subsection \ref{S4}). Both of these frames analyze images on local patches. Thus, for the sake of simplicity, we also reduce our model of the cartoon part to a local cartoon model. This is done in two steps.

Cartoon images are locally close to piecewise constant functions with discontinuity along an edge curve, which is locally close to a line. We thus consider the \emph{windowed step function}
\begin{equation}
    w \mathcal{S}_0(x) = \begin{cases} w(x_1) & \quad x_2 \leq 0   \\ 0 & \quad x_2 >0,  \end{cases} \label{EQ22}
\end{equation} 
where $w \in C^{\infty}(\mathbb{R})$ is a \emph{weighted function} satisfying $w \not\equiv 0,  0 \leq w(u) \leq 1 $ and $\supp w \subset [-\rho, \rho].$
In our work, the cartoon part is analyzed using universal shearlets, which is a shear invariant system up to the cone adaptation. Since shearing is a way of changing orientation, our analysis also applies to the more general case where the discontinuity in $(\ref{EQ22})$ is along a general line in any orientation.
We note that it is possible to extend our results from the local cartoon model (\ref{EQ22}) to the global cartoon model (\ref{EQ:Cart}) by using a \emph{tubular neighborhood} argument as shown  in \cite{13}. 

The cartoon part $(\ref{EQ22})$ can be seen as a \emph{tempered distribution} acting on \emph{Schwartz functions} by
\begin{equation} \label{CT4}
\langle w \mathcal{S}_0 , f \rangle = \int_{-\rho}^{\rho} w(x_1) \int_{-\infty}^0 f(x_1, x_2) dx_2 dx_1, f \in \mathbb{S}(\mathbb{R}^2).
\end{equation} 
%\textcolor{blue}{The following lemma is used for computing the Fourier transform of $w \mathcal{S}_0$.
%\begin{Lemma}
%We have $$\int_{- \infty}^0 e^{-2 \pi i \omega x} d x= \frac{1}{2} \Big [ \delta(\omega) + \frac{i}{\pi  \omega} \Big ]. $$ \label{LM30}
%\end{Lemma}
%By Lemma \ref{LM30}, we can compute the Fourier transform of $w \mathcal{S}_0$ by 
%\begin{align}
% \langle \widehat{ w \mathcal{S}_0}, f  \rangle = \langle w \mathcal{S}_0, \hat{f} \rangle & =  \int_\rho^\rho w(x_1) \int_{-\infty}^0 \Big( \int_{\mathbb{R}^2} f(\xi) e^{-2 \pi i (\xi_1 x_1+ \xi_2 x_2) } d \xi \Big )dx_2 dx_1  \nonumber \\
%  &= \int_{\mathbb{R}^2} \hat{w}(\xi_1) \Big ( \int_{- \infty}^0 e^{-2 \pi i \xi_2 x_2} d x_2\Big )f(\xi) d\xi  \nonumber \\
%  &= \int_{\mathbb{R}^2} \hat{w}(\xi_1)   \Big ( \delta(\xi_2) + \frac{i}{\pi  \xi_2} \Big ) f(\xi) d\xi, \quad f \in \mathcal{L}(\mathbb{R}^2 ) . \label{CT5}
%\end{align}
%By $(\ref{CT5})$, we obtain $\widehat{w \mathcal{S}_0}(\xi) =  \hat{w} (\xi_1)   \Big ( \delta(\xi_2) + \frac{i}{\pi  \xi_2} \Big ). $ }
\textcolor{black}{By the fact that
$\widehat{\mathds{1}_{(-\infty,0]}}= \frac{1}{2} \Big [ \delta(\omega) + \frac{i}{\pi  \omega} \Big ]$, 
we have $\widehat{w \mathcal{S}_0}(\xi) =  \hat{w} (\xi_1)   \Big ( \delta(\xi_2) + \frac{i}{\pi  \xi_2} \Big )$.
}

\begin{figure}[H] \label{FG2-3}
  \centering
  \includegraphics[width=0.35\textwidth]{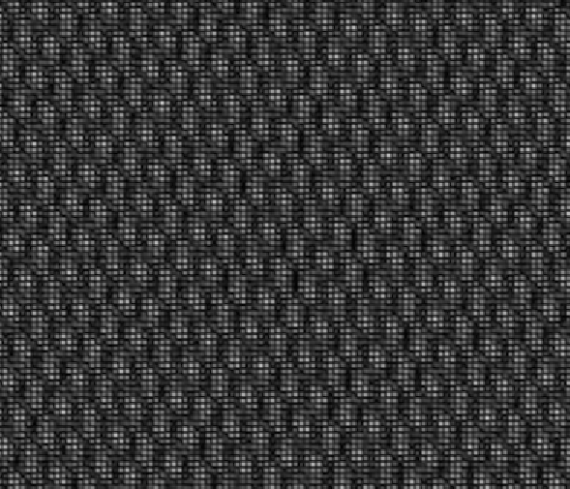}
  \quad
  \includegraphics[width=0.46\textwidth]{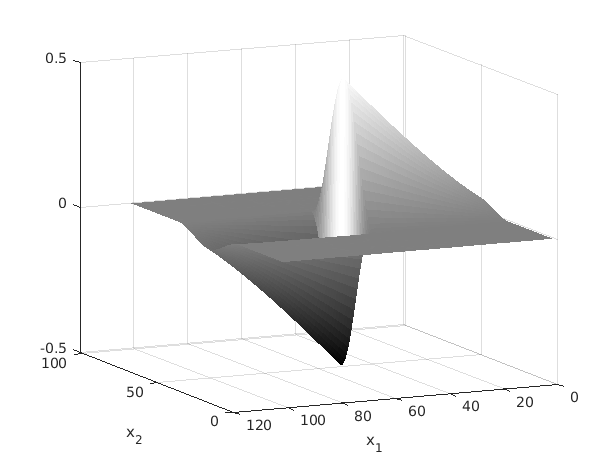}
  \caption{ \textcolor{black}{Left: Real valued} texture sample produced by randomly choosing a sparse set of Fourier coefficients with random values \textcolor{black}{and an appropriate symmetry in the frequency domain}, with the rest of the Fourier coefficients set to zero.  \textcolor{black}{Right:  A local cartoon patch.} }
\end{figure}

Now, we modify the local cartoon part to get an image in $L_2(\mathbb{R}^2)$.  We note that we are interested in the local behaviour of $w\mathcal{S}_0$ about $(x_1,x_2)=0$. Thus, the ``DC part'' of $\widehat{w\mathcal{S}_0}$ is not of interest to us, as it models some global asymptotic behaviour of the patch in $\mathbb{R}^2$. We thus filter our the band $|\xi_2|<r$ from  $\widehat{w\mathcal{S}_0}(\xi_1,\xi_2)$ for some arbitrarily small $r>0$, to obtain
\begin{equation} \label{EQ153}
    \widehat{w \mathcal{S}} = \mathds{1}_{\{ \lvert \xi_2 \rvert \geq r \} } \widehat{w \mathcal{S}_0}  .
\end{equation} 
It is easy to see that 

\begin{eqnarray}
\widehat{w\mathcal{S}} &=& \mathds{1}_{\{ |\xi_2| \geq r \}}(\xi_2)  \widehat{\mathds{1}_{(-\infty,0]}} (\xi_2)\hat{w}(\xi_1)  \nonumber \\
&= & \frac{i}{2 \pi \xi_2} \mathds{1}_{\{ |\xi_2| \geq r\}}(\xi_2) \hat{w}(\xi_1). \label{Revise_1}
\end{eqnarray} 
The main change incurred by this filtering of $w\mathcal{S}_0$ is a ``translation along the $y$ axis $w\mathcal{S}_0 \mapsto w\mathcal{S}_0-{\rm DC}$'', which does not affect the analysis via shearlet frames.
This filtering changes the behaviour of $w\mathcal{S}_0$ about $(x_1,x_2)=0$ negligibly, since $r$ is small, and retains the quality of $w \mathcal{S}_0$ being approximately piecewise constant locally with a line discontinuity (see Figure \ref{FG2-3}-right).

\textcolor{black}{To give more intuition about how to interpret the model as local, we note that in the multi-scale approach of Subsection \ref{sub_01}, the resulting models of cartoon are indeed spatially localized. This can be seen through Lemma \ref{LM100}, where the shearlet coefficients of the cartoon models are shown to decay fast in the spatial parameter. }

%{we take the coefficient space approach to localization.  In our case, and in the multi-scale approach,  there is actually a fast decay of the coefficients in the spatial direction of $\xi_2$ for each scale $j$.  In this sense, our cartoon model is local - it can be described to high accuracy with shearlet elements in some local domain.}

%\textcolor{black}{\Large You did not include a few words about why this is local.}
Henceforth,  we call the local cartoon part $w \mathcal{S}$ simply a cartoon part.

 \subsection{The missing part}
In our analysis, the shape of the missing region is chosen to be
   \begin{equation}
   \label{EQ:missing}
    \mathcal{M}_h= \{ x= (x_1, x_2) \in \mathbb{R}^2 \mid |x_1| \leq h \}, \quad \text{for } h >0 ,    
   \end{equation}
 and the orthogonal projection associated with this missing part is $P_M = \mathds{1}_{\mathcal{M}_h}.$
 
 %\begin{figure}[H] 
 %\centerline{\psfig{file=Images/anh3.png,width=5cm}}
%\vspace*{8pt}
%\caption{ \label{FG4} The missing part (gray) at scale $j$.} 
%\end{figure} 

Note that this models a local and axis aligned missing part, corresponding to the local model of texture (\ref{CT8}) and the local and axis aligned cartoon model (\ref{EQ153}). When combining and re-orienting the local patches, we can obtain many missing stripes at various orientations.
Moreover, in practice, the missing part can be any domain contained in $\mathcal{M}_h$ of (\ref{EQ:missing}). 

One practical example where the missing part is of the form (\ref{EQ:missing}) is seismic data, where the image is commonly incomplete due to missing or faulty sensor or land development causing white strips \cite{16,17}. 

\section{Sparsifying systems for texture and cartoon \label{S4}}

In this section, we choose frames that sparsely represent the texture and  cartoon parts.
To represent texture, it is clear from Definition \ref{Df3} that a  Gabor frame is a natural choice.
For the cartoon part,
among popular sparsifying representation systems such as wavelets, curvelets, and shearlets,  it was shown in \cite{15} that shearlets outperform not only wavelets, but also most other directional sparsifying systems  for inpainting larger gap size \cite{19, 25}.   

We hence choose the following sparse representation systems.
\begin{itemize}
\item  Gabor frame - a tight frame with time-frequency balanced elements
\item Universal shearlet frame -  a directional tight frame. 
\end{itemize}

\subsection{Gabor frame}

Gabor frames are defined as follows.

\begin{Def} \label{Def2}
Denote for  $x, y \in \mathbb{R}^2$ and  $g \in  L^2(\mathbb{R}^2)$
$$ g(x,y)(t)=g(t-x)e^{2 \pi i t^Ty}. $$
For $a, b> 0$ we call the collection of functions $\{g(ma,nb) \mid m,n \in \mathbb{Z}^2 \}$ a \emph{Gabor system}. Such a system is called a \emph{Gabor frame} if it forms a frame.
\end{Def}

For our Gabor system,  we use a window  $\tilde{g}$ that satisfies the assumptions of the window $g$ of texture (Definition \ref{Df3}). We note that  $\tilde{g}$ and $g$  may be different in general. However, the analysis for $\tilde{g}\neq g$ is almost identical to the analysis in case $g=\tilde{g}$. Without loss of generality, we assume henceforth that both Definitions \ref{Df3} and \ref{Def2} are based on the same window $g$, satisfying the assumptions of Definition \ref{Df3}.

For each scaling factor $s>0$, we consider the  Gabor tight frame $\mathbf{G_s}= \{ (g_s)_\lambda(x) \}_\lambda$. This frame is represented in the frequency domain as
\begin{equation} \label{EQ0}
   (\hat{g})_\lambda(\xi) = \hat{g}_s(\xi - sn) e^{2 \pi i  \xi^T \frac{m}{2s}},  
\end{equation}
where $\lambda=(m,n)$ is the spatial and frequency position and the parameter $s$ denotes the band-size.  This system constitutes a tight frame for $L^2(\mathbb{R} ^2)$, see \cite{11} for more details.

\subsection{The universal-scaling shearlet frame}

Shearlets were first introduced in \cite{6} as a representation system extending the wavelet framework. They are a directional representation system with similar optimal approximation properties as curvelets, but with the advantage of allowing for a unified treatment of the continuum and digital domains.  Extending the shearlet system,  $\alpha$-shearlets, first introduced in  \cite{20},  can be regarded as a parametrized family ranging from wavelets to shearlets. 
In \cite{52, 53} it was shown that $\alpha-$shearlets provide optimally sparse approximations for cartoon images defined as piecewise $C^{1/\alpha}-$functions, separated by a $C^{1/\alpha}$ singularity curve.
A further extension is universal shearlets \cite{25}, which is motivated by \cite{29}, with the aim at constructing a type of $\alpha-$shearlets that forms a Parseval frame. 
In our setting we choose universal shearlets as the sparsifying system of the cartoon part, since they form a Parseval frame, a necessary condition in our theory. 

Let $\phi$ be a function in $\mathbb{S}(\mathbb{R})$ satisfying $0 \leq \hat{\phi}(u) \leq 1$ for $u\in\mathbb{R}$, $\hat{\phi}(u) =1$ for $u \in [-1/16, 1/16]$ and $\supp \hat{\phi} \subset [-1/8, 1/8]$. Define the \emph{low pass function} $\Phi(\xi)$ and the \emph{corona scaling functions} for $j \in \mathbb{N} $ and $\xi= (\xi_1,\xi_2 ) \in \mathbb{R}^2,$
$$ \hat{\Phi}(\xi):= \hat{\phi}(\xi_1) \hat{\phi}(\xi_2),$$
\begin{equation} \label{CT35}
    W(\xi):= \sqrt{\hat{\Phi}^2(2^{-2}\xi)-\hat{\Phi}^2(\xi)}, \quad W_j(\xi):= W(2^{-2j} \xi).
\end{equation} 
\begin{figure}[H]
\centering
\includegraphics[width=150pt, height=150pt]{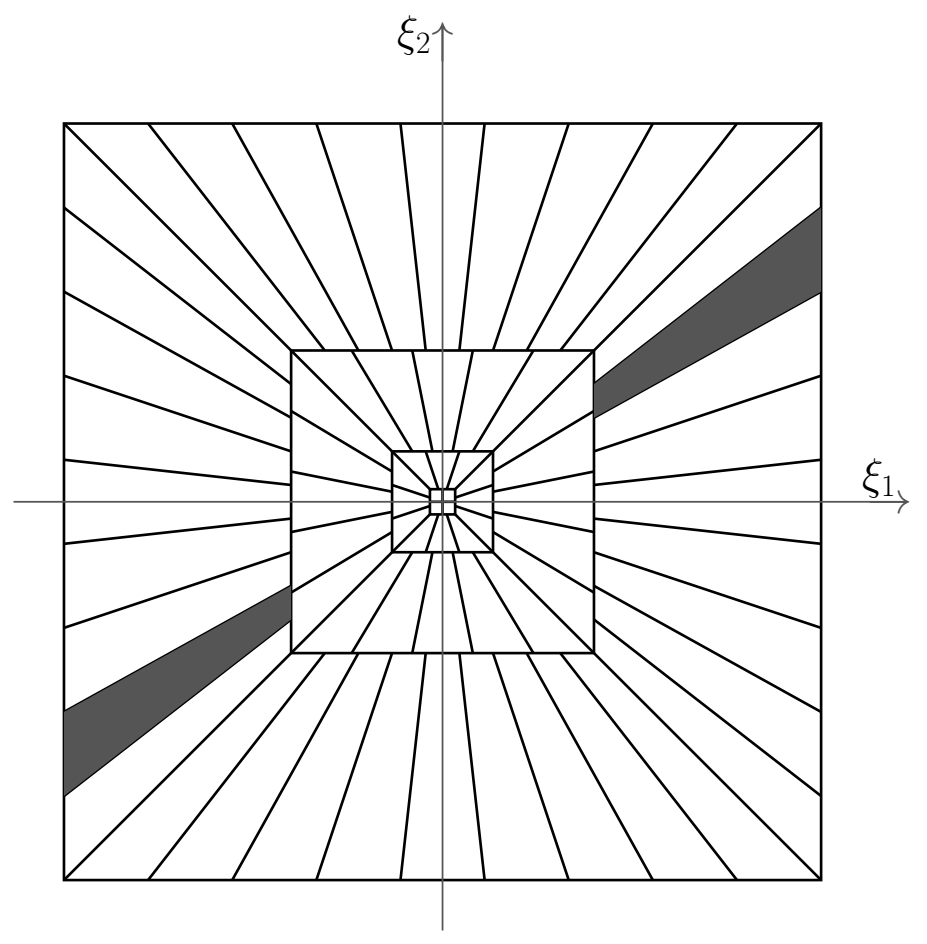}
\caption{\label{FG5} Frequency tiling of a cone-adapted shearlet. } 
\end{figure} 
It is easy to see that we have the partition of unity property
\begin{equation} \label{CT3}
 \hat{\Phi}^2(\xi) + \sum_{j \geq 0 } W_j^2(\xi) = 1, \quad \xi \in \mathbb{R}^2.
\end{equation}

Next, we use a bump-like function $v \in C^\infty (\mathbb{R})$ to produce the directional scaling feature of the system. Suppose $\supp(v) \subset [-1,1] $ and $| v(u-1)|^2 + |v(u)|^2 + |v(u+1)|^2 =1$ for $ u \in [-1,1]$. Define the \emph{horizontal frequency cone} and the \emph{vertical frequency cone}
\begin{equation} \label{EQ30}
\mathcal{C}_{(\rm{h})} := \Big \{ (\xi_1, \xi_2) \in \mathbb{R}^2 \mid \Big | \phantom{.}  \frac{\xi_2}{\xi_1} \Big | \leq 1 \Big \} 
\end{equation}
and
\begin{equation} \label{EQ31}
     \mathcal{C}_{\rm{(v)}} := \Big \{ (\xi_2, \xi_1) \in \mathbb{R}^2 \mid  \Big  | \frac{\xi_1}{\xi_2} \Big | \leq 1 \Big \}.
\end{equation}
Define  the cone functions $V_{(\rm{h})}, V_{(\rm{v})}$ by
 \begin{equation} \label{EQ32}
 V_{(\rm{h})}(\xi): = v\Big ( \frac{\xi_2}{\xi_1} \Big ), \quad V_{\rm{(v)}}(\xi): = v\Big ( \frac{\xi_1}{\xi_2} \Big ). 
 \end{equation}
The \emph{shearing} and \emph{scaling matrix} are defined by
\begin{equation} \label{EQ33}
A_{\alpha, (\rm{h})} := 
\begin{bmatrix}
    2^2   & 0 \\
      0  & 2^\alpha \\
   \end{bmatrix}, \quad S_{(\rm{h})} := \begin{bmatrix}
    1  & 1 \\
      0  & 1 \\
   \end{bmatrix} ,
\end{equation}

\begin{equation} \label{EQ34}
 A_{\alpha, \rm{(v)}} := 
\begin{bmatrix}
    2^\alpha & 0 \\
      0  & 2^2 \\
   \end{bmatrix}, \quad S_{\rm{(v)}} := \begin{bmatrix}
    1  & 0 \\
      1  & 1 \\
   \end{bmatrix},
\end{equation}
where $\alpha \in (-\infty, 2)$ is the \emph{scaling parameter}. 
The following two definitions are taken from \cite{25}.

\begin{Def} 
Let $\Phi, W, v$ be defined as above.  For $\alpha \in (-\infty, 2),  k \in \mathbb{Z}^2 $, we define 
\begin{enumerate}[label*=\arabic*.]
\item  \emph{Coarse scaling functions}: $\psi_{-1,k}(x):= \Phi(x-k),  k \in \mathbb{Z}^2, x \in \mathbb{R}^2.$
\item \emph{Interior shearlets}: let $ j \in \mathbb{N}_0, l \in \mathbb{Z}$ such that $|l| < 2^{(2-\alpha)j}, k \in \mathbb{Z}^2$ and $\rm{(\iota)} \in \{ \rm{(h)},\rm{(v)}\}$. Then we define $\psi_{j,l,k}^{\alpha, (\iota)}(x)$ by its Fourier transform
 $$ \hat{\psi}_{j,l,k}^{\alpha, (\iota)}(\xi): = 2^{-(2+\alpha)j/2}W_j(\xi) V_{(\iota)} \Big ( \xi^T A_{\alpha, (\iota)}^{-j} S^{-l}_{(\iota)} \Big ) e^{-2\pi i \xi^T A_{\alpha, (\iota)}^{-j} S^{-l}_{(\iota)}k},  \xi \in \mathbb{R}^2.$$
 \item  \emph{Boundary shearlets}: let $ j \in \mathbb{N},$ and $ l = \pm \lceil 2^{(2-\alpha)j} \rceil $ we define 
 $$ \hat{\psi}_{j,l,k}^{\alpha, \rm{(b)}}(\xi):  = \begin{cases} 2^{-(2+\alpha)j/2 -1/2}W_j( \xi) V_{(\rm{h})} \Big ( \xi^T A_{\alpha, (\rm{h})}^{-j} S^{-l}_{(\rm{h})} \Big ) e^{-\pi i \xi^T A_{\alpha, (\rm{h})}^{-j}  S^{-l}_{(\rm{h})}k},&\xi \in \mathcal{C}_{(\rm{h})} \\ 2^{-(2+\alpha)j/2 -1/2}W_j(\xi) V_{\rm{(v)}} \Big ( \xi^T A_{\alpha, \rm{(v)}}^{-j} S^{-l}_{\rm{(v)}} \Big ) e^{-\pi i \xi^T A_{\alpha, \rm{(v)}}^{-j}  S^{-l}_{\rm{(v)}}k},&\xi \in \mathcal{C}_{\rm{(v)}} \end{cases} $$ and in the case $j=0, l = \pm 1,$ we define
  $$ \hat{\psi}_{0,l,k}^{\alpha, \rm{(b)}}(\xi):  = \begin{cases} W(\xi) V_{(\rm{h})} \Big ( \xi^T  S^{-l}_{(\rm{h})} \Big ) e^{-2\pi i \xi^T k}, &  \xi \in \mathcal{C}_{(\rm{h})} \\ W(\xi) V_{\rm{(v)}} \Big ( \xi^T  S^{-l}_{\rm{(v)}} \Big ) e^{-2\pi i \xi^T k}, &  \xi \in \mathcal{C}_{\rm{(v)}}.\end{cases} $$
\end{enumerate}
 \end{Def}

\begin{Def}  \label{Df4}
Let $(\alpha_j)_{j \in \mathbb{N}_0} \subset \mathbb{R}$ be a \emph{scaling sequence}, i.e, 
\begin{equation} \label{EQ37}
    \alpha_j \in Z_j:= \Big \{ \frac{m}{j} \mid  m \in \mathbb{Z}, m \leq 2j-1 \Big \} .
\end{equation}
We define the associated \emph{universal-scaling shearlet system}, or shorter, \emph{universal shearlet system}, by
$$\boldsymbol{\Psi}=\textup{SH}(\phi, v, (\alpha_j)_j):= \textup{SH}_{\textup{Low}} (\phi) \cup \textup{SH}_{\textup{Int}}(\phi, v, (\alpha_j)_j) \cup \textup{SH}_{\textup{Bound}}(\phi, v, (\alpha_j)_j), $$
where
\begin{eqnarray*}
&& \textup{SH}_{\textup{Low}}(\phi)  := \{ \psi_{-1,k}(x) \mid  k \in \mathbb{Z}^2 \}  \\
&& \textup{SH}_{\textup{Int}} (\phi, v, (\alpha_j)_j):= \Big \{ \psi^{\alpha_j,  \rm{(\iota)}}_{j,l,k}(x) \mid  j \in \mathbb{N}, l \in \mathbb{Z}, |l| < 2^{(2-\alpha_j)j}, k \in \mathbb{Z}^2, \rm{\iota} \in \{ \rm{h},\rm{v}\} \Big \}  \\
&& \textup{SH}_{\textup{Bound}} (\phi, v, (\alpha_j)_j):= \Big \{ \psi^{\alpha_j, \rm{(b)}}_{j,l,k}(x) \mid  j \in \mathbb{N}, l \in \mathbb{Z},  |l| =\pm  2^{(2-\alpha_j)j}, k \in \mathbb{Z}^2\Big \}.
\end{eqnarray*}

\end{Def}
In \cite[Theorem 2.4]{25} it was shown that $\{ \boldsymbol{\Psi} \}_\eta$,   $\eta=(j, l, k; \alpha_j, \rm{(\iota)}) $ constitutes a Parseval frame for $L^2(\mathbb{R}^2)$, \textcolor{black}{and all shearlets (coarse scaling functions, interior shearlets, and boundary shearlets) are smooth in the frequency domain.}  
The universal shearlet system is \textcolor{black}{restricted to} discrete $\alpha_j \in Z_j$ in order for the set of sheared cone functions to exactly cover the horizontal and vertical frequency cones, \textcolor{black}{and for the boundary shearlet to be smooth in the frequency domain}. 
For a real $\alpha \in (-\infty, 2)$, we can choose $(\alpha_j)_j$ that provide the best possible approximation of $\alpha$, i.e, 
 \begin{equation} \label{EQ120}
   \alpha_j:= \argmin_{\tilde{\alpha}_j \in Z_j} |\tilde{\alpha}_j-\alpha|, \quad j \geq 1.   
 \end{equation}
This implies that $2^{\alpha_j j } = \Theta(2^{\alpha j })$ and $ \lim \alpha_j = \alpha$ as $j \rightarrow \infty$. Indeed, we can show that 
\begin{equation} \label{EQ5}
\alpha_j = \frac{\lfloor j \alpha + 0.5 \rfloor}{j}.
\end{equation}

For our model, we consider universal shearlets with an arbitrary scaling sequence $(\alpha_j)_{j \in \mathbb{N}_0} $ chosen as in (\ref{EQ5}) with $\alpha \in (0, 2)$.

\section{Inpainting and separation of cartoon and texture} 

In this section we formulate our main result on inpainting and separation of cartoon and texture.

\subsection{Multi-scale separation and inpainting \label{sub_01}}

In our analysis, the separation problem is analyzed in a multi-scale setting. What we show is that the inpainting and separation method is more accurate for smaller scale components of the image, as long as the size of the missing part gets smaller in scale.  \textcolor{black}{Here, we adopt the multi-scale asymptotic analysis approach which was introduced  in \cite{10} for the problem of image separation. This approach was then adopted in \cite{19} for image inpainting. The premise is that in the separation and inpainting problem  there is no direct way to do asymptotic analysis, as the problem is fixed. One simple way to do asymptotic analysis is letting the missing band become smaller, and computing the rate of convergence of the reconstructed parts to the ground true parts. However, such a setting is somewhat trivial. The idea in \cite{19} is to modify the model by breaking the cartoon part (and also the texture part in our case) to different scale components. By considering this sequence of models, we can formulate the accuracy of the problem as a function of the scale, thus deriving an asymptotic analysis. This approach was also adopted in \cite{25,55}. }

%\textcolor{black}{\Large Here you should write the philosophy of this analysis, and say that this asymptotic analysis approached was introduced in [*25*]}

%\textcolor{magenta}{The philosophy of the analysis should go later. Since I understand in Martin paper that, we just introduce the general problem of Inpainting with a certain model of missing part. But instead of considering this general problem we just apply it for each of sub-image to study assymptotic behaviour. In other words, we restrict our model of missing part to each sub-image to study assymptotic behaviour of the error.  } I think I wrote it in the right place.
%\textcolor{magenta}{After stating this we will focus on our model for each sub-image. But after (6.3) you stated the general problem. I think the philosophy should be around (6.5).}

\textcolor{black}{Since texture is a stationary phenomenon, it mainly has middle and high frequency components.  Thus, in our approach,  the  low  frequency  part  of  the  signal  is  directly  assigned  to  the cartoon part.  Namely, for the low frequency band, we only use the shearlet frame in the $l_1$ minimization problem, and treat the resulting signal as a cartoon part.  This is essentially the problem which is solved in \cite{25}, where the only difference there is that instead of cartoon the paper studies curve singularities.  Note that our analysis extends to the one component analysis.} 

%\textcolor{magenta}{Fo more detail, we can specify here Proposition 8.1 and Proposition 9.3 guarantee for the sucess of this restricted problem.}
Consider the window function $W$ defined in (\ref{CT35}). We construct a family of \emph{frequency filters} $F_j$ with the Fourier representation 
\begin{equation} \label{CT25}
   \hat{F}_j(\xi):=W_j(\xi) = W(\xi/2^{2j}), \quad \forall j \geq 1, \xi \in \mathbb{R}^2, 
\end{equation}
and in the case $j=0$ 
$$ \hat{F}_0(\xi):= \Phi(\xi), \quad \xi \in \mathbb{R}^2.$$ Notice that $W_j$ is compactly supported in the corona
\begin{equation} \label{EQ11}
 \mathcal{A}_j := [-2^{2j-1}, 2^{2j-1}]^2 \setminus [-2^{2j-4}, 2^{2j-4}]^2 , \quad  \forall j \geq 1.
\end{equation}

We now consider the texture and cartoon parts at different scales by filtering them with $F_j$.
For the patch size of texture, we typically use a fixed size $s$ for all scales $j$. However, we also allow $s$ to depend on the scale $j$ to make the setting more flexible. We hence consider a variable $s_j$, and associated with $s_j$ the domain $\mathcal{A}_{s, j}, \forall j \geq 1 $, defined by
\begin{equation} \label{EQ12}
 \mathcal{A}_{s, j} =  \{ \xi \in \mathbb{R}^2: s_j \xi \in \mathcal{A}_j \}.
 \end{equation}
% \textcolor{black}{We construct our problem as follows.  The goal is to recover the image $f=w \mathcal{S}+ \mathcal{T}_s$ from its known part $P_K f$, where the components $w\mathcal{S}$ and $\mathcal{T}_s$ are unknown to us.}
 %
 %\textcolor{magenta}{I think we should not state the general problem here, for this part we just focus on each sub-image to study assymptotic behaviour of the error. Or the philosophy should be here to restrict our model to each sub-image.}
 %
  %\textcolor{black}{NOW EXPLAIN WHY THIS MULTI-SCALE APPROACH MAKES SENSE, SINCE THIS IS NOT DIRECTLY THE FULL INPAINTING PROBLEM. Explain that we study this way the asymptotic behavior of the separation problem.} \textcolor{black}{Then write that one usecase of our construction is: Suppose I have my known image $(1-P_M)f$. Now, suppose that we want to inpaint and separate it all the way to some scale $J$. Now, take all $h_{j}=h_J$, namey constant. The theorey works, since for $j<J$ $h_J=o(2^{-\alpha_j j})$.} \textcolor{blue}{\Large Does the analysis now work if you replace the order of filtering and projecting? If so, this is a use case.}
  %
  %\textcolor{magenta}{I am not sure in this case. This changes the theory a lot. This means the constraint now is: $F_j*(P_Kf_j)= F*(P_K(x+y))$. This depends on the filter. Does this concludes $P_Kf_j = P_K(x+y$  ?. 
  %\begin{eqnarray*}
  %F_j*(P_Kf_j) &=& F*(P_K(x+y)) \\
  %F_j*F_j*(P_Kf_j) &=& F_j*F_*(P_K(x+y))  \\
  %W_j^2.\widehat{P_Kf_j} &=&  W_j^2. %\widehat{(P_K(x+y)} \\
  %\widehat{P_Kf_j} &=& \widehat{(P_K(x+y)} 
  %\end{eqnarray*} }
  
  We split the components \textcolor{black}{$w\mathcal{S}$} and $\mathcal{T}_s$ into pieces of different scales  $\textcolor{black}{w\mathcal{S}_j = F_j * w\mathcal{S}}$ and $\mathcal{T}_{s,j} = F_j * \mathcal{T}_{s}$, where $\mathcal{T}_{s}$ is the texture of patch size $s=s_j$.
 The total signal at scale $j$ is defined to be
   $$ f_j=\textcolor{black}{w\mathcal{S}_j} + \mathcal{T}_{s,j}.$$
 As mentioned above, we typically pick a constant $s_j=s$. In the constant $s$ case,
 the pieces $f_j$ can be filtered directly from $f$, namely,
$f_j = F_j * f$. We moreover have by (\ref{CT3})
\begin{equation} \label{EQ97}
   f = \sum_jF_j * f_j.
\end{equation} 
In the general case, the Fourier transform $\hat{f}_{j}$ of $f_{j}$ is supported in the annulus with inner radius $2^{2j-4}$ and outer radius $2^{2j-1}$. 
 Now,  at each scale $j$ we consider the simultaneous inpainting and separation problem  of extracting $\textcolor{black}{w\mathcal{S}_j}$ and $\mathcal{T}_{s,j}$
 from
\begin{equation} \label{EQ98}
  (1-P_j)f_{j} = (1-P_j)\textcolor{black}{w\mathcal{S}_j} + (1-P_j)\mathcal{T}_{s,j}.   
\end{equation}

As in \cite{25}, we assume that at scale $j$ the size of the missing part $h_j$ depends on $j$. 
Denote by  $P_j$ the characteristic function of the missing part at scale $j$, i.e,
\begin{equation} \label{DN2}
P_j = \mathds{1}_{\{ |x_1| \leq h_j \}}.
\end{equation}
%
%\textcolor{green}{Now, the part about reconstructing the whole image from its components should be deleted.}
Since generally $s=s_j$ depends on $j$, we use different Gabor frames to solve $(\ref{EQ98})$ for different scales. We denote by $\mathbf{G_j}$ the  Gabor tight frame associated with $s=s_j$, i.e,
 \begin{equation} \label{DN1} \mathbf{G_j}:=\mathbf{G_{s_j}}=\{(g_{s_j})_{(m,n)}(x) \}_{(m,n)\in \mathbb{Z}^2}.
 \end{equation}
 
 \textcolor{black}{Under the above construction, in Theorem \ref{TR10} we show that algorithm \textsc{(Inp-Sep)}, applied to each scale model (\ref{EQ98}), is asymptotically accurate as $j\rightarrow\infty$.}

In the analysis of the next subsections, we need the following representation of $P_j$ in the frequency domain. 
\begin{Lemma} \label{LM51}
We have 
$\hat{P}_j(\xi)= 2h_j {\rm sinc}(2h_j \xi_1) \delta(\xi_2).$
\end{Lemma}  
\begin{proof}
By definition, we have
\begin{eqnarray*}
    \hat{P}_j(\xi)&=& \int_{\mathbb{R}^2}P_j(x) e^{-2\pi i x^T \xi} dx \\
    &=& \Big (\int_{-h_j}^{h_j} e^{-2\pi i x_1 \xi_1 } dx_1 \Big )\Big (\int_\mathbb{R} e^{-2\pi i x_2\xi_2} dx_2 \Big )\\
    &=&2h_j {\rm sinc}(2h_j \xi_1) \delta(\xi_2).
\end{eqnarray*} 
\end{proof}

\subsection{Balancing the texture and cartoon parts}
It is important to avoid the trivial case where one of the parts, cartoon or texture, is much larger than the other. This would lead to the trivial separation outcome in which the whole image is taken as the estimation of the larger part, and the smaller part is estimated as zero.
We thus suppose that the filtered components $\textcolor{black}{w\mathcal{S}_j}$ and $\mathcal{T}_{s,j}$ have comparable magnitudes at each scale. 

Consider the sub-band components
\begin{equation} \label{DN4}
    w\mathcal{S}_j = F_j * w \mathcal{S},
\end{equation}
and  
\begin{equation} \label{DN5}
   \mathcal{T}_{s,j} = F_j *  \mathcal{T}_s, 
\end{equation}
where $w \mathcal{S}$ is defined in (\ref{EQ153}) and texture is in Definition \ref{Df3} .
The following lemma formulates the \emph{energy balancing condition}.

\begin{Lemma} $\label{LM1}$ 
Consider the cartoon patch and texture with sub-bands defined in (\ref{DN4}), (\ref{DN5}), respectively. We have
$$  \| w \mathcal{S}_j \|_2^2 \sim 2^{-2j}, \quad j \rightarrow \infty, $$ and
$$\| \mathcal{T}_{s,j} \|_2^2 \sim \sum_{n \in \mathbb{Z}^2 \cap \mathcal{A}_{s,j} } |d_{n}|^2. $$
Therefore, energy balance is achieved for
\begin{equation} \label{EQ1}
c_0 2^{-2j} \leq \sum_{n \in I_T \cap \mathcal{A}_{s,j} } |d_{n}|^2 \leq c_12^{-2j}, \quad  c_0, c_1 >0.
\end{equation}
\end{Lemma}

\begin{proof}
We present the following sketch of the proof.
For $W$ sufficient nice we have
\begin{eqnarray*}
\| w \mathcal{S}_j \|_2^2 &=&  \int_{\substack{ \xi \in \mathbb{R}^2, \\  |\xi_2| \geq r}}  |\hat{w}(\xi_1)|^2 \cdot W^2(2^{-2j }\xi)  \cdot \frac{1}{4}\Big \lvert \delta(\xi_2)+\frac{i}{\pi \xi_2} \Big \rvert^2  d\xi \\
& \sim & \int_{\substack{ \xi \in \mathcal{A}_j, \\  |\xi_2| \geq r}} \frac{1}{\xi_2^2} \cdot  |\hat{w}(\xi_1)|^2   d\xi \\
& \sim & \int_{\substack{ \xi \in \mathcal{A}_j,  \\  |\xi_2| \geq 2^{2j-5}}}  \frac{1}{\xi_2^2}|\hat{w}(\xi_1)|^2  \textcolor{black}{d\xi} + \int_{\substack{ \xi \in \mathcal{A}_j,  \\ 2^{2j-5} \geq|\xi_2| \geq r}}  \frac{1}{\xi_2^2}|\hat{w}(\xi_1)|^2  \textcolor{black}{d\xi}.
\end{eqnarray*}
Note that by (\ref{EQ11}), $\mathcal{A}_j=  [-2^{2j-1}, 2^{2j-1}]^2 \setminus [-2^{2j-4}, 2^{2j-4}]^2.$ Thus, $\{ \xi \in \mathcal{A}_j, 2^{2j-5} \geq |\xi_2 | \geq r \} \subset \{ \xi \in \mathcal{A}_j, |\xi_1| \geq 2^{2j-4} \} $.  We now use the rapid decay of $\hat{w}$
$$\sup_{|\xi_1| \geq 2^{2j-5} }|\hat{w}(\xi_1)| \leq C_N\langle |2^{2j}|\rangle^{-N}, \quad \forall N \in \mathbb{N}.$$
This leads to 
\begin{equation}
    \int_{\substack{ \xi \in \mathcal{A}_j,  \\ 2^{2j-5} \geq|\xi_2| \geq r}}  \frac{1}{\xi_2^2}|\hat{w}(\xi_1)|^2  d\xi_2 \leq  C_N \cdot  \frac{ 2^{2j}}{r^2} \cdot 2^{2j} \langle | 2^{2j} | \rangle^{-N}, \quad \forall N \in \mathbb{N}. \label{F1}
\end{equation} 
\textcolor{black}{ Since $\xi \in \mathcal{A}_j=  [-2^{2j-1}, 2^{2j-1}]^2 \setminus [-2^{2j-4}, 2^{2j-4}]^2$, we have
   $$
    C2^{-2j}=\int_{2^{2j-4}}^{2^{2j-1}} \frac{1}{\xi_2^2} d\xi_2\cdot \int_{-2^{2j-1}}^{2^{2j-1}} |\hat{w}(\xi_1)|^2 d\xi_1\leq \int_{\substack{ \xi \in \mathcal{A}_j,  \\  |\xi_2| \geq 2^{2j-5}}}  \frac{1}{\xi_2^2}|\hat{w}(\xi_1)|^2  d\xi. $$
 Thus,  }
\begin{equation}
    \int_{\substack{ \xi \in \mathcal{A}_j,  \\  |\xi_2| \geq 2^{2j-5}}}  \frac{1}{\xi_2^2}|\hat{w}(\xi_1)|^2  d\xi_2 \sim C\cdot 2^{-2j}. \label{F2}
\end{equation}
Combining (\ref{F1}) and (\ref{F2}), we finally obtain 
\begin{equation} \label{} \label{D1}
    \| w \mathcal{S}_j \|_2^2 \sim 2^{-2j}.
\end{equation}

On the other hand, we have $$ \mathcal{T}_{s} = \sum_{n \in I_{T}} d_{n} g_{s_j}(x) e^{2 \pi i \langle s_jn, x \rangle}.$$
This leads to 
\begin{equation}
    \hat{\mathcal{T}}_{s} = \sum_{n \in I_T} d_n \hat{g}_{s_j}(\xi -s_jn). \label{I0}
\end{equation} 
Now, using the change of variable $\omega = \frac{\xi}{s_j}$ we get
\begin{eqnarray*}
\| \mathcal{T}_{s,j} \|_2^2 &= & \| \hat{\mathcal{T}}_{s,j} \|_2^2=\sum_{n, \tilde{n} \in I_T} \int d_n \overline{d_{\tilde{n}}} W^2(\xi/2^{2j}) \hat{g}_{s_j}(\xi -s_jn)\hat{g}_{s_j}(\xi -s_j\tilde{n})d\xi  \\
& = & s_j^2 \cdot \sum_{\substack{|n- \tilde{n}| \leq 1\\ n, \tilde{n} \in  I_T  } }  \int d_n \overline{d_{\tilde{n}}} W^2(s_j \omega/2^{2j}) \hat{g}_{s_j}(s_j(\omega -n))\hat{g}_{s_j}(s_j(\omega -\tilde{n}))d\omega  \\
& = &  \sum_{\substack{|n- \tilde{n}| \leq 1\\ n, \tilde{n} \in  I_T  } } 
 \int d_n \overline{d_{\tilde{n}}} W^2(s_j \omega/2^{2j}) \hat{g}(\omega -n)\hat{g}(\omega -\tilde{n})d\omega . \\
\end{eqnarray*}

Note that for each $n$, there exist a finite number of $\tilde{n}$'s independent of $j$ satisfying $|n- \tilde{n}| \leq 1,$ and since $W$ is sufficiently nice,
$$ \int  W^2(s_j \omega/2^{2j}) \hat{g}(\omega -n)\hat{g}(\omega -\tilde{n})d\omega \sim \int_{\mathcal{A}_{s,j}}   \hat{g}(\omega -n)\hat{g}(\omega -\tilde{n})d\omega.$$
Thus, we obtain
\begin{eqnarray*}
\| \mathcal{T}_{s,j} \|_2^2  & \sim & \sum_{\substack{|n- \tilde{n}| \leq 1\\ n \in  I_T \cap \mathcal{A}_{s,j}  } } \int_{\mathcal{A}_{s,j}} d_n \overline{d_{\tilde{n}}}  \hat{g}(\omega -n)\overline{\hat{g}(\omega -\tilde{n})}d\omega \\
& \sim &  \sum_{ n \in   I_T \cap \mathcal{A}_{s,j}}  \int_{\mathcal{A}_{s,j}} |d_n|^2   |\hat{g}(\omega -n) |^2 d\omega.
\end{eqnarray*}
Up to a small set of coefficients,  we assume that the support of $\hat{g}$ is always entirely contained in $\mathcal{A}_{s,j}$. Hence,
\begin{eqnarray}
\| \mathcal{T}_{s,j} \|_2^2  & \sim & \sum_{ n \in  I_T \cap \mathcal{A}_{s,j}  }  \int_{\mathbb{R}^2} |d_n|^2   |\hat{g}(\omega -n) |^2 d\omega \nonumber \\
& = & \sum_{ n \in I_T \cap  \mathcal{A}_{s,j}  } |d_n|^2   \int_{\mathbb{R}^2}  |\hat{g}(\omega ) |^2 d\omega  \sim   \sum_{ n \in I_T \cap \mathcal{A}_{s,j}   }  |d_n|^2. \label{D2}
\end{eqnarray}
Combining (\ref{D1}) and (\ref{D2}), we finish the proof.
\end{proof}

For later analysis, we also need an energy estimate of the cartoon patch in the missing region. We use the following useful notation in the whole paper
\begin{equation}
 \langle |x| \rangle = (1 + |x|^2)^{1/2}. \label{CT2}
\end{equation} 
\begin{Lemma} \label{LM52}
For $h_j = o(2^{-\alpha_j j}), \alpha_j \in (0,2) $, and $\omega (0) \neq 0$, there exists a constant $C_1 >0 $ such that  
\begin{equation}
     \| P_j w \mathcal{S}_j  \|_2^2\geq C_1\cdot h_j^22^{-2j} \label{EQ60}
\end{equation}
\end{Lemma}
\begin{proof}

By Plancherel’s theorem and Lemma \ref{LM51}, we have 
\begin{eqnarray}
 \| P_j w \mathcal{S}_j \|_2^2 
 & =& \int_{\substack{\xi \in \mathbb{R}^2 \\ |\xi_2| \geq r }}  \Big | \int_{\mathbb{R}^2}  \hat{P}_j(t) \widehat{w \mathcal{S}_j}(\xi- t) dt \Big |^2  d\xi  \nonumber  \\
& = &  \int_{\substack{\xi \in \mathbb{R}^2 \\ |\xi_2| \geq r }}  \Big | \int_{\mathbb{R}}2h_j{\rm sinc}(2h_j t_1) \hat{w}(\xi_1-t_1) \frac{1}{(\pi\xi_2)}W_j(\xi_1-t_1, \xi_2) dt_1 \Big |^2 d\xi \nonumber    \\
& = & \frac{4h^2_j}{\pi^2} \int_{\substack{\xi \in \mathbb{R}^2 \\ |\xi_2| \geq r }} \frac{1}{\xi_2^2} \Big | \int_{\mathbb{R}}{\rm sinc}(2h_j t_1) \hat{w}(\xi_1-t_1) W_j(\xi_1-t_1, \xi_2) dt_1 \Big |^2  d\xi. \nonumber \\  \label{E1}
\end{eqnarray}

Next, 
recalling the definition of $W_j(\xi)$ in (\ref{CT35}), we have
\begin{equation}
    W_j(\xi) = 1, \quad \forall \xi \textup{ such that }  \xi \in  [-2^{2j-2}, 2^{2j-2}]^2 \setminus [-2^{2j-3}, 2^{2j-3}]^2. \label{E2} 
\end{equation} 
For $\xi$ in the corona-shape $\{ \xi \in \mathbb{R}^2 \mid \xi \in [-2^{2j-2}, 2^{2j-2}]^2 \setminus [-2^{2j-3}, 2^{2j-3}]^2 \}$, if $\xi_1$ is small then $\xi_2$ must be large. More accurately, for $\xi$ with $|\xi_1| \leq 1 $, we are guaranteed to have $\xi \in [-2^{2j-2}, 2^{2j-2}]^2 \setminus [-2^{2j-3}, 2^{2j-3}]^2$ in case $\xi_2 \in (2^{2j-3}, 2^{2j-2}). $ 
Combining this observation with (\ref{E1}), the following holds 
\begin{equation}
    \| P_j \omega \mathcal{S}_j \|_2^2 
\geq  C h^2_j \int_{2^{2j-3}}^{2^{2j-2}} \frac{1}{\xi_2^2} \int_{-1}^1  \Big |   \int_{\mathbb{R}} {\rm sinc}(2h_j t_1) \hat{\omega}(\xi_1-t_1) W_j(\xi_1-t_1, \xi_2)dt_1 \Big |^2 d\xi_1 d\xi_2. \label{E10}
\end{equation} 

Next, we study the term $I_j(\xi) := \Big | \int_{\mathbb{R}} {\rm sinc}(2h_j t_1) \hat{\omega}(\xi_1-t_1) W_j(\xi_1-t_1, \xi_2)dt_1 \Big | $.
By the triangle inequality, we have
\begin{eqnarray}
 I_j(\xi) & \geq  & \Big | \int_{|t_1| \leq 2^{\alpha_j j-4}} {\rm sinc}(2h_j t_1) \hat{\omega}(\xi_1-t_1) W_j(\xi_1-t_1, \xi_2)dt_1 \Big |  -  \nonumber \\
 && \Big | \int_{|t_1| \geq 2^{\alpha_j j-4}} {\rm sinc}(2h_j t_1) \hat{\omega}(\xi_1-t_1) W_j(\xi_1-t_1, \xi_2)dt_1 \Big |.\label{E9}
\end{eqnarray}
We now use the rapid decay of $\hat{w}$ to obtain 
\begin{equation}
    \sup_{|t_1| \geq 2^{\alpha_j j} }|\hat{w}(\xi_1-t_1)| \leq C_N\langle |2^{\alpha_j j}|\rangle^{-N}, \quad  \textcolor{black}{\forall \xi_1 \in [-1, 1],} \; \forall N \in \mathbb{N}. \label{E4}
\end{equation}
This leads to 
\begin{equation}
    \Big | \int_{\textcolor{black}{|t_1| \geq 2^{\alpha_j j-4} }} {\rm sinc}(2h_j t_1) \hat{\omega}(\xi_1-t_1) W_j(\xi_1-t_1, \xi_2)dt_1 \Big | \leq  C_N^\prime \langle |2^{\alpha_j j}|\rangle^{-N}.\label{E8}
\end{equation} 
For the other term, we observe that for fine scale $j$, we have $(\xi_1-t_1, \xi_2)  \in [-2^{2j-2}, 2^{2j-2}]^2 \setminus [-2^{2j-3}, 2^{2j-3}]^2 $ for $\xi_1 \in [-1, 1], |t_1| \leq 2^{\alpha_j j -4}, \alpha_j \in (0,2)$ and $\xi_2 \in (2^{2j-3}, 2^{2j - 2})$ . Combining this with (\ref{E2}), we derive $W_j(\xi_1 - t_1, \xi_2) = 1$. Thus, we obtain
\begin{eqnarray*}
       \int_{|t_1| \leq 2^{\alpha_j j-4}} {\rm sinc}(2h_j t_1) \hat{\omega}(\xi_1-t_1) W_j(\xi_1-t_1, \xi_2)dt_1 &=&  \int_{|t_1| \leq 2^{\alpha_j j-4}} {\rm sinc}(2h_j t_1) \\ 
       && \qquad \cdot \hat{\omega}(\xi_1-t_1) dt_1. 
\end{eqnarray*}

We now prove that there exists a constant $C^\prime > 0$  such that for sufficiently large $j$, we have
\begin{equation}
    \Big | \int_{|t_1| \leq 2^{\alpha_j j-4}} {\rm sinc}(2h_j t_1) \hat{\omega}(\xi_1-t_1) dt_1 \Big | \geq C^\prime, \quad \forall \xi_1 \in [-1,1]. \label{E7}
\end{equation}
By triangle inequality, we obtain
\begin{eqnarray}
 && \Big | \int_{|t_1| \leq 2^{\alpha_j j-4}} {\rm sinc}(2h_j t_1) \hat{\omega}(\xi_1-t_1) dt_1 \Big |  \geq  \nonumber \\
&& \Big | \int_{|t_1| \leq 2^{\alpha_j j -4}} \hat{w}(\xi_1-t_1) dt_1 \Big | -  
  \Big |\int_{|t_1| \leq 2^{\alpha_j j-4}}\Big ( {\rm sinc}(2h_j t_1) -1 \Big )  \hat{w}(\xi_1 -t_1) dt_1\Big |. \nonumber \\  \label{E3}
\end{eqnarray}
For the first term, the following holds for $\forall \xi_1 \in [-1,1]$ 
\begin{eqnarray}
\Big | \int_{|t_1| \leq 2^{\alpha_j j -4}} \hat{w}(\xi_1-t_1) dt_1 \Big |& \geq & \Big |\int_\mathbb{R}\hat{w}(t_1) dt_1 \Big |- \Big | \int_{|t_1| \geq 2^{\alpha_j j -4}} \hat{w}(\xi_1-t_1) dt_1 \Big | \nonumber \\
 & \stackrel{\mathclap{\normalfont{ \textup{by} \; (\ref{E4})}} } \approx &  \quad |w(0)| = \text{const} > 0.  \label{E6}
\end{eqnarray}
For the second term of the right hand side of (\ref{E3}), by the assumption $h_j2^{\alpha_j j } \xrightarrow{j \rightarrow +\infty}{0},$ we have
\begin{eqnarray}
&& \Big |\int_{|t_1| \leq 2^{\alpha_j j-4}} ( {\rm sinc}(2h_j t_1) -1  )  \hat{w}(\xi_1 -t_1) dt_1 \Big | \nonumber \\
 &\leq& \sup_{|t_1| \leq 2^{\alpha_j j -4}}|{\rm sinc}(2h_j t_1)-1| \int_\mathbb{R}|\hat{w}(t_1| dt_1 \nonumber \\
&\leq & C \cdot  \sup_{|t_1| \leq 2^{\alpha_j j -4}}|{\rm sinc}(2h_j t_1)-1|  \xrightarrow{j \rightarrow + \infty}{0} . \label{E5}
\end{eqnarray} 
Thus, by (\ref{E3}), (\ref{E5}) and (\ref{E6}), we conclude the proof of (\ref{E7}).

Now, by (\ref{E9}), (\ref{E8}), and (\ref{E7}), we have
$I_j(\xi) \geq C^\prime = \textup{const.}$ Combining this with (\ref{E10}), we finally obtain
\begin{eqnarray*}
  \| P_j \omega \mathcal{S}_j \|_2^2 
&\geq&  C h^2_j \Big (\int_{2^{2j-3}}^{2^{2j-\frac{5}{2}}} \frac{1}{\xi_2^2} d\xi_2  \Big ) \Big ( \int_{-1}^1 C^{\prime 2} d\xi_1 \Big ) \nonumber \\
& = & C_1 \cdot h^2_j 2^{-2j}, \label{E11}
\end{eqnarray*}
which concludes the proof.
\end{proof}

\subsection{Sparsity assumptions on texture} 
\label{The separation and inpainting theorem for cartoon and texture}  

Next, we restrict the definition of texture by bounding the size of the index set $I_T$ of Definition \ref{Df3}.
Define the domain
\begin{equation}
 M_{j, l , \rm{(\iota)}}: = \Big [ s_j^{-1} \cdot S_{(\iota)}^{-l} \cdot ( \supp \hat{\psi}_{j, 0,0}^{\alpha_{j}, \rm{(\iota)}} + B_1(0,1) )  \Big ] \cap \mathbb{Z}^2, \label{CT50}
\end{equation}
where  $A+B:=\{ a+b \mid a\in A, b \in B \}$ is the Minkowski sum, and multiplication of a set of vectors by a matrix is done elementwise. We assume that at every scale $j$,  the number of non-zero Gabor elements with the same position, generating $\mathcal{T}_{s,j}$, is  not too large.
More accurately, 
denote $|A| = \# \{ n \in \mathbb{Z}^2 \cap A \}$, and define the neighborhood set $I_T^{\pm} $  of the index set $I_T$ by 
\begin{equation} \label{CT26}
    I_T^{\pm} = \{n^\prime \in \mathbb{Z}^2 \mid \exists n \in I_T: | n^\prime -n| \leq 1  \}.
\end{equation}
What we assume is that at every scale $j$,  for all $ l \in \mathbb{Z}$ satisfying $ |l| \leq 2^{(2-\alpha_j)j} $ we have
\begin{equation}
|I_T \cap M_{j, l , \rm{(\iota)}}| \leq 2^{(2-\alpha_j-\epsilon)j/2 }=o(2^{(2-\alpha_j)j/2} ),  \quad  \forall \rm{(\iota)} \in \{ \rm{(h)}, \rm{(v)}, \rm{(b)} \}, \label{EQ2} 
\end{equation}
and
\begin{equation}
    | I_T^{\pm} \cap \mathcal{A}_{s,j} | \leq \frac{2^{\alpha_j j}}{s_j}   \quad \textup{as} \, j \rightarrow \infty . \label{EQ20} 
\end{equation} 

\begin{figure}[H] 
 \centering
 \includegraphics[width=0.3\textwidth]{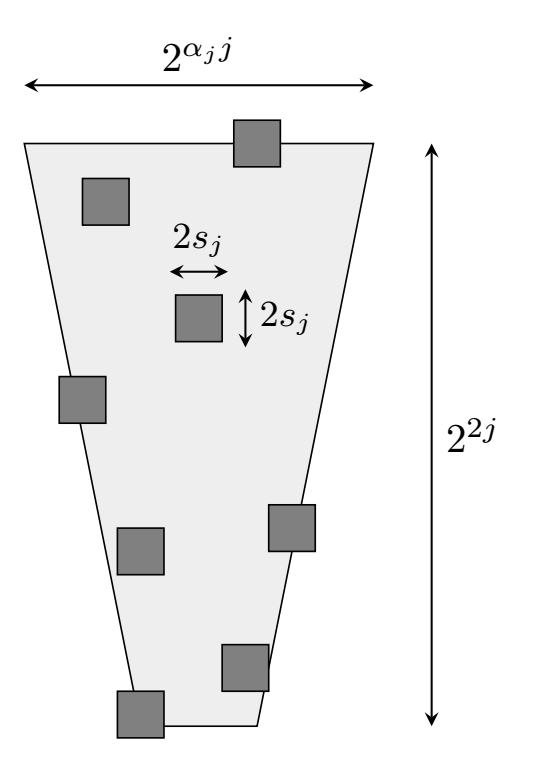}
\;
  \includegraphics[width=0.47\textwidth]{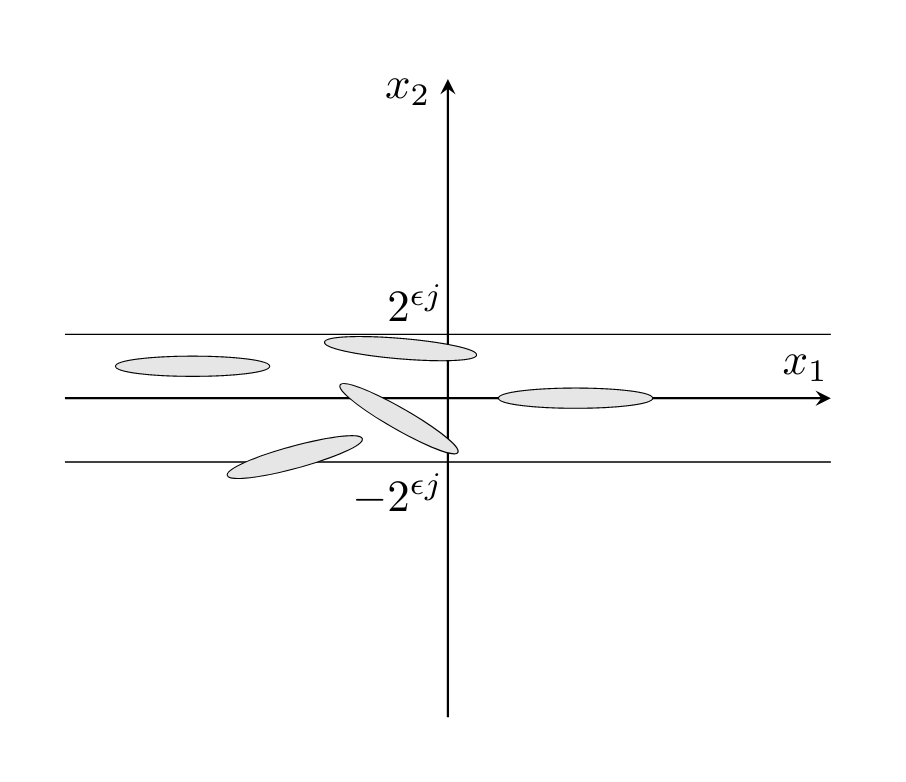}
  \caption{\label{FG6} Left: Interaction between a cluster of Gabor elements (black small squares) and a shearlet (gray) in the frequency domain; Right:  Some shearlet elements (gray) associated with $\Lambda_{1, j}$  in the spatial domain for $l=0$. }
\end{figure}

\subsection{Cluster sets for texture and cartoon \label{Re_sub05}} 
\label{The separation and inpainting theorem for cartoon and texture3}  

We define the cluster for texture at scale $j$ to be
\begin{equation} \label{CT7}
\Lambda_{2,j} := \Big ( \mathbb{Z}^2 \cap B(0,M_j) \Big )  \times \Big ( I_T^{\pm} \cap \mathcal{A}_{s,j} \Big ),
\end{equation}
where $M_j:= 2^{\epsilon j/6}$ and $B(0,r)$ denotes the closed $l_2$ ball around the origin in $\mathbb{R}^2$. The term $M_j= 2^{\epsilon j/6}$  controls the trade-off  between the relative sparsity and cluster coherence of the Gabor systems.

To represent the cartoon part, consider a universal shearlet system with $\alpha_j$ from (\ref{EQ120}), where $\alpha \in (0,2)$ is a global constant.
To define the cluster sets, we fix a constant $\epsilon$ satisfying 
\begin{equation} \label{EQ112}
     0 < \epsilon <   \frac{2-\alpha}{3} .
\end{equation}
Since $\alpha_j$ satisfies $\lim_{j \rightarrow \infty} \alpha_j = \alpha$,  we have  \textcolor{black}{$\lim_{j \rightarrow +\infty } \alpha_j \leq 2-3\epsilon$ and }
\begin{equation} \label{EQ130}
   0 < \epsilon <   \frac{2-\alpha_j}{3}, 
\end{equation} at fine enough scales $j$.
For the cartoon model $w\mathcal{S}_j$  we  define the set of significant coefficients of universal shearlet system by 
\begin{equation} \label{CT6}
\Lambda_{1,j}^{\pm}:= \Lambda_{1,j-1} \cup \Lambda_{1,j} \cup \Lambda_{1,j+1}, \quad  \forall j \geq 2,
\end{equation}
where
$\Lambda_{1,j}: = \Big \{ (j,l,k; \alpha_j, {\rm v}) \; | \; |l| \leq 1, k=(k_1, k_2) \in \mathbb{Z}^2, |k_2 - lk_1| \leq 2^{\epsilon j} \Big \}.$

\subsection{The separation and inpainting theorem for cartoon and texture}  
\label{The separation and inpainting theorem for cartoon and texture2}  

 We now present our main result. 
 In the following theorem, we prove the success of separation and inpaining via Algorithm \textsc{(Inp-Sep)}. Since we inpaint a band of width $h_j$ for each scale $j$, it is important to show that the inpainting error is asymptotically smaller than the energy that typical cartoon and texture parts have in the missing band. For the cartoon part, we can prove that the relative reconstruction error, restricted to the missing part, goes to zero as $j\rightarrow\infty$. For texture, we note that it is possible for $\mathcal{T}_{s,j}$ to be close to zero in the missing part, and thus the relative error restricted to the missing part need not go to zero in the general case. However, generic texture parts are not close to zero in the missing band, and typically 
$$ \| P_j \mathcal{T}_{s,j} \|_2 \propto h_j \| \mathcal{T}_{s,j} \|_2 \propto h_j 2^{-j}.$$
 We thus consider for texture a relative error of the form
$$ \frac{ \| P_j\mathcal{T}_j^\star - P_j\mathcal{T}_{s,j}  \|_2}{h_j 2^{-j}} .$$

\begin{th} \label{TR10} 
Consider the cartoon patch and texture with $w\mathcal{S}_j$ and $\mathcal{T}_{s,j}$  defined in (\ref{DN4}) and (\ref{DN5}) respectively. 
Suppose that the energy matching \textup{(\ref{EQ1})} holds. Suppose that the  index set of the texture $\mathcal{T}_{s,j}$  satisfies \textup{(\ref{EQ2})} and \textup{(\ref{EQ20})}.
Then, for $0<h_j = o(2^{-(\alpha_j + \epsilon)j})$ with $\alpha_j \in (0,2), \liminf \alpha_j > 0$ and $\epsilon$ satisfying (\ref{EQ130}), the recovery error provided by Algorithm \textsc{(Inp-Sep)} decays rapidly and we  have asymptotically perfect simultaneous separation and inpainting. Namely, for all $ N \in \mathbb{N}_0,$

\begin{equation} 
 \frac{\| \mathcal{C}_j^\star - w\mathcal{S}_j  \|_2+ \| \mathcal{T}_j^\star - \mathcal{T}_{s,j}  \|_2}{\| w\mathcal{S}_j \|_2 + \| \mathcal{T}_{s,j} \|_2} = o(2^{-Nj}) \rightarrow 0, \quad j \rightarrow \infty,  \label{EQ201}
\end{equation}
where $(\mathcal{C}_j^\star, \mathcal{T}_j^\star)$ is the solution of \textsc{(Inp-Sep)} and $(\mathcal{C}_j, \mathcal{T}_{s,j})$ are ground truth components.    In addition, 
 if $ \omega(0) \neq 0$, we have asymptotically accurate relative reconstruction error in the missing part
\begin{equation} 
\frac{\| P_j\mathcal{C}_j^\star - P_j\omega\mathcal{S}_j  \|_2}{\| P_j \omega \mathcal{S}_j \|_2} \xrightarrow  {j \rightarrow +\infty} 0 \quad \text{ and }  \quad
 \frac{ \| P_j \mathcal{T}_j^\star - P_j \mathcal{T}_{s,j}  \|_2}{h_j 2^{-j}} \xrightarrow  {j \rightarrow +\infty} 0.  \label{EQ21}
\end{equation}
\end{th}
We postpone the proof of Theorem \ref{TR10} to  Section \ref{S10}, after we discuss preliminary material.

Note that in Theorem \ref{TR10} there is no direct restriction on the texture patch sizes $s_j^{-1}$, and the theorem works even when the texture patch at each scale is smaller than the scale $2^{-j}$. However, the most useful case of Theorem \ref{TR10}, which appropriately models the relation between texture and cartoon, is when $s_j=s$ is constant for all $j$. This is formulated in the following corollary. 
\begin{Cor} Consider the conditions of Theorem \ref{TR10}, and suppose that 
the texture part $\mathcal{T}_{s,j}$  is based on a constant $s_j = s$. 
Then,  we have asymptotically perfect separation and inpainting, i.e, for $\forall N \in \mathbb{N}_0$,
\begin{equation}
    \frac{\| \mathcal{C}_j^\star - \mathcal{C}_j  \|_2+ \| \mathcal{T}_j^\star - \mathcal{T}_{s,j}  \|_2}{\| \mathcal{C}_j \|_2 + \| \mathcal{T}_{s,j} \|_2} =o(2^{-Nj}) \rightarrow 0, \quad j \rightarrow \infty,
\end{equation} 
where $(\mathcal{C}_j^\star, \mathcal{T}_j^\star)$ is the solution of \textsc{(Inp-Sep)} and $(\mathcal{C}_j, \mathcal{T}_{s,j})$ are ground truth components.   In addition, if $\omega(0) \neq 0$, we have
\begin{equation} 
\frac{\| P_j\mathcal{C}_j^\star - P_j\omega\mathcal{S}_j  \|_2}{\| P_j \omega \mathcal{S}_j \|_2} \xrightarrow  {j \rightarrow +\infty} 0 \quad \text{ and }  \quad 
 \frac{ \| P_j \mathcal{T}_j^\star - P_j \mathcal{T}_{s,j}  \|_2}{h_j2^{-j}} \xrightarrow  {j \rightarrow +\infty} 0.  
\end{equation}
\end{Cor}

%\textcolor{black}{Remark: Let us describe one way to use our analysis in a non-multi-scale manner. Let $(1-P)f$ be the the image in the known part. Since low frequencies are assigned automatically to the cartoon part, we first consider low frequency and high frequency filtered versions of the known image. This can be done by filtering $(1-P)f$ to low and high frequency bands $\tilde{f}_l=F_l*(1-P)f$ and $\tilde{f}_h=F_h*(1-P)f$. Here, we assume that the filter $F_l$ keeps the band $j\in(0,j_0]$ and the filter $F_h$ keeps the band $j\in(j_0,j_1)$ \textcolor{magenta}{I think we should use $j_1 = +\infty$} \textcolor{blue}{\large [If $j_1=\infty$ then $h_{j_1}=0$ and there is no missing part.]} 
%\textcolor{magenta}{missing part in high frequency now is limited by $j_o$ not $j_1$, we will make $j_o$ goes to infinity, i.e, instead of considering $f_0, f_1.....$, we consider two sub-images like $f_0 + f_1+...+ f_n$ and $f_{n+1} + f_{n+2}+....$ based on reconstruction formula.}
%[WRITE MORE PROPERLY] .}

%\vskip 2mm
%\hline
%\vskip 2mm

\section{Extensions and future directions \label{S8}}

In this section, we present potential extensions and future directions of our approach.
\begin{enumerate}[label*=\arabic*.]
\item \textbf{Global cartoon model}.
Using the technique introduced in \cite{10} and \cite{13}, we can localize a global cartoon part $\mathcal{C}_j$ by a \emph{partition of unity} $(w_Q)_{Q \in \mathcal{Q}} $. Denoting $ \mathcal{C}_{j,Q} = \mathcal{C}_j \cdot w_Q $, we have
$$ \sum_Q \mathcal{C}_{j,Q} = \mathcal{C}_{j}.$$
In this approach, using a \emph{tubular neighborhood theorem} (\cite[Sect. 6]{10} and \cite{13}), we apply a diffeomorphism $\Phi_Q$ to each piece $\mathcal{C}_{j,Q}$ to straighten out the local curve discontinuity. 
By combining the different neighborhoods, we can derive a convergence results for a global cartoon part.
\item \textbf{Other types of components}.
Our theoretical analysis holds for general representation systems which form  Parseval frames. Thus, our general technique can be applied to problems of image separation and image ipainting  of other types of image parts.
\item \textbf{Noisy case}. In future work we will  extend our theory to the case where the image contains noise.
\item \textbf{More components}. In our analysis we consider  the case where there are two components. In future work we will  extend our theoretical guarantee for separation and inpainting of more than two geometric components. 
\end{enumerate}

In the rest of the paper we prove Theorem \ref{TR10}. For that, in section \ref{S6} we study the relative sparsity of texture and cartoon, and in Section \ref{S7} we study the cluster coherence of the sparsifying systems.

\section{Relative sparsity of texture and cartoon} \label{S6}

In this section, we bound the relative sparsity of texture with respect to shearlet frame (Subsection \ref{Sub11}) and texture with respect to Gabor frame (Subsection \ref{Sub12}).

\subsection{Decay estimates of shearlet and Gabor elements}
First, we provide decay estimates of shearlet and Gabor elements.

\begin{Lemma} \label{LM10} 
Consider the shearlet frame $\boldsymbol{\Psi}$ of Definition  \ref{Df4} and the Gabor frame of scale $j, \mathbf{G_j}$  defined in (\ref{DN1}).
For any arbitrary integer $ N =1,2, \dots $ there exists a constant $C_N$ independent of $j$ such that the following estimates hold
\begin{enumerate}[label*=\arabic*.]
\item $ |(g_{s_j})_{m,n}(x)| \leq C_N \cdot s_j \cdot  \langle |s_j x_1 + \frac{m_1}{2}|\rangle^{-N}\langle |s_j x_2 + \frac{m_2}{2}|\rangle^{-N},$ 
\item $ | \psi_{j, l, k}^{\alpha_j, \rm{(v)}}(x)| \leq C_N \cdot 2^{(2+\alpha_j)j/2} \cdot \langle |2^{\alpha_j j}x_1 - k_1|\rangle^{-N} \langle |2^{2j}x_2+ l2^{\alpha_j j } x_1-k_2|\rangle^{-N}, $ \\
$ | \psi_{j, l, k}^{\alpha_j, \rm{(h)}}(x)| \leq C_N \cdot 2^{(2+\alpha_j) j/2} \cdot \langle |2^{2j}x_1+l2^{\alpha_j j}x_2-k_1)|\rangle^{-N} \langle |2^{\alpha_j j}x_2-k_2|\rangle^{-N}, $ 
\item  $ |\langle (g_{s_j})_{m,n},  \psi_{j, l, k}^{\alpha_j, \rm{(\iota)}} \rangle|  \leq C_N \cdot 2^{-(2-\alpha_j)j/2},$ \quad   $\forall 
 \rm{(\iota) \in \{ \rm{(v)}, \rm{(h)}, \rm{(b)} \}}.$
\end{enumerate}
\end{Lemma}

\begin{proof}
1) By the change of variable $\zeta =s_j^{-1}\xi - n$, we have 
\begin{eqnarray*}
|(g_{s_j})_{m,n}(x)| & =& \Big | \int_{\mathbb{R}^2} (\hat{g}_{s_j})_{m,n}(\xi) e^{2\pi i \xi^T x} d\xi\Big |  =  \Big | \int_{\mathbb{R}^2} s_j^{-1}\hat{g}( s_j^{-1} \xi -n) e^{2\pi i \xi^T (x+\frac{m}{2s_j})} d\xi\Big | \\
& = &  \Big |  s_j e^{2 \pi i n^T(s_j x+ \frac{m}{2})}\int_{\mathbb{R}^2}  \hat{g}(\zeta ) e^{2\pi i \zeta^T (s_j x+\frac{m}{2})} d\zeta \Big | \\
& \leq & \Big |  s_j \int_{\mathbb{R}^2}  \hat{g}(\zeta ) e^{2\pi i \zeta^T (s_j x+\frac{m}{2})} d\zeta \Big |.
\end{eqnarray*}
We now apply integration by parts for $N_1, N_2=1,2, \dots, $  with respect to $\zeta_1, \zeta_2$, respectively, we obtain
\begin{eqnarray*}
&& |(g_{s_j})_{m,n}(x)| \\
& =&  \Big |  \int_{\mathbb{R}^2} s_j(s_j x_1 + \frac{m_1}{2})^{-N_1} \frac{\partial^{N_1}}{\partial \zeta_1^{N_1}} [\hat{g}(\zeta )] e^{2\pi i \zeta^T (s_j x+\frac{m}{2})} d\zeta \Big | \\
& = &  \Big |  \int_{\mathbb{R}^2} s_j(s_j x_1 + \frac{m_1}{2})^{-N_1}(s_j x_2 + \frac{m_2}{2})^{-N_2} \frac{\partial^{N_1+N_2}}{\partial \zeta_1^{N_1} \partial \zeta_2^{N_2}} [\hat{g}(\zeta )] e^{2\pi i \zeta^T (s_j x+\frac{m}{2})} d\zeta \Big | \\
& \leq & s_j| s_j x_1 + \frac{m_1}{2}|^{-N_1} |s_j x_2 + \frac{m_2}{2}|^{-N_2}  \textcolor{black}{\int_{\mathbb{R}^2} \Big | \frac{\partial^{N_1+N_2}}{\partial \zeta_1^{N_1} \partial \zeta_2^{N_2}} [\hat{g}(\zeta )] \Big | d\zeta} ,
\end{eqnarray*}
and similarly
$$
|(g_{s_j})_{m,n}(x)| 
 \leq  s_j| s_j x_k + \frac{m_k}{2}|^{-N_k}  \textcolor{black}{\int_{\mathbb{R}^2} \Big | \frac{\partial^{N_k}}{\partial \zeta_k^{N_k} } [\hat{g}(\zeta )] \Big | d\zeta }
$$
for $k=1,2$.
Here, the boundary terms vanish due to the compact support of $(\hat{g}_{s_j})_{m,n}(\zeta)$.
Thus, \\
$ \phantom{1111} s_j^{-1}\Big ( 1+ |s_j x_1 + \frac{m_1}{2}|^{N_1}+ |s_j x_2 + \frac{m_2}{2}|^{N_2} + |s_j x_1 +  \frac{m_1}{2}|^{N_1}|s_j x_2 + \frac{m_2}{2}|^{N_2} \Big )  |(g_{s_j})_{m,n}(x)| $
\begin{eqnarray*}
 & = & s_j^{-1}(1+ |s_j x_1 + \frac{m_1}{2}|^{N_1})(1+ |s_j x_2 + \frac{m_2}{2}|^{N_2} )  |(g_{s_j})_{m,n}(x)|  \\
& \leq & \textcolor{black}{ \int_{\mathbb{R}^2}  | \hat{g}(\zeta )   | d\zeta +  \int_{\mathbb{R}^2} \Big | \frac{\partial^{N_1}}{\partial \zeta_1^{N_1}} [\hat{g}(\zeta )] \Big | d\zeta +  \int_{\mathbb{R}^2} \Big | \frac{\partial^{N_2}}{\partial \zeta_2^{N_2}} [\hat{g}(\zeta )] \Big | d\zeta } \\
&& \textcolor{black}{+  \int_{\mathbb{R}^2} \Big | \frac{\partial^{N_1+N_2}}{\partial \zeta_1^{N_1} \partial \zeta_2^{N_2}} [\hat{g}(\zeta )] \Big | d\zeta .}
\end{eqnarray*}
By the smoothness of $\hat{g}(\zeta)$ and $\supp \hat{g} = [-1,1],$ there exists a constant $C_{N_1, N_2}^\prime$ independent of $j$ such that \\
 $ \textcolor{black}{\int_{\mathbb{R}^2} |\hat{g}(\zeta ) | d\zeta + \int_{\mathbb{R}^2}  \Big | \frac{\partial^{N_1}}{\partial \zeta_1^{N_1}} [\hat{g}(\zeta )]  \Big | d\zeta+ \int_{\mathbb{R}^2}  \Big | \frac{\partial^{N_2}}{\partial \zeta_2^{N_2}} [\hat{g}(\zeta )] \Big | d\zeta +  \int_{\mathbb{R}^2} \Big | \frac{\partial^{N_1+N_2}}{\partial \zeta_1^{N_1} \partial \zeta_2^{N_2}} [\hat{g}(\zeta )] \Big | d\zeta } \leq C_{N_1, N_2}^\prime.$

We obtain
$$ |(g_{s_j})_{m,n}(x)| \leq C_{N_1,N_2}^\prime \cdot s_j \cdot \frac{1}{1+ |s_j x_1 + \frac{m_1}{2}|^{N_1}}\frac{1}{1+ |s_j x_2 + \frac{m_2}{2}|^{N_2}}.$$
Moreover, for each $N_1, N_2=1,2, \dots $ there exist $C_{N_1}^\prime, C_{N_2}^\prime $ such that 
$$\langle | s_j x_1 + \frac{m_1}{2} | \rangle^{N_1}= (1+ | s_j x_1 + \frac{m_1}{2} |^2 )^{N_1/2 } \leq C_{N_1}^\prime (1 + | s_j x_1 + \frac{m_1}{2} |)^{N_1},$$
$$\langle | s_j x_2 + \frac{m_2}{2} | \rangle^{N_2}= (1+ | s_j x_2 + \frac{m_2}{2} |^2 )^{N_2/2 } \leq C_{N_2}^\prime (1 + | s_j x_2 + \frac{m_2}{2} |)^{N_2}.$$
Thus,
$$ |(g_{s_j})_{m,n}(x)| \leq C_{N_1,N_2} \cdot s_j \cdot \langle |s_j x_1 + \frac{m_1}{2}|\rangle^{-N_1}\langle |s_j x_2 + \frac{m_2}{2}|\rangle^{-N_2}.$$
This proves the claim for any arbitrary integer $N$. \\
2) By the change of variable $\zeta^T =\xi^T A_{\alpha_j, \rm{(\iota)}}^{-j} S^{-l}_{\rm{(\iota)}},$ we have
\begin{eqnarray*}
| \psi_{j, l, k}^{\alpha_j, \rm{(\iota)}}(x)| &=& \Big | \int_{\mathbb{R}^2} \hat{\psi}_{j,l,k}^{\alpha_j, \rm{(\iota)}}(\xi) e^{2\pi i \xi^T x} d\xi\Big | \\
&= & \Big | \int_{\mathbb{R}^2}  2^{-(2+\alpha_j)j/2}W_j(\xi) V_{\rm{(\iota)}} \Big ( \xi^T A_{\alpha_j, \rm{(\iota)}}^{-j} S^{-l}_{\rm{(\iota)}} \Big ) e^{2\pi i \xi^T (x - A_{\alpha_j, \rm{(\iota)}}^{-j} S^{-l}_{\rm{(\iota)}}k )}  d\xi\Big | \\
&= & \Big | \int_{\mathbb{R}^2}  2^{(2+\alpha_j)j/2}W_j((S^{l}_{\rm{(\iota)}}A_{\alpha_j, \rm{(\iota)}}^{j} )^T \zeta) V_{\rm{(\iota)}} (\zeta) e^{2\pi i \zeta^T (S^{l}_{\rm{(\iota)}}A_{\alpha_j, \rm{(\iota)}}^{j} x - k )}  d\zeta \Big |. \\
\end{eqnarray*}
Similarly to in 1) we apply integration by parts for $N_1, N_2 = 1,2, \dots $ with respect to $\zeta_1, \zeta_2$. We obtain the decay estimate of universal shearlets
$$ | \psi_{j, l, k}^{\alpha_j, \rm{(v)}}(x)| \leq C_N \cdot 2^{(2+\alpha_j)j/2} \cdot \langle |2^{\alpha_j j}x_1 - k_1|\rangle^{-N} \langle |2^{2j}x_2+ l2^{\alpha_j j } x_1-k_2|\rangle^{-N}, $$

$$ | \psi_{j, l, k}^{\alpha_j, \rm{(h)}}(x)| \leq C_N \cdot 2^{(2+\alpha_j)j/2} \cdot \langle |2^{2j}x_1+l2^{\alpha_j j}x_2-k_1)|\rangle^{-N} \langle |2^{\alpha_j j}x_2-k_2|\rangle^{-N}. $$
3) We consider three cases. \\
Case 1: $\rm{(\iota)}= \rm{(v)}$. 
Applying 1) and 2) and the change of variables $(y_1, y_2)=(s_j x_1, 2^{2j}x_2+ l2^{\alpha_j j } x_1) $, we obtain 
\begin{eqnarray*}
 |\langle (g_{s_j})_{m,n},  \psi_{j, l, k}^{\alpha_j, \rm{(v)}} \rangle| & = & \Big | \int_{\mathbb{R}^2} (g_{s_j})_{m,n}(x) \overline{\psi_{j, l, k}^{\alpha_j, \rm{(v)}}(x)} dx \Big |\\
 & \leq &  C_N^\prime \cdot   2^{-(2-\alpha_j)j/2} \int_{\mathbb{R}^2} \langle |y_1 + \frac{m_1}{2}|\rangle^{-N}  \langle |y_2-k_2|\rangle^{-N}  dy_1 dy_2\\
  & \leq &  C_N \cdot   2^{-(2-\alpha_j)j/2}.
\end{eqnarray*}
Case 2: $\rm{(\iota)}= \rm{(h)}$. Similarly, we can use 1), 2) and the change of varibles $(y_1, y_2)=(2^{2j} x_1+ l2^{\alpha_j j}x_2, s_j x_2) $ to obtain  $ |\langle (g_{s_j})_{m,n},  \psi_{j, l, k}^{\alpha_j, \rm{(v)}} \rangle| \leq C_N \cdot 2^{-(2-\alpha_j)j/2}.$ \\
Case 3: $\rm{(\iota)}= \rm{(b)}.$ By the definition of boundary shearlets, we can verify this estimate by combining two cases above.  
\end{proof}
 
 %\textcolor{blue}{[Since you delete the proof, don't say that we prove. Say that we state a lemma. Fix this for every claim for which you deleted the proof.]}

%\textcolor{black}{Next, we state a useful lemma for estimating bounds of the type of Proposition \ref{PR4} and Proposition \ref{PR1}.}

Next, we prove a useful lemma for estimating bounds of the type of Proposition \ref{PR4} and Proposition \ref{PR1}.  
\begin{Lemma} \label{LM01}
For each \textcolor{black}{$N \in\mathbb{N}$} there exists a constant $C_N>0$ such that
 $$\int_{\mathbb{R}} \langle |z| \rangle^{-N} \langle | z + t | \rangle^{-N} dz 
\leq  C_N \langle |t| \rangle^{-N}, \quad  \forall t \in \mathbb{R}.$$ 
\end{Lemma}

\begin{proof}
We have
\begin{eqnarray*}
 \int_{\mathbb{R}} \langle |z| \rangle^{-N} \langle | z + t | \rangle^{-N} dz 
&=& \int_{\mathbb{R}} \max \{ \langle |z| \rangle^{-N},  \langle | z + t | \rangle^{-N}\}  \min \{ \langle |z| \rangle^{-N} , \langle | z + t | \rangle^{-N}\} dz \nonumber \\
& \leq & \int_{\mathbb{R}} \Big ( \langle |z| \rangle^{-N} + \langle | z + t | \rangle^{-N} \Big ) \cdot \langle |t/2| \rangle^{-N} dz \nonumber \\
& \leq & C_N \langle |t| \rangle^{-N}.
\end{eqnarray*}
\end{proof} 

\subsection{Cartoon patch} \label{Sub11}

For the sake of brevity, we use some indexing sets for universal shearlets 
\begin{eqnarray} 
\Delta &:=& \Big \{(j,l,k; \alpha_j, \rm{(\iota)}) \mid j \geq 0, |l| < 2^{(2-\alpha_j )j}, k \in \mathbb{Z}^2, \rm{(\iota)} \in \{ \rm{(h)}, \rm{(v)} \}  \Big \}  \nonumber \\
& & \bigcup \Big \{(j,l,k; \alpha_j, \rm{(b)}) \mid  j \geq 0, |l| = 2^{(2-\alpha_j)j}, k \in \mathbb{Z}^2 \Big \},  \label{CT60} \\
\Delta_j &:=& \{ (j^\prime,k,l; \alpha_j, \rm{(\iota)}) \in \Delta \mid j^\prime = j \}, \quad j \geq 0,  \label{CT61} \\
\Delta_j^{\pm} &:=& \Delta_{j-1} \cup \Delta_j \cup \Delta_{j+1}, \label{CT62}
\end{eqnarray} 
where $\Delta_{-1}=\emptyset$. We now have the following result.
\begin{Prop} 
\label{PR4}
Consider the shearlet frame $\boldsymbol{\Psi}$ of Definition  \ref{Df4} and the cartoon patch $w \mathcal{S}_j $ as defined in (\ref{DN4}). We assume that $\liminf_{j \rightarrow \infty} \alpha_j >0$. Then, we have the following decay estimate of the cluster approximate error $\delta_{1,j}$ 
\begin{equation} \label{CT9}
\delta_{1,j}:= \sum_{\eta \in \Delta, \eta \notin \Lambda_{1,j}^\pm} | \langle \psi_{j^\prime,l,k}^{\alpha_{j^\prime}, \rm{(\iota)}}, w \mathcal{S}_j \rangle | = o(2^{-Nj}),  \quad \forall N \in \mathbb{N},
\end{equation}
where $\eta=(j^\prime, l, k; \alpha_{j'}, \rm{(\iota)}).$    
\end{Prop} 
 The proof of this proposition roughly follows the lines of the Proposition 5.2 in \cite{25}. For the sake of brevity we denote 
\begin{equation} \label{CT40}
    t^{(\rm{v})} = (t_1^{(\rm{v})}, t_2^{(\rm{v})}): = A^{-j}_{\alpha_j, (\rm{v})} S^{-l}_{(\rm{v})} k= (2^{-\alpha_j j}k_1, 2^{-2 j} (k_2-lk_1)), 
\end{equation} 
\begin{equation}  \label{CT41}
  t^{(\rm{h})} = (t_1^{(\rm{h})}, t_2^{(\rm{h})}): = A^{-j}_{\alpha_j, (\rm{h})} S^{-l}_{(\rm{h})} k= (2^{-2j}(k_1-lk_2), 2^{-\alpha_j j} k_2).  
\end{equation}
The proof also relies on the following  two  lemmas.

\begin{Lemma} \label{LM100}
Consider the shearlet frame $\boldsymbol{\Psi}$ of Definition  \ref{Df4} and the cartoon patch with sub-band $w\mathcal{S}_j$ defined in (\ref{DN4}). Then, for $\alpha_j \in (0,2), j^\prime \in \{ j-1, j, j+1  \}, j \geq 2, $ the following estimates hold for arbitrary integers $ M \geq 0$ 
\begin{enumerate}[label*=\arabic*.]
\item
If $(\rm{\iota}) = (\rm{v})$ and $|l|>1$, we have 
$$\Big | \langle w\mathcal{S}_j, \psi_{j^\prime,l,k}^{\alpha_{j^\prime}, \rm{(v)}} \rangle \Big | \leq C_{ M } \cdot 2^{-(2+\alpha_{j^\prime})j^\prime/2} \cdot \langle |t_1^{\rm{(v)}} | \rangle^{- 2} \langle |t_2^{\rm{(v)}} |\rangle^{-2}   \langle | 2^{\alpha_{j^\prime} j^\prime} | \rangle^{-M}.  $$
\item
If $(\rm{\iota}) = (\rm{h})$, we have 
$$\Big | \langle w\mathcal{S}_j, \psi_{j^\prime,l,k}^{\alpha_{j^\prime}, \rm{(h)}} \rangle \Big | \leq C_{M } \cdot \frac{2^{2j^\prime}}{r^3} \cdot \langle |t_1^{\rm{(h)}} | \rangle^{- 2} \langle|t_2^{\rm{(h)}} |\rangle^{-2}   \langle | 2^{2 j^\prime} | \rangle^{-M},   $$
where $r$ is defined in (\ref{EQ153}).
\item
If $(\rm{\iota}) =\rm{(b)}$ and $|l| = 2^{(2-\alpha_{j^\prime}) j^\prime},$ we have
$$\Big | \langle w\mathcal{S}_j, \psi_{j^\prime,l,k}^{\alpha_{j^\prime}, \rm{(b)}} \rangle \Big | =  C_{ M } \cdot 2^{-(2+\alpha_{j^\prime})j^\prime/2} \cdot \langle|t_1^{\rm{(b)}} |\rangle^{- 2} \langle |t_2^{\rm{(b)}} |\rangle^{-2}    \langle | 2^{\alpha_{j^\prime} j^\prime} | \rangle^{-M} . $$
\end{enumerate}
\end{Lemma} 

\begin{proof} 
Without loss of generality, we  prove for $j^\prime = j$. The other cases are treated similarly.

1) By the definition of $w \mathcal{S}_j$ and Plancherel's theorem, we obtain
\begin{eqnarray*}
\Big | \langle w\mathcal{S}_j, \psi_{j,l,k}^{\alpha_j, \rm{(v)}} \rangle \Big | & =& \Big | \langle \widehat{w\mathcal{S}}_j, \hat{\psi}_{j,l,k}^{\alpha_j, \rm{(v)}} \rangle \Big | = \Big | \int_{\mathbb{R}^2} \hat{w}(\xi_1) \frac{1}{2\pi\xi_2}W_j(\xi) \overline{\hat{\psi}^{\alpha_j, \rm{(v)}}_{j,l,k}(\xi)} \Big |\\
& = & \frac{1}{2\pi}\Big | \int_\mathbb{R}e^{2\pi i t_2^{\rm{(v)}}\xi_2} \int_\mathbb{R}\hat{w}(\xi_1)\frac{1}{\xi_2}W_j(\xi)\hat{\psi}^{\alpha_j, \rm{(v)}}_{j,l,0}(\xi)e^{2\pi it_1^{\rm{(v)}}\xi_1}d\xi_1 d\xi_2 \Big | \\
&=&\frac{1}{2\pi} \Big |  \int_\mathbb{R}e^{2\pi i t_2^{\rm{(v)}}\xi_2} \int_\mathbb{R} \hat{w}(\xi_1) \tau_{j,l}(\xi)e^{2\pi it_1^{\rm{(v)}}\xi_1}d\xi_1 d\xi_2 \Big |,
\end{eqnarray*}
where $\tau_{j,l}(\xi):= \frac{1}{\xi_2} W_j(\xi)\hat{\psi}^{\alpha_j, \rm{(v)}}_{j,l,0}(\xi) $.

We now apply repeated integration by parts  with respect to $\xi_i$, $i=1,2,$. We get
\begin{eqnarray}
  \Big | \langle w\mathcal{S}_j, \psi_{j,l,k}^{\alpha_j, \rm{(v)}} \rangle \Big |  &=& \frac{1}{2\pi} \Big |  \int_\mathbb{R}e^{2\pi i t_2^{\rm{(v)}}\xi_2} \int_\mathbb{R} \frac{\partial^2}{\partial \xi_1^{2}} [\hat{w}(\xi_1) \tau_{j,l}(\xi)]\frac{1}{(2\pi i t_1^{\rm{(v)}})^{2}}e^{2\pi it_1^{\rm{(v)}}\xi_1}d\xi_1 \Big ] d\xi_2 \Big |  \nonumber \\
  & \leq & \frac{1}{|t_1^{\rm{(v)}}|^2} \Big | \int_\mathbb{R}  e^{2\pi i t_2^{\rm{(v)}}\xi_2} \int_\mathbb{R} \frac{\partial^2}{\partial \xi_1^{2}} [\hat{w}(\xi_1) \tau_{j,l}(\xi)]e^{2\pi it_1^{\rm{(v)}}\xi_1} d\xi_1 d\xi_2 \Big | \nonumber \\ 
 & \leq & |t_1^{\rm{(v)}}|^{-2} \int_\mathbb{R} \int_\mathbb{R} \Big | \frac{\partial^2}{\partial \xi_1^{2}} [\hat{w}(\xi_1) \tau_{j,l}(\xi)] \Big | d\xi_1 d \xi_2  \label{CT70} 
\end{eqnarray}
 and similarly
\begin{eqnarray}
 \Big | \langle w\mathcal{S}_j, \psi_{j,l,k}^{\alpha_j, \rm{(v)}} \rangle \Big | & \leq &  |t_2^{\rm{(v)}}|^{-2} \int_\mathbb{R} \int_\mathbb{R} \Big | \frac{\partial^2}{\partial \xi_2^{2}} [\hat{w}(\xi_1) \tau_{j,l}(\xi)] \Big | d\xi_1 d \xi_2   \label{CT71}
 \end{eqnarray}
\begin{eqnarray}
 \Big | \langle w\mathcal{S}_j, \psi_{j,l,k}^{\alpha_j, \rm{(v)}} \rangle \Big | & \leq & |t_1^{\rm{(v)}}|^{-2} |t_2^{\rm{(v)}}|^{-2} \int_\mathbb{R} \int_\mathbb{R} \Big | \frac{\partial^4}{\partial \xi_1^{2}\partial \xi_2^{2}} [\hat{w}(\xi_1) \tau_{j,l}(\xi)] \Big | d\xi_1 d \xi_2.  \label{CT72}
\end{eqnarray}
Thus,
$$
 |t_1^{\rm{(v)}}|^{2}  \Big | \langle w\mathcal{S}_j, \psi_{j,l,k}^{\alpha_j, \rm{(v)}} \rangle \Big | \stackrel{\mathclap{\normalfont{(\ref{CT70})}}}{ \leq }   \int_\mathbb{R} \int_\mathbb{R} \Big | \frac{\partial^2}{\partial \xi_1^{2}} [\hat{w}(\xi_1) \tau_{j,l}(\xi)] \Big | d\xi_1 d \xi_2  $$
 
 $$ |t_2^{\rm{(v)}}|^{2}  \Big | \langle w\mathcal{S}_j, \psi_{j,l,k}^{\alpha_j, \rm{(v)}} \rangle \Big |  \stackrel{\mathclap{\normalfont{(\ref{CT71})}}}{ \leq }    \int_\mathbb{R} \int_\mathbb{R} \Big | \frac{\partial^2}{\partial \xi_2^{2}} [\hat{w}(\xi_1) \tau_{j,l}(\xi)] \Big | d\xi_1 d \xi_2  $$
 
  $$|t_1^{\rm{(v)}}|^{2} |t_2^{\rm{(v)}}|^{2}  \Big | \langle w\mathcal{S}_j, \psi_{j,l,k}^{\alpha_j, \rm{(v)}} \rangle \Big | \stackrel{\mathclap{\normalfont{(\ref{CT72})}}}{ \leq }   \int_\mathbb{R} \int_\mathbb{R} \Big | \frac{\partial^4}{\partial \xi_1^{2}\partial \xi_2^{2}} [\hat{w}(\xi_1) \tau_{j,l}(\xi)] \Big | d\xi_1 d \xi_2. $$
  This leads to 
  $$ (1+ |t_1^{\rm{(v)}}|^{2} + |t_2^{\rm{(v)}}|^{2}+ |t_1^{\rm{(v)}}|^{2} |t_2^{\rm{(v)}}|^{2})  \Big | \langle w\mathcal{S}_j, \psi_{j,l,k}^{\alpha_j, \rm{(v)}} \rangle \Big | =   \langle  |t_1^{\rm{(v)}}| \rangle^{2}  \langle | t_2^{\rm{(v)}}| \rangle^{2} \Big | \langle w\mathcal{S}_j, \psi_{j,l,k}^{\alpha_j, \rm{(v)}} \rangle \Big |  $$
 \begin{eqnarray}
    &\leq&  \int_{\mathbb{R}^2} \Big |  \hat{w}(\xi_1) \tau_{j,l}(\xi) \Big |+ \Big | \frac{\partial^2}{\partial \xi_1^{2}} [\hat{w}(\xi_1) \tau_{j,l}(\xi)] \Big |+\Big | \frac{\partial^2}{\partial \xi_2^{2}} [\hat{w}(\xi_1) \tau_{j,l}(\xi)] \Big | \qquad  \qquad \qquad \nonumber \\ && + \Big | \frac{\partial^4}{\partial \xi_1^{2}\partial \xi_2^{2}} [\hat{w}(\xi_1) \tau_{j,l}(\xi)] \Big | d\xi \label{EQ70} \\
    &= & \int_{\mathbb{R}} \Gamma(\xi_2) d\xi_2, \nonumber
 \end{eqnarray}
 
where  
$\Gamma(\xi_2):= \\
\int_{\mathbb{R}} \Big |  \hat{w}(\xi_1) \tau_{j,l}(\xi) \Big |+ \Big | \frac{\partial^2}{\partial \xi_1^{2}} [\hat{w}(\xi_1) \tau_{j,l}(\xi)] \Big |+\Big | \frac{\partial^2}{\partial \xi_2^{2}} [\hat{w}(\xi_1) \tau_{j,l}(\xi)] \Big |+ \Big | \frac{\partial^4}{\partial \xi_1^{2}\partial \xi_2^{2}} [\hat{w}(\xi_1) \tau_{j,l}(\xi)] \Big | d\xi_1.$
Next, we bound $\Gamma$.
By the definition of the universal shearlets,  $\hat{\psi}_{j,l,k}^{\alpha_j, \rm{(v)}} $ has compact support in the trapezoidal region 
 \begin{equation}
     {\rm supp}(\hat{\psi}_{j,l,k}^{\alpha_j, \rm{(v)}})=\Big \{  \xi \in \mathbb{R}^2 \ | \ \xi_2 \in [-2^{2j-1}, 2^{2j-1}] \setminus [-2^{2j-4}, 2^{2j-4}], \Big | \frac{\xi_1}{\xi_2} - l2^{-(2-\alpha_j )j} \Big | \Big \}.
     \label{eq:supp1}
 \end{equation}
 This implies that for any $\xi \in \supp \tau_{j,l}$ we have
 $$ (l-1)2^{(\alpha_j -2)j } \leq \frac{\xi_1}{\xi_2} \leq (l+1)2^{(\alpha_j -2)j} \text{ and } 2^{2j-4} \leq |\xi_2| \leq 2^{2j-1}.$$ 
Using the assumption $|l| > 1$ follows that there exist constants $C_1, C_2 >0$ such that
$$ \xi_1 \in I_{j,l}:= [-C_2l2^{\alpha_jj },-C_1l2^{\alpha_jj }] \cup  [C_2l2^{\alpha_jj },C_1l2^{\alpha_jj }].$$

We now give some further estimates. By the fact that $\tau_{j,l}(\xi)$ is zero for $\xi_2 < 2^{2j-4}$, we obtain
\begin{eqnarray}
2^{(2+ \alpha_j)j/2}\Big |\frac{ \partial \tau_{j,l}}{\partial \xi_1} (\xi) \Big | & = & \Big | \frac{\partial}{\partial \xi_1} \Big [ \frac{1}{\xi_2} W^2_j( \xi) v\Big ( 2^{(2-\alpha_j )j} \frac{\xi_1}{\xi_2} - l \Big ) \Big ] \Big |  \nonumber \\
& \leq &  \frac{1}{2^{2j-4}} \Big [  \frac{2}{2^{2j}} W(\frac{\xi_1}{2^{2j}}, \frac{\xi_2}{2^{2j}})\frac{\partial W}{\partial \xi_1}(\frac{\xi_1}{2^{2j}}, \frac{\xi_2}{2^{2j}})  v\Big ( 2^{(2-\alpha_j )j} \frac{\xi_1}{\xi_2} - l \Big ) +   \nonumber \\
& &
\frac{2^{(2-\alpha_j)j}}{\xi_2} \frac{\partial v}{\partial \xi_1}{\Big ( 2^{(2-\alpha_j )j} \frac{\xi_1}{\xi_2} - l \Big )} W^2(\frac{\xi_1}{2^{2j}}, \frac{\xi_2}{2^{2j}}) \Big ]  \nonumber \\
& \leq & \frac{1}{2^{2j-4}} \Big ( C_1 2^{-2j} + C_2 \frac{2^{(2-\alpha_j)j}}{|\xi_2|} \Big ) \stackrel{\mathclap{\normalfont{0< \alpha_j \leq 2 }}}{ \quad\leq  \quad} C2^{-2 j }. \label{EQ73}
\end{eqnarray}
Furthermore, we have
$$ 2^{(2+ \alpha_j)j/2}\Big |\frac{ \partial \tau_{j,l}}{\partial \xi_2} (\xi) \Big | = 2^{(2+ \alpha_j)j/2} \Big | \frac{ \partial \frac{\sigma_{j,l}}{\xi_2}}{\partial \xi_2} \Big |, $$
where $\sigma_{j,l}: = W^2_j( \xi) v\Big ( 2^{(2-\alpha_j )j} \frac{\xi_1}{\xi_2} - l \Big ).$
Thus, in the support (\ref{eq:supp1}) we have
\begin{eqnarray*}
  && 2^{(2+ \alpha_j)j/2} \Big | \frac{ \partial \sigma_{j,l}}{\partial \xi_2} \Big |  =  \Big | \frac{\partial}{\partial \xi_2} \Big [ W^2_j( \xi) v\Big ( 2^{(2-\alpha_j )j} \frac{\xi_1}{\xi_2} - l \Big ) \Big ] \Big | \\
& = &  \frac{2}{2^{2j}} W(\frac{\xi_1}{2^{2j}}, \frac{\xi_2}{2^{2j}})\frac{\partial W}{\partial \xi_2}(\frac{\xi_1}{2^{2j}}, \frac{\xi_2}{2^{2j}})  v\Big ( 2^{(2-\alpha_j )j} \frac{\xi_1}{\xi_2} - l \Big ) -  \frac{2^{(2-\alpha_j)j} \xi_1}{\xi_2^2} \frac{\partial v}{\partial \xi_1} W^2(\frac{\xi_1}{2^{2j}}, \frac{\xi_2}{2^{2j}})\\
& \leq & C_1 2^{-2j} + C_2 \frac{2^{(2-\alpha_j)j}|\xi_1|}{|\xi_2|^2} \stackrel{\mathclap{\normalfont{\alpha_j \leq 2 }}}{\quad \leq \quad} C2^{-\alpha_j j }.
\end{eqnarray*}
Consequently, we get
\begin{equation}
     2^{(2+ \alpha_j)j/2}\Big |\frac{ \partial \tau_{j,l}}{\partial \xi_2} (\xi) \Big | = 2^{(2+ \alpha_j)j/2} \Big | \frac{ \partial \sigma_{j,l}}{\partial \xi_2} \frac{1}{\xi_2} - \frac{\sigma_{j, l}(\xi)}{\xi_2^2} \Big |  \leq \frac{C2^{-\alpha_j j } }{2^{2j-4}} +  \frac{C^\prime}{2^{2j-4}} \leq  C^{\prime\prime}2^{-2j}. \label{EQ74}
\end{equation}
Similarly, we obtain
\begin{equation}
    2^{(2+ \alpha_j)j/2}\Big |\frac{ \partial^2 \tau_{j,l}}{\partial \xi_1^2} (\xi) \Big | \leq  C2^{-2 j }, \label{EQ75}
\end{equation}  
and 
\begin{equation}
2^{(2+ \alpha_j)j/2}\Big |\frac{ \partial^2 \tau_{j,l}}{\partial \xi_2^2} (\xi) \Big | \leq  C2^{-2 j }. \label{EQ76}
\end{equation} 
 
\begin{figure}[H]
  \includegraphics[width=0.4\textwidth]{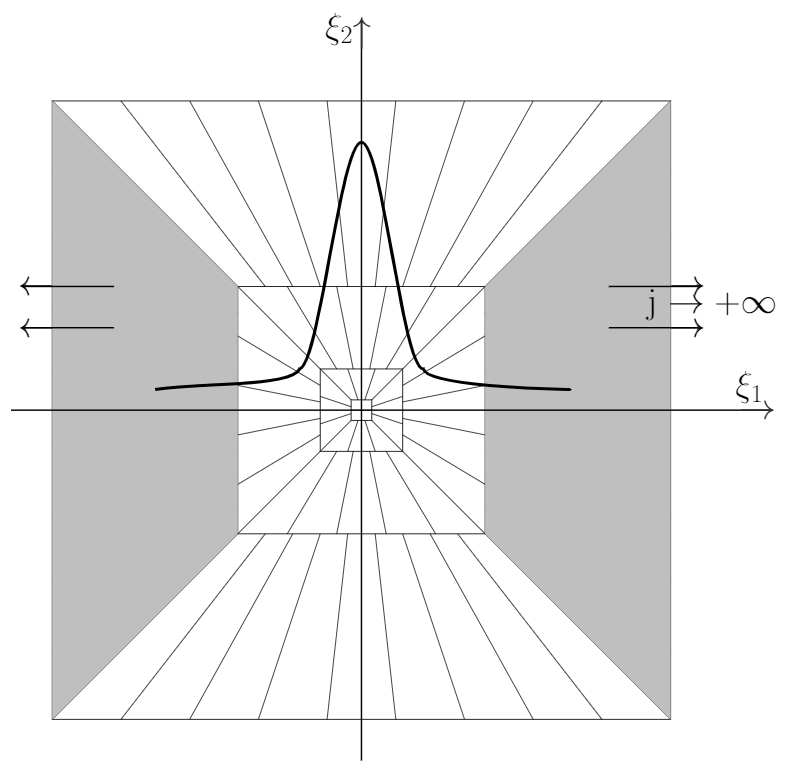}
  \hfill
  \includegraphics[width=0.4\textwidth]{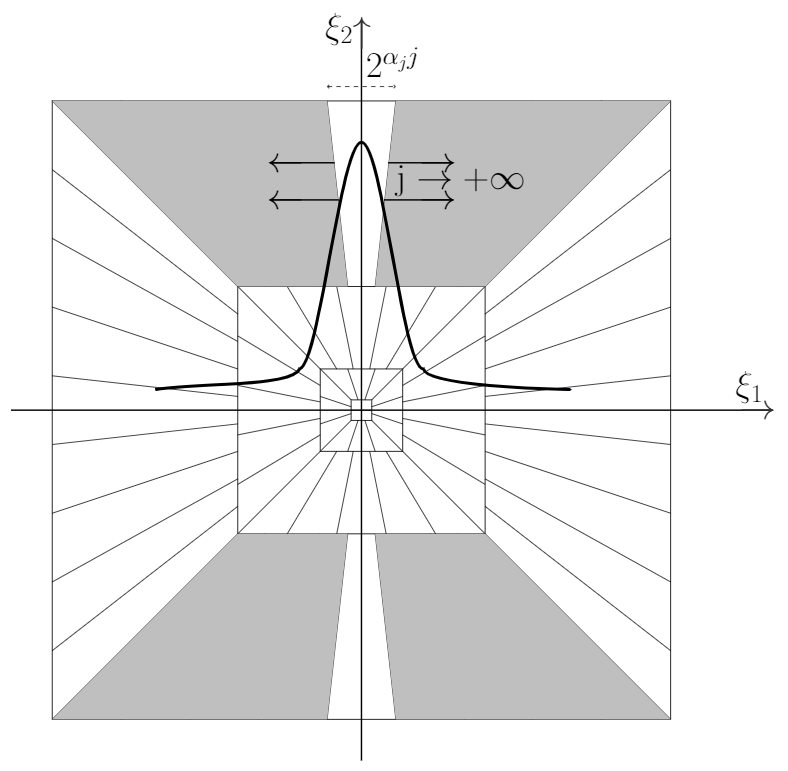}
  \caption{\label{FG20} Interactions between horizontal shearlets, vertical shearlets  and the sub-image of cartoon part $w\mathcal{S}_j$. Left: Frequency support of horizontal shearlets at scale $j$ (gray), $w\mathcal{S}_j$ (bold). Right: Frequency support of vertical shearlets with $|l| >1$ (gray), $w\mathcal{S}_j$ (bold). }
\end{figure}

 We now exploit  the specific form of $I_{j,l}$  as well as the rapid decay of $\hat{w}$ to study  $\hat{w}^{(n)}$ for $n=0,1,2$. These lead to (see Figure \ref{FG20} for illustration) 
\begin{eqnarray}
 \|  \hat{w}^{(n)} \|_{L^1(I_{j,l})} &\leq& vol(I_{j,l}) \sup_{\xi_1 \in I_{j,l}} | \hat{w}^{(n)}(\xi_1)| \leq C_{M, n} | l2^{\alpha_j j}| \langle |l2^{\alpha_j j} | \rangle^{-(M+1)} \nonumber \\
 & \leq & C_{M, n} \langle |  l2^{\alpha_j j } |\rangle^{-M} \stackrel{\mathclap{\normalfont{l > 1}}}{ \; \leq \;}  C_{M, n} \langle |  2^{\alpha_j j } |\rangle^{-M},  \; \forall M \in \mathbb{N}, \forall n = 1,2, \dots \qquad  \; \label{EQ77}
 \end{eqnarray}
Combining   (\ref{EQ73}), (\ref{EQ74}), (\ref{EQ75}), (\ref{EQ76}), (\ref{EQ77}),  with (\ref{EQ70}), we have
\begin{eqnarray}
 \Gamma(\xi_2) &\leq & C\cdot 2^{-(2+ \alpha_j)j/2}  \cdot 2^{-2j} \max \Big \{ \|\hat{w} \|_{L^1(I_{j,l})}, \|\frac{\partial \hat{w}}{\partial \xi_1} \|_{L^1(I_{j,l})}, \|\frac{\partial^2 \hat{w}}{\partial \xi_1^2} \|_{L^1(I_{j,l})} \Big \} \qquad \\
 &\leq & C_{M} \cdot 2^{-(2+ \alpha_j)j/2} \ \cdot 2^{-2j} \cdot \langle |  2^{\alpha_j j } |\rangle^{-M}.
\end{eqnarray}
In addition, since $\supp \Gamma(\xi_2)= [ 2^{2j-4}, 2^{2j-1}]$,
we finally obtain
\begin{eqnarray*}
   \Big | \langle w\mathcal{S}_j, \psi_{j,l,k}^{\alpha_j, \rm{(v)}} \rangle \Big |  &\leq&  C_{M} \cdot 2^{-(2+\alpha_j)j/2} \cdot  \langle |t_1^{\rm{(v)}}| \rangle^{-2}  \langle | t_2^{\rm{(v)}}| \rangle^{-2}  \langle |  2^{\alpha_j j } |\rangle^{-M} .
\end{eqnarray*}

2)  This claim is similar to 1) but we need to modify the intervals $I_{j,l} $ to be 
$ [ -C_2 2^{2j}, -C_1 2^{2j}] \cup [ C_1 2^{2j}, C_2 2^{2j}]$, where $C_1, C_2 >0$.  This leads to the additional terms $\langle | 2^{2j}| \rangle^{-M}$ and $r$ in 2) when comparing with 1). Next, we use $r \leq |\xi_2| \leq 2^{2j-1}$
instead of using $2^{2j-4} \leq |\xi_2| \leq 2^{2j-1}$ for the support of $\tau_{j,l}$. Thus, 
redefining $\tau_{j,l}(\xi)= \frac{1}{\xi_2} W_j(\xi)\hat{\psi}^{\alpha_j, \rm{(h)}}_{j,l,0}(\xi) $, we modify
$(\ref{EQ73}), (\ref{EQ74}), (\ref{EQ75}),(\ref{EQ76})$ as follows.

\begin{eqnarray}
2^{(2+ \alpha_j)j/2}\Big |\frac{ \partial \tau_{j,l}}{\partial \xi_1} (\xi) \Big | & = & \Big | \frac{\partial}{\partial \xi_1} \Big [ \frac{1}{\xi_2} W^2_j( \xi) v\Big ( 2^{(2-\alpha_j )j} \frac{\xi_2}{\xi_1} - l \Big ) \Big ] \Big |  \nonumber \\
& \leq &  \frac{1}{r} \Big |  \frac{2}{2^{2j}} W(\frac{\xi_1}{2^{2j}}, \frac{\xi_2}{2^{2j}})\frac{\partial W}{\partial \xi_1}(\frac{\xi_1}{2^{2j}}, \frac{\xi_2}{2^{2j}})  v\Big ( 2^{(2-\alpha_j )j} \frac{\xi_2}{\xi_1} - l \Big )  \Big | +   \nonumber \\
& & \Big |
\frac{2^{(2-\alpha_j)j}}{\xi_1^2} \frac{\partial v}{\partial \xi_1} W^2(\frac{\xi_1}{2^{2j}}, \frac{\xi_2}{2^{2j}}) \Big | \nonumber \\
& \leq & \frac{1}{r} \Big ( C_1 2^{-2j} + C_2 \frac{2^{(2-\alpha_j)j}}{|\xi_1|^2} \Big ) \stackrel{\mathclap{\normalfont{0< \alpha_j \leq 2 }}}{ \quad \leq \quad }  \; \frac{C}{ r \cdot 2^{2j}}. \label{EQ83}
\end{eqnarray}
Furthermore, we have
$$ 2^{(2+ \alpha_j)j/2}\Big |\frac{ \partial \tau_{j,l}}{\partial \xi_2} (\xi) \Big | = 2^{(2+ \alpha_j)j/2} \Big | \frac{ \partial \frac{\sigma_{j,l}}{\xi_2}}{\partial \xi_2} \Big |, $$
where $\sigma_{j,l}: = W^2_j( \xi) v\Big ( 2^{(2-\alpha_j )j} \frac{\xi_2}{\xi_1} - l \Big ).$
Thus, 
\begin{eqnarray*}
  && 2^{(2+ \alpha_j)j/2} \Big | \frac{ \partial \sigma_{j,l}}{\partial \xi_2} \Big |  =  \Big | \frac{\partial}{\partial \xi_2} \Big [ W^2_j( \xi) v\Big ( 2^{(2-\alpha_j )j} \frac{\xi_2}{\xi_1} - l \Big ) \Big ] \Big | \\
& \leq &  \Big | \frac{2}{2^{2j}} W(\frac{\xi_1}{2^{2j}}, \frac{\xi_2}{2^{2j}})\frac{\partial W}{\partial \xi_2}(\frac{\xi_1}{2^{2j}}, \frac{\xi_2}{2^{2j}})  v\Big ( 2^{(2-\alpha_j )j} \frac{\xi_2}{\xi_1} - l \Big )  \Big | + \Big |   \frac{2^{(2-\alpha_j)j}}{\xi_1} \frac{\partial v}{\partial \xi_1} W^2(\frac{\xi_1}{2^{2j}}, \frac{\xi_2}{2^{2j}}) \Big | \\
& \leq & C_1 2^{-2j} + C_2 \frac{2^{(2-\alpha_j)j}}{|\xi_1|} \stackrel{\mathclap{\normalfont{0< \alpha_j \leq 2 }}}{ \quad \; \leq \quad \; } C .
\end{eqnarray*}
Consequently, we have
\begin{equation}  \label{EQ84}
     \Big |\frac{ \partial \tau_{j,l}}{\partial \xi_2} (\xi) \Big | =  \Big | \frac{ \partial \sigma_{j,l}}{\partial \xi_2} \frac{1}{\xi_2} - \frac{\sigma_{j, l}(\xi)}{\xi_2^2} \Big |  \leq \frac{C^\prime }{r^2}. 
\end{equation} 

Similarly, we obtain
\begin{equation}
    2^{(2+ \alpha_j)j/2}\Big |\frac{ \partial^2 \tau_{j,l}}{\partial \xi_1^2} (\xi) \Big | \leq  \frac{C}{r \cdot 2^{2j}}, \label{EQ85}
\end{equation}  
and
\begin{equation}
\Big |\frac{ \partial^2 \tau_{j,l}}{\partial \xi_2^2} (\xi) \Big | \leq  \frac{C^\prime}{r^3}. \label{EQ86}
\end{equation} 
These bounds lead to
$$\Gamma(\xi_2) \leq \frac{C}{r^3}  \cdot \max \Big \{ \|\hat{w} \|_{L^1(I_{j,l})}, \|\frac{\partial \hat{w}}{\partial \xi_1} \|_{L^1(I_{j,l})}, \|\frac{\partial^2 \hat{w}}{\partial \xi_1^2} \|_{L^1(I_{j,l})} \Big \}. $$

Thus, by the rapid decay of $\hat{w}$, as before, 
\begin{eqnarray*}
   \Big | \langle w\mathcal{S}_j, \psi_{j,l,k}^{\alpha_j, \rm{(v)}} \rangle \Big |  &\leq&  C_{M} \cdot  \frac{2^{2j}}{r^3} \cdot \langle |t_1^{\rm{(v)}}| \rangle^{-2}  \langle | t_2^{\rm{(v)}}| \rangle^{-2}  \langle |  2^{2 j } |\rangle^{-M} . 
\end{eqnarray*}
 For an illustration, we refer to Figure \ref{FG20}. Finally, we obtain
3). Recalling the definition of the shearlets, we can easily verify this estimate since boundary shearlets is just a particular case of horizontal and vertical shearlets.
\end{proof}

\begin{Lemma} \label{LM101} Consider the shearlet frame $\boldsymbol{\Psi}$ of Definition  \ref{Df4} and the cartoon patch with sub-band $w\mathcal{S}_j$ defined in (\ref{DN4}). For $j\in \{ j-1, j , j+1 \}, j \geq 2, $ and for any arbitrary $N \geq 2$, there exists a constant $C_N>0$ such that 
\begin{eqnarray*}
  \Big | \langle w\mathcal{S}_j, \psi_{j^ \prime,l,k}^{\alpha_{j^\prime}, \rm{(v)}} \rangle \Big | &\leq&  C_N 2^{(2-\alpha_{j^\prime})j^\prime/2}\int_\mathbb{R} \tilde{w}_{N,j}(2^{-\alpha_{j^\prime} j^\prime}(x_1 + k_1)) \langle | x_1 | \rangle^{-N}  \\ 
  && \langle |lx_1 + lk_1 - k_2 | \rangle^{-N} dx_1,
\end{eqnarray*}
where $\tilde{w}_{N,j^\prime} := | w | * \langle |2^{2j^\prime} [\cdot ] |\rangle^{-N}$, and $\langle | \cdot | \rangle$ is defined from $\rm{(\ref{CT2})}$ .
\end{Lemma}
\begin{proof}
Without loss of generality, we can prove only for $j^\prime = j.$ We first note that in $\supp \psi_{j,l,k}^{\alpha_j, \rm{(v)}} $ we have $w \mathcal{S}_0 \equiv w \mathcal{S},$ where $wS$ is defined in (\ref{EQ153}). We now
consider the line distribution $w \mathcal{L}$ introduced in  \cite{25} by
\begin{equation} \label{EQ13}
  \langle w \mathcal{L} , f \rangle = \int_{-\rho}^{\rho} w(x_1) f(x_1, 0) dx_1, \quad f \in \mathbb{S}(\mathbb{R}^2),
\end{equation} 
 with Fourier transform 
\begin{equation} \label{EQ14}
\langle \widehat{w \mathcal{L}} , f \rangle = \int_{\mathbb{R}^2} \hat{w}(\xi_1) f(\xi) d\xi.
\end{equation}

Next, we observe that
  \begin{align*}
  \langle w \mathcal{L} , f \rangle & = \int_{-\rho}^{\rho} w(x_1) f(x_1, 0) dx_1   =  \int_{-\rho}^{\rho} w(x_1) \int_\mathbb{R} f(x_1, x_2) \delta(x_2) dx_2 dx_1 \\
  & = \iint_{\mathbb{R}^2} w(x_1) \delta(x_2)f(x_1, x_2) dx_1 dx_2 = \langle w(x_1) \delta(x_2), f \rangle, 
  \end{align*}
 where $\delta(x) $ denotes the usual Dirac delta functional. Thus, 
 \begin{equation} \label{EQ15}
 w \mathcal{L}(x) = w(x_1) \delta(x_2). 
\end{equation} 
In addition, 
\begin{align*}
 \langle w \mathcal{S} , f \rangle & = \int_{-\rho}^{\rho} w(x_1) \int_{-\infty}^0 f(x) dx_2 dx_1   = \int_{-\rho}^{\rho} w(x_1) \int_{-\infty}^{+\infty} f(x) \mathds{1}_{ \{ 0 \geq x_2 \} }(0) dx_2 dx_1   \\
  & = \int_{-\rho}^{\rho} w(x_1) \int_{\mathbb{R}} f(x) \Big ( \int_{\mathbb{R}}\mathds{1}_{ \{ y \geq x_2 \} }(y) \delta(y) dy \Big ) dx_2 dx_1  \\
   & =  \int_{-\rho}^{\rho} w(x_1) \int_{\mathbb{R}} f(x) \Big ( \int_{x_2}^{+\infty} \delta(y) dy \Big ) dx_2 dx_1  = \Big \langle \int_{x_2}^{+\infty} w(x_1) \delta(y) dy, f(x) \Big \rangle.
\end{align*}
Hence, we obtain 
\begin{equation} \label{EQ16}
w \mathcal{S }(x)= \int_{x_2}^{+ \infty } w\mathcal{L}(x_1, y) dy.
\end{equation}  

Since integration commutes with convolution, we also have $$w\mathcal{S}_j = w \mathcal{S} * F_j = \int_{x_2}^{+\infty} w\mathcal{L}_j(x_1,y) dy.$$
where $w \mathcal{L}_j$ is defined by $w \mathcal{L} * F_j$.
In addition, by integration by parts with respect to variable $x_2$, we obtain
$$ \langle w \mathcal{S}_j, \psi_{j, l, k}^{\alpha_j, \rm{(v)}} \rangle = \langle \int_{x_2}^{+ \infty } w\mathcal{L}_j(x_1, y) dy, \psi_{j, l, k}^{\alpha_j, \rm{(v)}} \rangle \rangle = \langle  w\mathcal{L}_j(x) , \int_{x_2}^{+ \infty } \psi_{j, l, k}^{\alpha_j, \rm{(v)}}(x_1,y) dy \rangle, $$ where the boundary terms vanish due to the compact support of $ \psi_{j, l, k}^{\alpha_j, \rm{(v)}} .$ 
Now we put
$$ \Xi_{j, l, k}^{\alpha_j, \rm{(v)}}(x):= \int_{x_2}^{+\infty} \psi_{j, l, k}^{\alpha_j, \rm{(v)}}(x_1,y) dy,$$
and note that $\hat{\Xi}_{j, l, k}^{\alpha_j, \rm{(v)}}(\xi) = \frac{i}{2\pi\xi_2}\cdot \hat{\psi}_{j, l, k}^{\alpha_j, \rm{(v)}}(\xi).$ Similarly to  Lemma \ref{LM10}, which was based on $\psi_{j, l, k}^{\alpha_j, \rm{(v)}}$, we can prove a  lemma for $\Xi_{j, l, k}^{\alpha_j, \rm{(v)}}$, and prove the rapid decay property
\begin{eqnarray}
|\Xi_{j, l, k}^{\alpha_j, \rm{(v)}}(x) | & \leq &  C_N^\prime \cdot \frac{1}{2^{2j-4}} \cdot  2^{(2+ \alpha_j)j/2} \cdot \langle |2^{\alpha_jj} x_1 -k_1  | \rangle^{-N} \langle | l2^{\alpha_jj} x_1 + 2^{2j}x_2- k_2  |\rangle^{-N} \nonumber \\
& \leq &  C_N \cdot 2^{(\alpha_j -2)j/2} \cdot \langle |2^{\alpha_jj} x_1 -k_1  | \rangle^{-N} \langle | l2^{\alpha_jj} x_1 + 2^{2j}x_2- k_2  |\rangle^{-N}. \label{CT100}
\end{eqnarray}

Now we can use the decay estimate of the line singularities $w \mathcal{L}_j $ to obtain
\begin{eqnarray*}
| w \mathcal{L}_j (x) | &=& | [w \mathcal{L} * F_j](x) | = \Big | \int_\mathbb{R} w(y_1) F_j(x-(y_1, 0)) dy_1 \Big |  \\
& \leq & \Big | \int_\mathbb{R} |w(y_1) |2^{4j}| \check{W}( 2^{2j}(x-(y_1, 0))) | dy_1 \Big | \\
& \leq &   \int_\mathbb{R} |w(y_1) |2^{4j} C_N \langle | 2^{2j}x_2|\rangle^{-N}  \langle | 2^{2j}(y_1 -x_1)|\rangle^{-N}  dy_1  \\
&=& C_N2^{4j} \langle | 2^{2j} x_2 | \rangle^{-N} [ |w * \langle 2^{2j}[\cdot]| \rangle^{-N} ](x_1) \\
& = & C_N2^{4j} \langle | 2^{2j} x_2 | \rangle^{-N}  \tilde{w}_{N, j}(x_1), 
\end{eqnarray*}
where $\tilde{w}_{N, j}(x_1): =   [ |w * \langle 2^{2j}[\cdot]| \rangle^{-N} ](x_1), \quad x=(x_1, x_2) \in \mathbb{R}^2.$  
Combining this with the rapid decay properties of  $\Xi_{j, l, k}^{\alpha_j, \rm{(v)}}(x)$ from $(\ref{CT100})$ and \textup{Lemma \ref{LM01}}, we have
\begin{eqnarray*}
&& \Big | \langle  w \mathcal{S}_j , \psi_{j, l, k}^{\alpha_j, \rm{(v)}} \rangle  \Big |   \\
&\leq & C_N \int_{\mathbb{R}^2} 2^{2j} \langle | 2^{2j}x_2 |\rangle^{-N} \tilde{w}_{N,j}(x_1) 2^{(\alpha_j +2 )j/2} \langle |2^{\alpha_j j} x_1 - k_1 | \rangle^{-N}  \\
&& \qquad \quad   \cdot \langle | l 2^{\alpha_j j }x_1 + 2^{2j} x_2 - k_2 | \rangle^{-N} dx \\
& = &  C_N 2^{(2-\alpha_j ) j/2}  \int_{\mathbb{R}^2}  \langle |x_2 |\rangle^{-N} \tilde{w}_{N,j}(2^{-\alpha_j j}(x_1 + k_1)) \langle |x_1| \rangle^{-N} \\
&& \qquad \quad \qquad \quad \quad \cdot \langle | lx_1 + lk_1+ x_2 - k_2 | \rangle^{-N} dx \\
&  \leq & C_N 2^{(2-\alpha_j)j/2} \int_{\mathbb{R}} \tilde{w}_{N,j}(2^{-\alpha_j j}(x_1 + k_1)) \langle |x_1| \rangle^{-N} \langle | lx_1 + lk_1 -k_2 | \rangle^{-N} dx_1 .
\end{eqnarray*}
\end{proof}

Now we are ready to prove Proposition  \ref{PR4}
\begin{proof}[Proof of Proposition  \ref{PR4}]
By the definition of $\delta_{1,j}$ and the support of the window function $F_j$,  only shearlet coefficients in $\Delta_j^{\pm}$ have nonzero inner products with $wS_j$, so we have 
\begin{eqnarray*}
  \delta_{1,j} &=& \sum_{\eta \in \Delta, \eta \notin \Lambda_{1,j}^\pm} | \langle \psi_{j^\prime,l,k}^{\alpha_{j^\prime}, \rm{(\iota)}}, w \mathcal{S}_j \rangle | = \sum_{\eta \in \Delta_j^{\pm}, \eta \notin \Lambda_{1,j}^\pm} | \langle \psi_{j^\prime,l,k}^{\alpha_{j^\prime}, \rm{(\iota)}}, w \mathcal{S}_j \rangle |  \\
  & = & \sum_{\substack{k \in \mathbb{Z}^2, |l| \leq 1 \\ |k_2 -  lk_1 | > 2^{\epsilon j^\prime} }} | \langle \psi_{j^\prime,l,k}^{\alpha_{j^\prime}, \rm{(v)}}, w\mathcal{S}_j \rangle |  +  \sum_{k \in \mathbb{Z}^2, |l| >1 } | \langle \psi_{j^\prime,l,k}^{\alpha_{j^\prime}, \rm{(v)}}, w \mathcal{S}_j \rangle |  \\
  && +  \sum_{k \in \mathbb{Z}^2, l \in \mathbb{Z}} | \langle \psi_{j^\prime,l,k}^{\alpha_{j^\prime}, \rm{(h)}}, w \mathcal{S}_j \rangle |  +  \sum_{k \in \mathbb{Z}^2  } | \langle \psi_{j^\prime,\pm 2^{(2-\alpha_{j^\prime})j^\prime},k}^{\alpha_{j^\prime}, \rm{(b)}}, w \mathcal{S}_j \rangle |  \\
  & = & T_1 + T_2 + T_3 + T_4,
\end{eqnarray*}
where 
$$T_1: = \sum_{\substack{k \in \mathbb{Z}^2, |l|\leq 1 \\ |k_2 -  lk_1 | > 2^{\epsilon j^\prime} }} | \langle \psi_{j^\prime,l,k}^{\alpha_{j^\prime}, \rm{(v)}}, w \mathcal{S}_j \rangle | , \quad T_2: = \sum_{k \in \mathbb{Z}^2, |l| >1 } | \langle \psi_{j^\prime,l,k}^{\alpha_{j^\prime}, \rm{(v)}}, w \mathcal{S}_j \rangle | , $$
$$T_3= \sum_{k \in \mathbb{Z}^2, l \in \mathbb{Z}} | \langle \psi_{j^\prime,l,k}^{\alpha_{j^\prime}, \rm{(h)}}, w \mathcal{S}_j \rangle |, \quad T_4:= \sum_{k \in \mathbb{Z}^2  } | \langle \psi_{j^\prime,\pm 2^{(2-\alpha_{j^\prime})j^\prime},k}^{\alpha_{j^\prime}}, w \mathcal{S}_j \rangle |.  $$
Without loss of generality we restrict to the scale index $j^\prime = j$. The respective arguments for the other cases $ j^\prime = j \pm 1 $ are similar. 

We start with an estimation of $T_1$. By Lemma \ref{LM101}, for $N \geq 2,$ we obtain
\begin{eqnarray*}
&& T_1  = \sum_{\substack{k \in \mathbb{Z}^2, |l|\leq 1 \\ |k_2 -  lk_1 | > 2^{\epsilon j} }} | \langle \psi_{j,l,k}^{\alpha_{j}, \rm{(v)}}, w \mathcal{S}_j \rangle |  \\
& \leq & C_N 2^{(2-\alpha_j)j/2} \sum_{\substack{k \in \mathbb{Z}^2, |l| \leq 1 \\ |k_2 -  lk_1 | > 2^{\epsilon j} }} \int_\mathbb{R} \tilde{w}_{N,j}(2^{-\alpha_j j}(x_1 + k_1)) \langle | x_1 | \rangle^{-N} \langle |lx_1 + lk_1 - k_2 | \rangle^{-N} dx_1 \\
& = & C_N 2^{(2-\alpha_j)j/2} \sum_{\substack{k \in \mathbb{Z}^2, |l|\leq 1 \\ |k_2 | > 2^{\epsilon j} }} \int_\mathbb{R} \tilde{w}_{N,j}(2^{-\alpha_j j}(x_1 + k_1)) \langle | x_1 | \rangle^{-N} \langle |lx_1 + k_2 | \rangle^{-N} dx_1  \\
& = & C_N 2^{(2-\alpha_j)j/2} \sum_{\substack{k_2 \in \mathbb{Z}, |l|\leq 1 \\ |k_2 | > 2^{\epsilon j} }} \int_\mathbb{R} \Big ( \sum_{k_1 \in \mathbb{Z}}\tilde{w}_{N,j}(2^{-\alpha_j j}(x_1 + k_1)) \Big ) \langle | x_1 | \rangle^{-N} \langle |lx_1 + k_2 | \rangle^{-N} dx_1.            
\end{eqnarray*}
Furthermore, we have
\begin{eqnarray*}
\sum_{k_1 \in \mathbb{Z}}\tilde{w}_{N,j}(2^{-\alpha_j j}(x_1 + k_1)) & = & \sum_{k_1 \in \mathbb{Z}} \int_{\mathbb{R}}|w(y_1) | \langle | 2^{2j}(y_1-2^{-\alpha_j j }(x_1+ k_1))| \rangle^{-N}dy_1 \\
& = & \sum_{k_1 \in \mathbb{Z}} \int_{\mathbb{R}} |w(y_1) | \langle | 2^{(2-\alpha_j)j}(k_1 + x_1 - 2^{\alpha_j j }y_1)| \rangle^{-N}dy_1 \\
& \stackrel{\mathclap{\normalfont{\alpha_j \leq 2 }}} \leq & \; \int_{\mathbb{R}}  |w(y_1) | \Big (  \sum_{k_1 \in \mathbb{Z}}\langle |k_1 + x_1 - 2^{\alpha_j j }y_1| \rangle^{-N} \Big ) dy_1 \\
& \leq &  C'\int_{\mathbb{R}}  |w(y_1) | \Big (  \int_\mathbb{R} \langle |t + x_1 - 2^{\alpha_j j }y_1| \rangle^{-N} dt \Big ) dy_1 \\
& = &  C'\int_{\mathbb{R}}  |w(y_1) | \Big (  \int_\mathbb{R} \langle |t| \rangle^{-N} dt \Big ) dy_1  \leq C_N.
\end{eqnarray*}
This implies
\begin{eqnarray*}
 T_1 & \leq &  C_N 2^{(2-\alpha_j)j/2} \sum_{\substack{k_2 \in \mathbb{Z}, |l| \leq 1 \\ |k_2 | > 2^{\epsilon j^\prime} }} \int_\mathbb{R} \langle | x_1 | \rangle^{-N} \langle |lx_1 + k_2 | \rangle^{-N} dx_1. \\
& \stackrel{\mathclap{\normalfont{\text{Lemma } \ref{LM01}}}} \leq &   \quad C_N 2^{(2-\alpha_j)j/2} \sum_{k_2 \in \mathbb{Z},  |k_2 | > 2^{\epsilon j^\prime} }  \langle | k_2 | \rangle^{-N}   dx_1.  \\
& \leq &  C_N 2^{(2-\alpha_j)j/2} \int_{t > 2^{\epsilon j}} \langle | t | \rangle^{-N} \leq  C_N 2^{(2-\alpha_j)j/2} 2^{-(N-1)\epsilon j}.
\end{eqnarray*}

We now bound the decay rate for  the term $T_2$. First we fix $j \geq 0$ and $|l| >1$. For $N \geq 2$, we obtain
\begin{eqnarray}
 \sum_{k \in \mathbb{Z}^2} \langle |t_1^{\rm{(v)}}|\rangle^{-2} | t_2^{\rm{(v)}}|^{-N} &=&  \sum_{k \in \mathbb{Z}^2}  \langle |2^{-\alpha_J   j} k_1 |\rangle^{-2} \langle|2^{-2j}(k_2 - l k_1)|\rangle^{-2} \nonumber \\
 & \leq &  C' 2^{(2+\alpha_j)2j}  \int_{\mathbb{R}} \int_{\mathbb{R}} \langle |x_1 |\rangle^{-2} \langle |x_2|\rangle^{-2}dx_1 dx_2 \leq C  2^{(2+ \alpha_j )2j}. \nonumber \\ \label{CT102}
\end{eqnarray}

Now, by Lemma \ref{LM100}, we have
\begin{eqnarray*}
T_2 &= & \sum_{\substack{k \in \mathbb{Z}^2\\ 1 < |l| \leq 2^{(2-\alpha_j) j}}}| \langle \psi_{j,l,k}^{\alpha_{j}, \rm{(v)}}, w \mathcal{S}_j \rangle |  \\
&\leq& C_{M}   \sum_{\substack{k \in \mathbb{Z}^2\\ 1 < |l| \leq 2^{(2-\alpha_j) j}}} 2^{-(2+\alpha_{j})j/2} \cdot \langle |t_1^{\rm{(v)}} | \rangle^{- 2} \langle |t_2^{\rm{(v)}} |\rangle^{-2}   \langle | 2^{\alpha_{j} j} | \rangle^{-M}  \\
& \leq & C_{M} \cdot  2^{(2 -3\alpha_j)j/2}  \cdot  \langle | 2^{\alpha_j j} | \rangle^{-M} \sum_{\substack{k \in \mathbb{Z}^2\\ }} \langle| t_1^{\rm{(v)}}| \rangle^{-2} \langle| t_2^{\rm{(v)}}|\rangle^{-2}    \\
& \stackrel{\mathclap{\normalfont{ (\ref{CT102})}} }\leq & \; C_{M}^\prime \cdot 2^{(2 -3\alpha_j)j/2}  \cdot  2^{(2+ \alpha_j )2j} \cdot \langle | 2^{\alpha_j j} | \rangle^{-M}.
\end{eqnarray*}
Using the assumption $\liminf_{j \rightarrow \infty} \alpha_j \ > 0$ and choosing $M$ sufficiently large we obtain the desired decay rates. By using 2) and 3) of the Lemma \ref{LM100} the estimates for $T_3$ and $T_4$ are done similarly. Here, note that $\frac{1}{r^3} \cdot o(2^{-Nj})= o(2^{-Nj})$  for fixed  $r > 0$. For an illustration of Lemmas \ref{LM100}, \ref{LM101}, and the main proposition,  we refer to Figure \ref{FG20}.
\end{proof}
\subsection{Texture} \label{Sub12} We now provide a relative sparsity error (Definition \ref{Def4} ) for the texture part.
\begin{Prop} \label{PR1}
Consider the  Gabor frame of scale $j, \mathbf{G_j}$,  defined in (\ref{DN1}), and the texture  defined in Definition \ref{Df3}.
Then, for every $N \in \mathbb{N}$ the sequence $(\delta_{2,j})_{j \in \mathbb{N}}$ decays rapidly in the sense
\begin{equation} \label{CT11} 
\delta_{2,j}:=\sum_{\lambda \notin B(0,M_j) \times I_T^{\pm} } | \langle \mathcal{T}_{s,j}, (g_{s_j})_{\lambda} \rangle |= o(2^{-Nj}).
\end{equation}
\end{Prop}
\begin{proof}
First, note that 
$$ \hat{\mathcal{T}}_{s} =\int_{\mathbb{R}^2} \sum_{n \in I_T} d_n g_{s_j}(x) e^{-2 \pi i \langle \xi - s_jn, x \rangle} dx= \sum_{n \in I_T } d_n \hat{g}_{s_j}(\xi - s_jn).$$
Using the support condition of $\hat{g}$, denoting $\lambda = (\tilde{m}, \tilde{n})$, we obtain 
\begin{eqnarray*}
 | \langle \mathcal{T}_{s,j}, (g_{s_j})_{\lambda} \rangle | & = &| \langle \hat{\mathcal{T}}_{s,j}, (\hat{g}_{s_j})_{\lambda} \rangle |   \\
 &=& \Big |  \sum_{ n \in I_T}   d_n \int W_j(\xi) \hat{g}_{s_j}(\xi -s_jn) \cdot \overline{\hat{g}_{s_j}(\xi - s_j \tilde{n})} e^{-2 \pi i \frac{\langle \tilde{m}, \xi \rangle}{2s_j}} d \xi \Big | \\
  &=& \Big |  \sum_{\substack{ n \in I_T, \\  |n - \tilde{n}| \leq 1}}   d_n \int W_j(\xi) \hat{g}_{s_j}(\xi -s_jn) \cdot \overline{\hat{g}_{s_j}(\xi - s_j \tilde{n})} e^{- \pi i \frac{\langle \tilde{m}, \xi \rangle}{s_j}} d \xi\Big |. \\
\end{eqnarray*}
Indeed, $| \langle \mathcal{T}_{s,j}, (g_{s_j})_{\lambda} \rangle |  = | \langle \hat{\mathcal{T}}_{s,j}, (\hat{g}_{s_j})_{\lambda} \rangle | =0, $ for $\tilde{n} \notin I_T^{\pm} \cap \mathcal{A}_{s,j}$. Moreover, we have \\
$\{ (\tilde{m}, \tilde{n}) \notin B(0,M) \times I_T^{\pm} \} = \{ (\tilde{m}, \tilde{n})  \in \mathbb{Z}^2 \times \mathbb{Z}^2, |\tilde{m}| > M_j \} \cup \{(\tilde{m}, \tilde{n}) \in \mathbb{Z}^2 \times \mathbb{Z}^2, \tilde{n} \notin I_T^{\pm} \}.$
Thus, since the Gabor coefficients of $\mathcal{T}_{s,j}$ are zero for frequencies outside $I_T^{\pm}$,
\begin{eqnarray*}
 \delta_{2,j} &=& \sum_{\lambda \notin B(0,M_j) \times I_T^{\pm} } | \langle \mathcal{T}_{s,j}, (g_{s_j})_{\lambda} \rangle | =  \sum_{|\tilde{m}| \geq M_j, \tilde{n} \in I_T^{\pm}\cap \mathcal{A}_{s,j} } | \langle \mathcal{T}_{s,j}, (g_{s_j})_{\lambda} \rangle |   \\
  &= &   \sum_{|\tilde{m}| \geq M_j, \tilde{n} \in I_T^{\pm} \cap \mathcal{A}_{s,j} } \Big |  \int_{\mathbb{R}^2} \sum_{\substack{ n \in I_T, \\  |n - \tilde{n}| \leq 1}}  d_{n} (g_{s_j})_{(0,n)}(x)  \overline{(g_{s_j})_{(\tilde{m}, \tilde{n})}(x)} dx\Big |  \\
  &\leq &  C \sum_{|\tilde{m}| \geq M_j, \tilde{n} \in I_T^{\pm} \cap \mathcal{A}_{s,j} }  |d_{\tilde{n}}| \int_{\mathbb{R}^2}   |(g_{s_j})_{(0,\tilde{n})}(x)| \cdot | (g_{s_j})_{(\tilde{m}, \tilde{n})}(x)| dx,
\end{eqnarray*}
where $d_{\tilde{n}}= \max\{ |d_n| \mid n \in I_T, |n-\tilde{n}| \leq 1\}$.
Note that the last inequality  holds since for each $\tilde{n}$, there exist a finite number of $n$  satisfying $|n- \tilde{n}| \leq 1.$
Moreover, by (\ref{EQ2}), there are  $o(2^{(2-\alpha_j)j/2})$ points of $I_T^{\pm}$ in each support of a shearlet and there are $o(2^{(2-\alpha_j)j})$ shearlets in $A_{j}$, so 
$$ \# \{ n  \in I_T^{\pm} \cap \mathcal{A}_{s,j} \} \leq 9 \# \{ n \in I_T \cap \mathcal{A}_{s,j} \} = o(2^{(2-\alpha_j)j/2}) \cdot  2^{(2-\alpha_j)j} = o(2^{(2-\alpha_j)3j/2}).$$
Combining this  with the energy matching condition
$$\sum_{n \in  I_T \cap \mathcal{A}_{s,j}   } |d_{n}|^2 \leq  c_1 2^{-2j},$$
which implies $d_n \leq \sqrt{c_1}, \forall n \in I_T \cap  \mathcal{A}_{s,j} ,$ 
we have
\begin{eqnarray*}
 \delta_{2,j} &\leq & C \cdot  o(2^{(2-\alpha_j)3j/2}) \sum_{|\tilde{m}| > M_j } \int_{\mathbb{R}^2}   |(g_{s_j})_{(0,  \tilde{n})}(x)| \cdot | (g_{s_j})_{\tilde{m}, \tilde{n}}(x)| dx. \\
& \stackrel{\mathclap{\normalfont{\textup{Lemma }  \ref{LM10} }}} \leq&  \quad C_N \cdot o(2^{(2-\alpha_j)3j/2})   \int_{\mathbb{R}^2} \sum_{|\tilde{m}| > M_j } s_j^2 \cdot  \langle |s_j x_1|\rangle^{-2N}\langle |s_j x_2|\rangle^{-2N} \langle |s_j x_1 + \frac{\tilde{m}_1}{2}|\rangle^{-N} \\
&& \quad \cdot \langle |s_j x_2 + \frac{\tilde{m}_2}{2}|\rangle^{-N}  dx.
\end{eqnarray*}

By Lemma \ref{LM01}, we now obtain
\begin{eqnarray*}
 \delta_{2,j} &\leq &  C_N \cdot o(2^{(2-\alpha_j)3j/2})  \cdot  \sum_{|\tilde{m}| > M_j }\langle |\frac{\tilde{m}_1}{2}|\rangle^{-N}  \langle |\frac{\tilde{m}_2}{2}|\rangle^{-N}  \\
  &\leq & C_N \cdot o(2^{(2-\alpha_j)3j/2})  \cdot  \int_\mathbb{R} \int_{\frac{1}{2}M_j} \langle \frac{t_1}{2}|\rangle^{-N}  \langle | \frac{t_2}{2}|\rangle^{-N} dt_1 dt_2 \\
& \stackrel{\mathclap{\normalfont{M_j = 2^{\epsilon j/6} }}} \leq & C_N \cdot o(2^{(2-\alpha_j)3j/2})  \cdot 2^{-(N-1)\epsilon j/6}.
\end{eqnarray*}
This proves the claim since we can choose an  arbitrarily large $N \in \mathbb{N}.$
\end{proof}

\section{Cluster coherence for cartoon and texture} \label{S7}
In the previous section we proved relative sparsity for the texture and cartoon parts. For the success of separation and inpainting, by  Theorem \ref{TR1}, we need to prove that the cluster coherence terms are less than $1/2$ when summed together. To guarantee this asymptotically, we prove that the cluster coherence terms decay to zero when the scale goes to infinity. 

\subsection{Cluster coherence of the un-projected frames} We now analyze the cluster coherence of the un-projected frames $\mu_c (\Lambda_{1,j}^{\pm},  \boldsymbol{\Psi}; \mathbf{G_j} ) $ and
$\mu_c(\Lambda_{2,j},  \mathbf{G_j};\boldsymbol{\Psi}).$ 

\begin{Prop} \label{PR3}
Consider the shearlet frame $\boldsymbol{\Psi}$ of Definition  \ref{Df4} and the Gabor frame of scale $j, \mathbf{G_j}$,  defined in (\ref{DN1}). 
We have
$$\mu_c (\Lambda_{1,j}^{\pm},  \boldsymbol{\Psi}; \mathbf{G_j} ) \rightarrow 0, \quad j \rightarrow \infty. $$
where $\Lambda_{1,j}^{\pm}$ is defined in (\ref{CT6}).
\end{Prop}

\begin{proof}
Let us  recall the definition of the cluster of shearlets
$$ \Lambda_{1,j} =\Big \{ (j,l,k; \alpha_j, {\rm(v)}) \mid \lvert l \rvert \leq 1, k \in \mathbb{Z}^2, \lvert k_2-lk_1\rvert \leq 2^{\epsilon j} \Big \}, $$ 
and
$\Lambda_{1,j}^{\pm} = \Lambda_{1,j-1} \cup \Lambda_{1,j} \cup \Lambda_{1,j+1}$ as defined in (\ref{CT6}). Without loss of generality,  we prove only $\mu_c (\Lambda_{1,j},  \boldsymbol{\Psi}   ;  \mathbf{G_j}) \rightarrow 0$ as $j \rightarrow \infty$, and note that summing over $j'=j-1,j,j+1$ does not change the asymptotics.  By the definition of the cluster coherence, we have
$$ \mu_c (\Lambda_{1,j},  \boldsymbol{\Psi}   ; \mathbf{G_j}) = \max_{(m,n)} \sum_{\eta \in \Lambda_{1,j}} |\langle  \psi_{j, l, k}^{\alpha_j, \rm{(\iota)}}, (g_{s_j})_{m,n} \rangle|,$$ where $\eta=(j, l, k; \alpha_j, \rm{(\iota)}).$ 
Suppose that the maximum is attained for  $\bar{m}, \bar{n} \in \mathbb{Z}^2$ (the maximum exists, since both frames are translation invariant, and the shearlet elements at scale $j$ are compactly supported in frequency).  Applying the change of variables $(y_1, y_2) = (s_j x_1, 2^{2j}x_2+ l2^{\alpha_j j } x_1 -lk_1)$ and  Lemma \ref{LM10} 1) and 2) , we obtain
\begin{eqnarray*}
 && \mu_c (\Lambda_{1,j},  \boldsymbol{\Psi}   ;  \mathbf{G_j})  =  \sum_{\eta \in \Lambda_{1,j}} |\langle  \psi_{j, l, k}^{\alpha_j, \rm{(\iota)}}, (g_{s_j})_{\bar{m},\bar{n}} \rangle| \\
 & \leq &  \sum_{\substack{|l| \leq 1, k \in \mathbb{Z}^2 \\ |k_2 - lk_1| \leq 2^{\epsilon j} }} \int_{\mathbb{R}^2}C_N \cdot s_j \cdot 2^{(2+\alpha_j)j/2} \cdot \langle |2^{\alpha_j j}x_1 - k_1|\rangle^{-N} \langle |2^{2j}x_2+ l2^{\alpha_j j } x_1-k_2|\rangle^{-N} \\
&& \qquad \qquad
 \cdot \langle |s_j x_1 + \frac{m_1}{2}|\rangle^{-N}\langle |s_j x_2 + \frac{m_2}{2}|\rangle^{-N}dx_1 dx_2. \\
 & =&   \sum_{\substack{|l| \leq 1, k \in \mathbb{Z}^2 \\ |k_2 - lk_1| \leq 2^{\epsilon j} }} \int_{\mathbb{R}^2}C_N \cdot 2^{-(2-\alpha_j)j/2} \cdot \langle |\frac{2^{\alpha_j j}}{s_j}y_1 - k_1|\rangle^{-N} \langle |y_2+ lk_1-k_2|\rangle^{-N}  \\
&& \qquad \qquad \quad \cdot \langle |y_1 + \frac{m_1}{2}|\rangle^{-N} \langle |\frac{s_j}{2^{2j}} y_2 -l2^{-(2-\alpha_j)j} y_1+s_j2^{-2j}lk_1 +\frac{m_2}{2}|\rangle^{-N}dy_1 dy_2 \\ 
& \leq &   \sum_{\substack{|l| \leq 1, k \in \mathbb{Z}^2\\ |k_2| \leq 2^{\epsilon j} }} \int_{\mathbb{R}^2}C_N \cdot 2^{-(2-\alpha_j)j/2}  \langle |\frac{2^{\alpha_j j}}{s_j}y_1 - k_1|\rangle^{-N} \langle |y_2-k_2|\rangle^{-N} \\ 
&& \qquad \qquad \quad \cdot \langle |y_1 + \frac{m_1}{2}|\rangle^{-N} dy_1dy_2 \\
 & =&  \sum_{\substack{|l| \leq 1, k_2 \in \mathbb{Z}\\ |k_2| \leq 2^{\epsilon j} }} \int_{\mathbb{R}^2}C_N2^{-(2-\alpha_j)j/2}  \Big ( \sum_{k_1 \in \mathbb{Z}}\langle |\frac{2^{\alpha_j j}}{s_j}y_1 - k_1|  \rangle^{-N} \Big ) \langle |y_2-k_2|\rangle^{-N} \\
 && \qquad \qquad 
 \cdot \langle |y_1 + \frac{m_1}{2}|\rangle^{-N}dy_1 dy_2 \\ 
 & \leq & C_N \cdot 2^{-(2-\alpha_j)j/2} \cdot 2^{\epsilon j} = C_N \cdot 2^{-(2-\alpha_j-2\epsilon)j/2} \xrightarrow{j \rightarrow +\infty} 0,
\end{eqnarray*}
where $(2-\alpha_j-2\epsilon)> \epsilon $ at each scale $j$ by (\ref{EQ130}).
\end{proof}
Next, we prove the following cluster coherence decays.
\begin{Prop} \label{PR2} Consider the shearlet frame $\boldsymbol{\Psi}$ of Definition  \ref{Df4} and the Gabor frame of scale $j, \mathbf{G_j}$,  defined in (\ref{DN1}).
We have 
$$\mu_c(\Lambda_{2,j},  \mathbf{G_j};\boldsymbol{\Psi}) \rightarrow 0, \quad j \rightarrow \infty,$$
where $\Lambda_{2,j}$ is defined in (\ref{CT7}).
\end{Prop}
\begin{proof}

Suppose that the maximum in the mutual coherence is attained at $\bar{l} \in \mathbb{Z}, \bar{k} \in \mathbb{Z}^2$ and $\rm{(\iota)} \in \{ \rm{(h)}, \rm{(v)}, \rm{(b)} \}$ (there is a maximum since the systems are translation invariant and the inner product goes to zero as $n$ goes to infinity). For each $N=1, 2 \dots $ using (\ref{EQ2}) and Lemma \ref{LM10} 3), we have
\begin{eqnarray*}  
\mu_c(\Lambda_{2,j},  \mathbf{G_j}; \boldsymbol{\Psi}) & = & \sum_{m \in B(0, M_j), n \in I_T^{\pm}}  | \langle (g_{s_j})_{m,n}, \psi_{j, \bar{l}, \bar{k}}^{\alpha_j, \iota} \rangle | \\
& \leq &  C_N \cdot 2^{-(2-\alpha_j) j/2 }  \# \{ m \in B(0,M_j) \} \cdot \#  \{ n \in I_T^{\pm} \cap M_{j, l, \rm{(\iota)}}\} \\
& \leq & C_N \cdot  2^{-(2-\alpha_j) j/2}  \cdot M_j^2 \cdot 2^{(2-\alpha_j-\epsilon)j/2}  \\
& \stackrel{\mathclap{\normalfont{M_j = 2^{\epsilon j/6} }}} =& \quad C_N \cdot 2^{-\epsilon j/6}
\xrightarrow{j \rightarrow +\infty} 0. 
\end{eqnarray*}
\end{proof}

\subsection{Cluster coherence of the projected frames}

\begin{figure}[H] 
  \centering
  \includegraphics[width=0.35\textwidth]{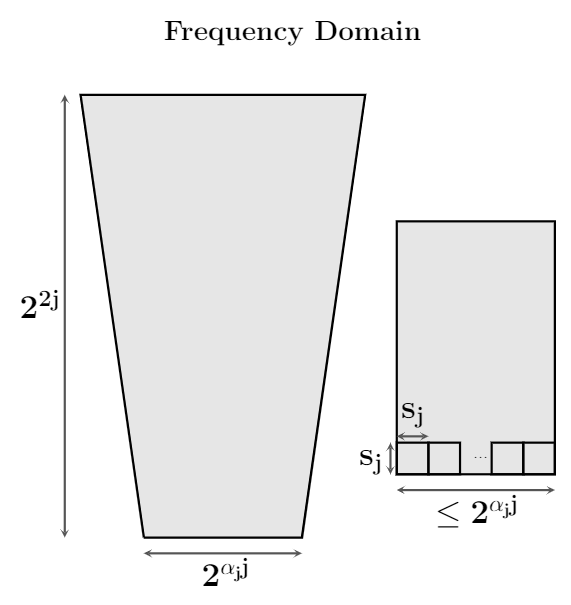}
  \quad
  \includegraphics[width=0.45\textwidth]{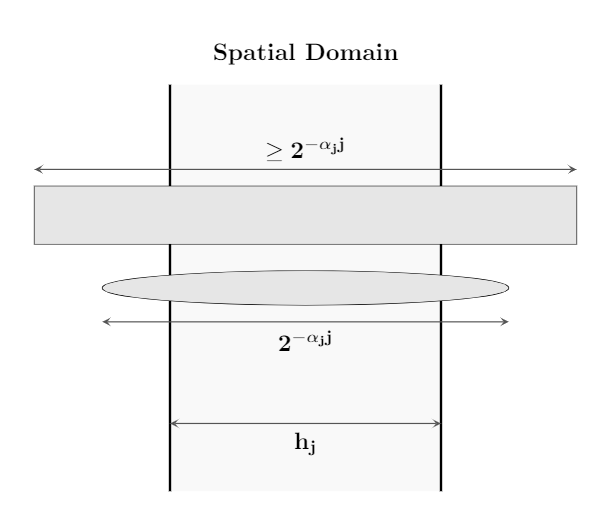}
  \caption{\label{FG10} Left: Frequency support of a shearlet (trapezoidal) and  Gabor elements (small squares).  Right: Missing part (vertical band), essential spatial support of a shearlet (ellipse) and a cluster of Gabor elements (rectangle).   The gap width does not exceed the essential length of the shearlet. }
\end{figure}

In this subsection we compute the cluster coherence terms corresponding to the missing stripe. To be able to inpaint the missing part, the following propositions rely of the fact that the gap size $h_j$ at each scale $j$ is  smaller than the essential length of the shearlet elements in this scale. For illustration, see Figure \ref{FG10}.

In \cite{25} the problem of inpainting a missing strip from a line singularity via shearlet frames was studied. As part of their construction, they prove the following proposition.

\begin{Prop} \textup{( \cite{25}, \textup{Prop.}  5.6)  }  \label{PR6} 
Consider the shearlet frame $\boldsymbol{\Psi}$ of Definition  \ref{Df4}. 
For $h_j= o(2^{(\alpha_j + \epsilon)j})$ with $\alpha_j \in (0,2),   \liminf \alpha_j >0 $  and $\epsilon$ satisfying (\ref{EQ130}), we have
$$\mu_c (\Lambda_{1,j}^{\pm}, P_j \boldsymbol{\Psi} ; \boldsymbol{\Psi}  ) \rightarrow 0, \quad j \rightarrow \infty, $$ 
where $\Lambda_{1,j}^{\pm}$ is defined in (\ref{CT6}), $P_j$ is defined in (\ref{DN2}).
\end{Prop}

Next, we prove the following result.

\begin{Prop} \label{PR7} Consider the the Gabor frame of scale $j, \mathbf{G_j}$  defined in (\ref{DN1}).
Suppose that $h_j = o(2^{-\alpha_j j })$  with $\alpha_j \in (0,2), \liminf \alpha_j >0,$ and $I_T$ satisfies $(\ref{EQ2})$ and $(\ref{EQ20})$.
Then we have
$$ \mu_c (\Lambda_{2,j}, P_j   \mathbf{G_j};   \mathbf{G_j} ) \rightarrow 0, \quad j \rightarrow \infty , $$
where $\Lambda_{2,j}$ is defined in (\ref{CT7}) and $P_j$ is defined in  (\ref{DN2}).
\end{Prop}

For an illustration of this Proposition, we refer to Figure \ref{FG10} (gray part). Next, we prove the following result.

\begin{proof}
First, from the definition of the cluster, we may assume that the maximum is attained for $(m^\prime, n^\prime) \in \Lambda_{2, j}$. Namely, we have
\begin{eqnarray}  \mu_c (\Lambda_{2,j}, P_j   \mathbf{G_j};   \mathbf{G_j} ) &=& \sum_{(m,n) \in M_j \times (I_T^{\pm} \cap \mathcal{A}_{s, j})}| \langle P_j (g_{s_j})_{m,n}, (g_{s_j})_{m^\prime,n^\prime}\rangle | \nonumber \\
& =& \sum_{(m,n) \in M_j \times (I_T^{\pm} \cap \mathcal{A}_{s,j})} \Big | \int_{-h_j}^{h_j} \int_\mathbb{R} (g_{s_j})_{m,n}(x) \overline{(g_{s_j})_{m^\prime,n^\prime}(x)} dx \Big | \nonumber \\
& \leq &  \sum_{(m,n) \in M_j \times (I_T^{\pm} \cap \mathcal{A}_{s, j})} \int_{-h_j}^{h_j} \int_\mathbb{R} |(g_{s_j})_{m,n}(x)| |(g_{s_j})_{m^\prime,n^\prime}(x)| dx. \qquad \label{A1}
\end{eqnarray}

By Lemma \ref{LM10}, we obtain the following decay estimate of $(g_{s_j})_{m,n}(x)$ for any $N\in\mathbb{N}$
$$ |(g_{s_j})_{m,n}(x)| \leq  C_N \cdot s_j \cdot  \langle |s_j x_1 + \frac{m_1}{2}|\rangle^{-N}\langle |s_j x_2 + \frac{m_2}{2}|\rangle^{-N} .$$
Thus, by changing variable $y = s_jx$, we have
\begin{eqnarray}
 &&  \int_{-h_j}^{h_j} \int_\mathbb{R}  |(g_{s_j})_{m,n}(x)|  |(g_{s_j})_{m^\prime,n^\prime}(x)| dx \nonumber \\
 &\leq&  \int_{-h_j}^{h_j} \int_\mathbb{R}  C_N s_j^2  \langle|s_j x_1 + \frac{m_1}{2} |\rangle^{-N} \langle |s_j x_2 + \frac{m_2}{2} |\rangle^{-N}  \langle|s_j x_1 + \frac{m_1^\prime}{2} |\rangle^{-N} 
\nonumber \\
&& \qquad \quad  \cdot \langle |s_j x_2 + \frac{m_2^\prime}{2} |\rangle^{-N} dx  \nonumber \\
 & = & \int_{-h_js_j}^{h_js_j} \int_\mathbb{R}  C_N \langle |y_1 + \frac{m_1}{2} |\rangle^{-N} \langle | y_2 + \frac{m_2}{2} |\rangle^{-N}  \langle |y_1 + \frac{m_1^\prime}{2} |\rangle^{-N} \langle |y_2 + \frac{m_2^\prime}{2} |\rangle^{-N} dy. \nonumber \\ 
 \label{A2}
 \end{eqnarray}
 
 Next, 
\begin{equation}
    \sum_m \langle |y_1 + \frac{m_1}{2} |\rangle^{-N} \langle | y_2 + \frac{m_2}{2} |\rangle^{-N} \leq C_N^\prime, \, \textup{ and } \int_\mathbb{R} \langle |y_2 + \frac{m_2^\prime}{2} |\rangle^{-N} dy_2 \leq C_N^{\prime \prime}. \qquad \label{A3}
\end{equation} 
 Thus, by (\ref{EQ20}) and (\ref{A1}) combined with (\ref{A2}) and (\ref{A3}) , we obtain
 \begin{eqnarray*}
  \mu_c (\Lambda_{2,j}, P_j \mathbf{G_j}; \mathbf{G_j} ) &\leq & \sum_{n \in I_T^{\pm} \cap \mathcal{A}_{s,j}} C_N \cdot h_js_j   \\
  & \leq & C_N \cdot \frac{2^{\alpha_j j}}{s_j}h_js_j = C_N \cdot 2^{\alpha_j j }h_j\xrightarrow  {j \rightarrow +\infty} 0. \\
 \end{eqnarray*}
\end{proof}

\begin{Prop} \label{PR9} 
Consider the shearlet frame $\boldsymbol{\Psi}$ of Definition  \ref{Df4} and the Gabor frame of scale $j, \mathbf{G_j}$,  defined in (\ref{DN1}).
Assuming that $I_T$ satisfies $(\ref{EQ2})$ and $(\ref{EQ20})$, we have
$$ \mu_c (\Lambda_{2,j}, P_j   \mathbf{G_j}, \boldsymbol{\Psi}) \rightarrow 0, \quad j \rightarrow \infty ,$$
where $\Lambda_{2,j}$ is defined in (\ref{CT7}) and $P_j$ is defined in (\ref{DN2}) .
\end{Prop}
\begin{proof}
The maximum of the cluster coherence is attained for some $(j,\bar{l},\bar{k}; \alpha_j,  \rm{(\iota)}) \in \Delta_j$. Thus,
\begin{eqnarray*} \mu_c (\Lambda_{2,j}, P_j   \mathbf{G_j}, \boldsymbol{\Psi}) &=& \sum_{(m,n) \in M_j \times (I_T^{\pm} \cap \mathcal{A}_{s, j})}| \langle P_j (g_{s_j})_{m,n},  \psi_{j, \bar{l}, \bar{k}}^{\alpha_j, \rm{(\iota)}}\rangle |\\
& =& \sum_{(m,n) \in M_j \times (I_T^{\pm} \cap \mathcal{A}_{s, j})} \Big | \int_{-h_j}^{h_j} \int_\mathbb{R} (g_{s_j})_{m,n}(x) \overline{\psi_{j, \bar{l}, \bar{k}}^{\alpha_j, \rm{(\iota)}}(x)} dx \Big | \\
& \leq &  \sum_{(m,n) \in M_j \times (I_T^{\pm} \cap \mathcal{A}_{s, j})} \int_{-h_j}^{h_j} \int_\mathbb{R} |(g_{s_j})_{m,n}(x)| |\psi_{j, \bar{l}, \bar{k}}^{\alpha_j, \rm{(\iota)}}(x)| dx. 
\end{eqnarray*}

Now we consider three cases 

1) Case $\rm{(\iota)} = \rm{(v)}.$ For each $N=1,2, \dots$, by Lemma \ref{LM10}, we derive the following decay estimates of $(g_{s_j})_{m,n}(x)$ and $\psi_{j, \bar{l}, \bar{k}}^{\alpha_j, \rm{(\iota)}}(x)$
$$ |(g_{s_j})_{m,n}(x)| \leq  C_N \cdot s_j \cdot  \langle |s_j x_1 + \frac{m_1}{2}|\rangle^{-N}\langle |s_j x_2 + \frac{m_2}{2}|\rangle^{-N} ,$$
$$ | \psi_{j, \bar{l}, \bar{k}}^{\alpha_j, \rm{(v)}}(x)| \leq C_N \cdot 2^{(2+\alpha)j/2} \cdot \langle |2^{\alpha_j j}x_1 - \bar{k}_1|\rangle^{-N} \langle |2^{2j}x_2+ \bar{l}2^{\alpha_j j } x_1-\bar{k}_2|\rangle^{-N}.$$
Thus, by the  change of variable $y =S^{\bar{l}}_{\rm{(v)}}A_{\alpha_j, \rm{(v)}}^{j} x= (2^{\alpha_j j}x_1, 2^{2j}x_2 + \bar{l}2^{\alpha_j j}x_1)$, we have
\begin{eqnarray}
 &&  \int_{-h_j}^{h_j} \int_\mathbb{R}  |(g_{s_j})_{m,n}(x)| | \psi_{j, \bar{l}, \bar{k}}^{\alpha_j, \rm{(v)}}(x)|   dx \nonumber \\
 &\leq &  C_N\cdot 2^{-(2+\alpha_j)j/2}\cdot s_j   \int_{-2^{\alpha_j j}h_j}^{2^{\alpha_j j}h_j} \int_\mathbb{R}     \langle|s_j 2^{-\alpha_j j}y_1 + \frac{m_1}{2} |\rangle^{-N} \nonumber \\
 && \qquad \cdot \langle |s_j 2^{-2j}(y_2-\bar{l}y_1) + \frac{m_2}{2}|\rangle^{-N} \langle |y_1-\bar{k}_1|\rangle^{-N}  \langle |y_2- \bar{k}_2|\rangle^{-N} dy_1dy_2 \nonumber \\
  &\leq &  C_N\cdot 2^{-(2+\alpha_j)j/2}\cdot s_j  \int_{\mathbb{R}} \int_\mathbb{R}     \langle|s_j 2^{-\alpha_j j}y_1 + \frac{m_1}{2} |\rangle^{-N} \langle |s_j 2^{-2j}(y_2-\bar{l}y_1) + \frac{m_2}{2}|\rangle^{-N}  \nonumber \\
  && \qquad \cdot \langle |y_1-\bar{k}_1|\rangle^{-N}  \langle |y_2- \bar{k}_2|\rangle^{-N} dy_1dy_2. \label{B1}
 \end{eqnarray}
 Furthermore, 
\begin{equation} \label{B2}
     \sum_m \langle |s_j 2^{-\alpha_j j}y_1 + \frac{m_1}{2} |\rangle^{-N} \langle |s_j 2^{-2j}(y_2-\bar{l}y_1)  + \frac{m_2}{2} |\rangle^{-N} \leq C_N^\prime 
\end{equation}
 and
\begin{equation} \label{B3}
    \int_\mathbb{R} \langle |y_1 - \bar{k}_1 |\rangle^{-N} dy_1 \leq C_N^{\prime \prime}, \quad \int_\mathbb{R} \langle |y_2 - \bar{k}_2 |\rangle^{-N} dy_2 \leq C_N^{\prime \prime \prime}.
\end{equation}
Thus, by (\ref{B1}), (\ref{B2}), (\ref{B3})  and (\ref{EQ20}), we obtain
\begin{eqnarray*}
 \mu_c (\Lambda_{2,j}, P_j   \mathbf{G_j}; \boldsymbol{\Psi})  &\leq & \sum_{n \in I_T^{\pm} \cap \mathcal{A}_{s, j}} C_N \cdot  2^{-(2+\alpha_j )j/2} \cdot s_j  \\
& \leq & C_N \cdot  2^{-(2+\alpha_j )j/2} \cdot s_j \cdot \frac{2^{\alpha_j j}}{s_j}   =   C_N \cdot  2^{-(2-\alpha_j )j/2} \xrightarrow  {j \rightarrow +\infty} 0. 
\end{eqnarray*}  

2) Case $\rm{(\iota)} = \rm{(h)} $. Similarly, by the change of variable $y = S^{\bar{l}}_{\rm{(h)}}A_{\alpha_j, \rm{(h)}}^{j} x= (2^{2 j}x_1+ \bar{l}2^{\alpha_j j}x_2, 2^{\alpha_j j}x_2)$, we have
\begin{eqnarray*}
 && 2^{(2+\alpha_j)j/2}\cdot s_j^{-1} \cdot \int_{-h_j}^{h_j} \int_\mathbb{R}  |(g_{s_j})_{m,n}(x)| | \psi_{j, \bar{l}, \bar{k}}^{\alpha_j, \rm{(h)}}(x)|dx \\
   &\leq&  C_N  \int_{\mathbb{R}^2}     \langle|s_j 2^{-2 j}(y_1-\bar{l}y_2) + \frac{m_1}{2} |\rangle^{-N} \langle |s_j 2^{-\alpha_j j}y_2 + \frac{m_2}{2}|\rangle^{-N} \\
   && \qquad  \cdot  \langle |y_1-\bar{k}_1|\rangle^{-N} 
  \langle |y_2- \bar{k}_2|\rangle^{-N} dy.
\end{eqnarray*}  

By using similar argument as in case $\rm{(\iota)}=\rm{(v)}$, we finally obtain
\begin{eqnarray*}
  \mu_c (\Lambda_{2,j}, P_j \mathbf{G_j}; \boldsymbol{\Psi}) & \leq &  \sum_{n \in I_T^{\pm} \cap \mathcal{A}_{s, j}} C_N \cdot  2^{-(2+\alpha_j )j/2} \cdot s_j  \\
& \stackrel{\mathclap{\normalfont{(\ref{EQ20}) }}} \leq& \; C_N \cdot2^{-(2- \alpha_j)j/2}   \xrightarrow  {j \rightarrow +\infty} 0.
\end{eqnarray*}
1) Case $\rm{(\iota)} = \rm{(b)}.$ Recalling the definition of boundary shearlets, the decay estimate can be done similarly.
\end{proof}

\begin{Prop} \label{PR10}
Consider the shearlet frame $\boldsymbol{\Psi}$ of Definition  \ref{Df4} and the Gabor frame of scale $j, \mathbf{G_j}$,  defined in (\ref{DN1}).
Suppose that $s_j \leq 2^{\alpha_j j}, \,  h_j = o(2^{\alpha_j  j }),$ and $\alpha_j \in (0, 2), \liminf \alpha_j >0$. Then we have
$$  \mu_c (\Lambda_{1,j}^{\pm}, P_j \boldsymbol{\Psi};  \mathbf{G_j}) \rightarrow 0, \quad j \rightarrow \infty ,$$
where $\Lambda_{1,j}^{\pm}$ is defined in (\ref{CT6}) and $P_j$ is defined in (\ref{DN2}).
\end{Prop}
\begin{proof}
Without loss of generality,  we prove only $\mu_c (\Lambda_{1,j}, P_j \boldsymbol{\Psi}   ;  \mathbf{G_j}) \rightarrow 0, \; j \rightarrow \infty$. First, we assume that the maximum in the definition of the cluster is attained for some $(\bar{m}, \bar{n}) \in \mathbb{Z}^2 \times \mathbb{Z}^2  $. We have
\begin{eqnarray*} \mu_c (\Lambda_{1,j}, P_j \boldsymbol{\Psi}; \mathbf{G_j} ) &=&  \sum_{(j, l, k; \alpha_j,  \iota) \in \Lambda_{1,j}} | \langle P_j \psi_{j, l, k}^{\alpha_j, \rm{(\iota)}} , (g_{s_j})_{\bar{m},\bar{n}}\rangle |\\
& =& \sum_{\substack{|l| \leq 1, k \in \mathbb{Z}^2 \\ |k_2 - lk_1| \leq 2^{\epsilon j} }}  \Big | \int_{-h_j}^{h_j} \int_\mathbb{R} \psi_{j, l, k}^{\alpha_j, \rm{(v)}}(x) \overline{(g_{s_j})_{\bar{m},\bar{n}}(x)}  dx \Big | \\
& \leq &  \sum_{\substack{|l| \leq 1, k \in \mathbb{Z}^2 \\ |k_2 - lk_1| \leq 2^{\epsilon j} }}   \int_{-h_j}^{h_j} \int_\mathbb{R} |\psi_{j, l, k}^{\alpha_j, \rm{(v)}}(x)| |(g_{s_j})_{\bar{m},\bar{n}}(x)|  dx|. 
\end{eqnarray*}
For each $N=1,2, \dots$, by Lemma \ref{LM10}, we derive the following decay estimates of $\psi_{j, l, k}^{\alpha_j, \rm{(v)}}(x)$  and  $(g_{s_j})_{\bar{m},\bar{n}}(x)$
$$ | \psi_{j, l, k}^{\alpha_j, \rm{(v)}}(x)| \leq C_N \cdot 2^{(2+\alpha)j/2} \cdot \langle |2^{\alpha_j j}x_1 - k_1|\rangle^{-N} \langle |2^{2j}x_2+ l2^{\alpha_j j } x_1-k_2|\rangle^{-N},$$
$$ |(g_{s_j})_{\bar{m},\bar{n}}(x)| \leq  C_N \cdot s_j \cdot  \langle |s_j x_1 + \frac{\bar{m}_1}{2}|\rangle^{-N}\langle |s_j x_2 + \frac{\bar{m}_2}{2}|\rangle^{-N} .$$
Thus, by changing variable $y =S^{l}_{\rm{(v)}}A_{\alpha_j, \rm{(v)}}^{j} x= (2^{\alpha_j j}x_1, 2^{2j}x_2 + l2^{\alpha_j j}x_1)$, we have
\begin{eqnarray}
 &&  \int_{-h_j}^{h_j} \int_\mathbb{R}  | \psi_{j, l, k}^{\alpha_j, \rm{(v)}}(x)| \cdot |(g_{s_j})_{\bar{m},\bar{n}}(x)|   dx \nonumber \\
 &\leq &  C_N \cdot 2^{-(2+\alpha_j)j/2}\cdot s_j \cdot \int_{-2^{\alpha_j j}h_j}^{2^{\alpha_j j}h_j} \int_\mathbb{R}   \langle |y_1-k_1|\rangle^{-N} \langle |y_2- k_2|\rangle^{-N}    \nonumber \\
 && \qquad \cdot \langle|s_j 2^{-\alpha_j j}y_1 + \frac{\bar{m}_1}{2} |\rangle^{-N}
 \langle |s_j 2^{-2j}(y_2-ly_1) + \frac{\bar{m}_2}{2}|\rangle^{-N}  dy_1dy_2. \label{C1}
 \end{eqnarray}

Next, 
\begin{equation} \label{C2}
    \langle |s_j 2^{-\alpha_j j}y_1 + \frac{\bar{m}_1}{2} |\rangle^{-N} \langle |s_j 2^{-2j}(y_2-ly_1)  + \frac{\bar{m}_2}{2} |\rangle^{-N} \leq 1
\end{equation} 
 and
 \begin{eqnarray}
&& \int_\mathbb{R} \Big ( \sum_{\substack{|l| \leq 1, k \in \mathbb{Z}^2 \\ |k_2 - lk_1| \leq 2^{\epsilon j} }}\langle |y_1-k_1|\rangle^{-N}  \langle |y_2- k_2|\rangle^{-N} \Big ) dy_2  \nonumber \\
& = &  \sum_{\substack{|l| \leq 1, k_1, k_3 \in \mathbb{Z} \\ |k_3| \leq 2^{\epsilon j} }}\langle |y_1-k_1|\rangle^{-N} \cdot \Big ( \int_\mathbb{R} \langle |y_2- k_3-lk_1|\rangle^{-N} dy_2 \Big )\nonumber \\
 &\leq & C_N \cdot 2^{\epsilon j} \cdot \int_{\mathbb{R}} \langle |t-y_1|\rangle^{-N}  dt \leq C_N^\prime 2^{\epsilon j}. \label{C3}
 \end{eqnarray}

Thus, by (\ref{C1}), (\ref{C2}),  (\ref{C3}), we obtain
\begin{eqnarray*}
 \mu_c (\Lambda_{1,j}, P_j \boldsymbol{\Psi}; \mathbf{G_j} ) &\leq& C_N^\prime \cdot 2^{-(2+ \alpha_j)j/2} \cdot s_j \cdot 2^{\alpha_j j} h_j \cdot 2^{\epsilon j} \\
& \stackrel{\mathclap{\normalfont{s_j \leq  2^{\alpha_j j} }}} \leq& \;\; C_N^\prime \cdot 2^{-(2-\alpha_j-2\epsilon)j/2} \cdot 2^{\alpha_j j} h_j \xrightarrow{j \rightarrow +\infty} 0 .
\end{eqnarray*}

\end{proof}

\section{Proof of theorem \ref{TR10} } \label{S10}

We are now ready to present the proof of Theorem \ref{TR10}.
\begin{proof}
By Theorem \ref{TR1}, we have
\begin{equation}
    \| P_j\mathcal{C}_j^\star - P_j w \mathcal{S}_j  \|_2+ \| P_j \mathcal{T}_j^\star - P_j \mathcal{T}_{s,j}  \|_2  \leq \frac{2\delta_j}{1 - 2 \mu_{c,j}}, \label{Z1}
\end{equation} 
where $\delta_j= \delta_{1,j} + \delta_{2,j}$ and 
    $\mu_{c,j} = $ \\ $\max  \{ \mu_c(\Lambda_{1,j}^\pm, P_j \boldsymbol{\Psi};\boldsymbol{\Psi}) +  \mu_c(\Lambda_{2,j}, P_j \mathbf{G_j}; \boldsymbol{\Psi}), 
 \mu_c(\Lambda_{2,j}, P_j \mathbf{G_j}; \mathbf{G_j}) 
 +  \mu_c(\Lambda_{1,j}^\pm, P_j \boldsymbol{\Psi}; \mathbf{G_j})  \}  \\
+ \max  \{ \mu_c(\Lambda_{1,j}^\pm, \Phi_1; \mathbf{G_j}), \mu_c(\Lambda_2, \mathbf{G_j}; \boldsymbol{\Psi})  \}.$ \\
Moreover, it follows from Propositions \ref{PR4} and \ref{PR1} that
\begin{equation}
    \delta_j = \delta_1 + \delta_2 = o(2^{-Nj}), \; \forall N \in \mathbb{N}. \label{Z2}
\end{equation}
For the other term, the following estimate holds as a consequence of 
 Propositions   \ref{PR3}, \ref{PR2}, \ref{PR6}, \ref{PR7}, \ref{PR9} and  \ref{PR10}
 \begin{equation}
     \mu_{c,j} \longrightarrow 0, \quad j \rightarrow +\infty. \label{Z3}
 \end{equation}
Combining (\ref{Z1}),  (\ref{Z2}), and  (\ref{Z3}), we get
$$  \| \mathcal{C}_j^\star - w \mathcal{S}_j  \|_2+ \|  \mathcal{T}_j^\star - \mathcal{T}_{s,j}  \|_2  = o(2^{-Nj}), \quad \forall N \in \mathbb{N}.$$ This proves the first claim of the theorem.

The second claim is obtained since the following estimate holds for any arbitrarily large number $N \in \mathbb{N}$ 
\begin{equation}
    \| P_j\mathcal{C}_j^\star - P_j w \mathcal{S}_j  \|_2+ \| P_j \mathcal{T}_j^\star - P_j \mathcal{T}_{s,j}  \|_2 \leq \| \mathcal{C}_j^\star - \mathcal{C}_j  \|_2+ \| \mathcal{T}_j^\star - \mathcal{T}_{s,j}  \|_2 \leq 2^{-Nj}. \label{Z4}
\end{equation} 
By Lemma \ref{LM52} and (\ref{Z4}), for $N \in \mathbb{N}$ chosen sufficiently large, we obtain
$$ \frac{\| P_j\mathcal{C}_j^\star - P_j w \mathcal{S}_j  \|_2}{\|P_j w \mathcal{S}_j  \|_2} \leq \frac{o(2^{-Nj})}{h_j 2^{-j}} \xrightarrow{j \rightarrow +\infty}{0}, $$
and
$$ \frac{\| P_j \mathcal{T}_j^\star - P_j \mathcal{T}_{s,j}  \|_2}{h_j2^{-j}} =\frac{o(2^{-Nj})}{h_j 2^{-j}} \xrightarrow{j \rightarrow +\infty}{0},$$
which concludes the proof.
\end{proof}

\textcolor{black}{\section*{Acknowledgements}
VT.D. acknowledges support by the VIED-MOET Fellowship through Project 911. 
R.L. acknowledges support by the DFG SPP 1798 “Compressed Sensing in Information Processing” through Project Massive MIMO-II. G.K. acknowledges support by the Deutsche Forschungsgemeinschaft (DFG) through Project KU 1446/21-2 within SPP 1798.}

\end{document}